\documentclass[twoside,leqno,a4paper]{amsart}
\usepackage{amssymb,amsmath}

\newcommand{\CC}{{\mathbb{C}}}
\newcommand{\FF}{{\mathbb{F}}}
\newcommand{\GG}{{\mathbb{G}}}
\newcommand{\QQ}{{\mathbb{Q}}}
\newcommand{\RR}{{\mathbb{R}}}
\newcommand{\TT}{{\mathbb{T}}}
\newcommand{\ZZ}{{\mathbb{Z}}}
\newcommand{\NN}{{\mathbb{N}}}

\newcommand{\cA}{{\mathcal{A}}}

\newcommand{\cM}{{\mathcal{M}}}
\newcommand{\cR}{{\mathcal{R}}}
\newcommand{\cS}{{\mathcal{S}}}

\newcommand{\bS}{{\mathbf{S}}}
\newcommand{\bb}{{\mathbf{b}}}
\newcommand{\bs}{{\mathbf{s}}}
\newcommand{\bw}{{\mathbf{w}}}

\newcommand{\fA}{{\mathfrak{A}}}
\newcommand{\fh}{{\mathfrak{h}}}
\newcommand{\fg}{{\mathfrak{g}}}
\newcommand{\fS}{{\mathfrak{S}}}

\newcommand{\bbeta}{{\boldsymbol{\beta}}}
\newcommand{\bpi}{{\boldsymbol{\pi}}}

\newcommand{\GL}{{\operatorname{GL}}}
\newcommand{\SL}{{\operatorname{SL}}}
\newcommand{\OO}{{\operatorname{O}}}
\newcommand{\UU}{{\operatorname{U}}}
\newcommand{\id}{{\operatorname{id}}}
\newcommand{\tr}{{\operatorname{tr}}}
\newcommand{\regV}{{V^{\operatorname{reg}}}}
\newcommand{\Ind}{{\operatorname{Ind}}}
\newcommand{\Irr}{{\operatorname{Irr}}}

\newcommand{\eps}{\epsilon}

\renewcommand{\leq}{\leqslant}
\renewcommand{\geq}{\geqslant}
\renewcommand{\atop}[2]{\genfrac{}{}{0pt}{}{#1}{#2}}
\newcommand{\braid}[1]{\boldsymbol{#1}}

\newtheoremstyle{theorem}%
 {5pt plus 1pt minus 1pt}{5pt plus 1pt minus 1pt}{\itshape}%
 {0pt}{\bfseries\boldmath}{.}{ }%
 {\thmnumber{#2}\thmname{.~#1}\thmnote{ {\mdseries\upshape(#3)}}}
\newtheoremstyle{definition}%
 {5pt plus 1pt minus 1pt}{5pt plus 1pt minus 1pt}{\upshape}%
 {0pt}{\bfseries\boldmath\itshape}{.}{ }%
 {\thmnumber{#2}\thmname{.~#1}\thmnote{ {\mdseries\upshape(#3)}}}

\theoremstyle{theorem}
\newtheorem{thm}{Theorem}[section]

\newtheorem{cor}[thm]{Corollary}
\newtheorem{prop}[thm]{Proposition}

\theoremstyle{definition}

\makeatletter
\newenvironment{para}[1]{\refstepcounter{thm}
\smallskip \noindent {\bfseries\boldmath\thethm.~{#1}.}}{\par\smallskip}
\makeatother

\begin{document}

\title{Contribution to the Handbook of Algebra}

\author{Meinolf Geck}
\address{M.G.: Institut Girard Desargues, bat. Jean Braconnier, Universit\'e 
Lyon 1, 21 av Claude Bernard, F--69622 Villeurbanne cedex, France}
\makeatletter
\email{geck@desargues.univ-lyon1.fr}
\makeatother

\author{Gunter Malle}
\address{G.M.: FB Mathematik/Informatik, Universit\"at Kassel, 
Heinrich-Plett-Str.~40, D--34132 Kassel, Germany.}
\makeatother
\email{malle@mathematik.uni-kassel.de}
\makeatother

\date{}

\maketitle

\pagestyle{myheadings}

\bigskip
\begin{center} {\bf\Large Chapter: Reflection groups} \end{center}
\bigskip

This chapter is concerned with the theory of finite reflection groups, that 
is, finite groups generated by reflections in a real or complex vector space.
This is a rich theory, both for intrinsic reasons and as far as applications
in other mathematical areas or mathematical physics are concerned. The origin
of the theory can be traced back to the ancient study of symmetries of 
regular polyhedra. Another extremely important impetus comes from the theory 
of semisimple Lie algebras and Lie groups, where finite reflection groups
occur as ``Weyl groups''. In the last decade, Brou\'e's ``abelian defect
group conjecture'' (a conjecture concerning the representations of finite
groups over fields of positive characteristic) has lead to a vast research
program, in which complex reflection groups, corresponding braid groups and
Hecke algebras play a prominent role. Thus, the theory of reflection groups
is at the same time a well-established classical piece of mathematics and
still a very active research area.  The aim of this chapter (and a subsequent
one on Hecke algebras) is to give an overwiew of both these aspects. 

As far as the study of reflection groups as such is concerned, there are
(at least) three reasons why this leads to an interesting and  rich theory:

\medskip
{\bf Classification:} Given a suitable notion of ``irreducible'' reflection 
groups, it is possible to give a complete classification, with typically 
several infinite families of groups and a certain number of exceptional cases.
In fact, this classification can be seen as the simplest possible model for
much more complex classification results concerning related algebraic 
structures, such as complex semisimple Lie algebras, simple algebraic groups 
and, eventually, finite simple groups. Besides the independent interest of
such a classification, we mention that there is a certain number of results 
on finite reflection groups which can be stated in general terms but whose 
proof requires a case--by--case analysis according to the classification.
(For example, the fact that every element in a finite real reflection 
group is conjugate to its inverse.)

\medskip
{\bf Presentations:} Reflection groups have a highly symmetric ``Coxeter 
type presentation'' with generators and defining relations (visualized 
by ``Dynkin diagrams'' or generalisations thereof), which makes it possible
to study them by purely combinatorial methods (length function, reduced 
expressions and so on). From this point of view, the associated Hecke 
algebras can be seen as ``deformations'' of the group algebras of 
finite reflection groups, where one or several formal parameters are
introduced into the set of defining relations. One of the most important
developments in this direction is the discovery of the Kazhdan--Lusztig 
polynomials and the whole theory coming with them. (This is discussed
in more detail in the chapter on Hecke algebras.)

\medskip
{\bf Topology and geometry:} The action of a reflection group on the 
underlying vector space opens the possibility of using geometric methods. 
First of all, the ring of invariant symmetric functions on that vector space
always is a polynomial ring (and this characterizes finite reflection 
groups). Furthermore, we have a corresponding hyperplane arrangement which
gives rise to the definition of an associated braid group as the fundamental 
group of a certain topological space. For the symmetric group, we obtain 
in this way the classical Artin braid group, with applications in the 
theory of knots and links. 

\medskip
Furthermore, all these aspects are related to each other which---despite
being quite elementary taken individually---eventually leads to a highly
sophisticated theory.

We have divided our survey into four major parts.  The first part deals 
with finite complex or real reflection groups in general.  The second part 
deals with finite real reflection groups and the relations with the theory 
of Coxeter groups. The third part is concerned with the associated braid 
groups. Finally, in the fourth part, we consider complex irreducible 
characters  of finite reflection groups.

We certainly do not pretend to give a complete picture of all aspects of 
the theory of reflection groups and Coxeter groups. Especially, we will
not say so much about areas that we do not feel competent in; at our best 
knowledge, we try to give at least some references for further reading in 
such cases. This concerns, in particular, all aspects of infinite (affine, 
hyperbolic $\ldots$) Coxeter groups.

\section{Finite groups generated by reflections} \label{sec:1}

\begin{para}{Definitions} \label{subsec:1.Def} 
Let $V$ be a finite dimensional vector space over a field $K$. A
{\em reflection on $V$} is a non-trivial element $g\in\GL(V)$ of finite order
which fixes a hyperplane in $V$ pointwise. There are two types of reflections,
according to whether $g$ is semisimple (hence diagonalisable) or unipotent.
Often, the term reflection is reserved for the first type of elements, while
the second are called {\em transvections}. They can only occur in positive
characteristic. Here, we will almost exclusively be concerned with ground
fields $K$ of characteristic~0, which we may and will then assume to be
subfields of the field $\CC$ of complex numbers. Then, by our
definition, reflections are always semisimple and (thus) diagonalisable.
Over fields $K$ contained in the field $\RR$ of real numbers, reflections
necessarily have order~2, which is the case motivating their name. Some authors
reserve the term reflection for this case, and speak of pseudo-reflections in
the case of arbitrary (finite) order.  \par
Let $g\in\GL(V)$ be a reflection. The hyperplane $C_V(g)$ fixed point-wise by
$g$ is called the {\em reflecting hyperplane} of $g$. Then $V=C_V(g)\oplus V_g$
for a unique
$g$-invariant subspace $V_g$ of $V$ of dimension~1. Any non-zero vector
$v\in V_g$ is called a {\em root for $g$}. Thus, a root for a reflection is an
eigenvector with eigenvalue different from~1. Now assume in addition that
$V$ is Hermitean. Then conversely, given a vector $v\ne0$ in $V$ and a natural
number $n\ge2$ we may define a reflection in $V$ with root $v$ and of
order $n$ by $g.v:=\exp(2\pi i/n)v$, and $g|_{V^\perp}=\id$.
 \par
A {\em reflection group on $V$} is now a finite subgroup $W\leq\GL(V)$
generated by reflections. Note that any finite subgroup of $\GL(V)$ leaves
invariant a non-degenerate Hermitean form. Thus, there is no loss in assuming
that a reflection group $W$ leaves such a form invariant.
\end{para}

\begin{para}{Invariants} \label{subsec:1.Inv} 
Let $V$ be a finite-dimensional vector space over $K\leq\CC$. Let $K[V]$
denote the algebra of symmetric functions on $V$, i.e., the symmetric algebra
$S(V^*)$ of the dual space $V^*$ of $V$. So $K[V]$ is a commutative algebra
over $K$ with a grading $K[V]=\oplus_{d\ge0} K[V]^d$, where, for any $d\ge0$,
$K[V]^d$ denotes the $d$th symmetric power of $V^*$. If $W\leq\GL(V)$ then
$W$ acts naturally on $K[V]$, respecting the grading. Now reflection groups
are characterised by the structure of their invariant ring $K[V]^W$:
\end{para}

\begin{thm}[\rm Shephard--Todd \cite{ShTo}, Chevalley \cite{Ch}]
 \label{Gthm:shtochev}
 Let $V$ be a finite-dimen\-sional vector space over a field of
 characteristic~0 and $W\leq\GL(V)$ a finite group. Then the following are
 equivalent:
 \begin{itemize}
 \item[(i)] the ring of invariants $K[V]^W$ is a polynomial ring,
 \item[(ii)] $W$ is generated by reflections.
 \end{itemize}
\end{thm}

The implication from (i) to (ii) is an easy consequence of Molien's formula
$$P(K[V]^W,x)=\frac{1}{|W|}\sum_{g\in W}\frac{1}{\det_V(1-gx)}$$
for the Hilbert series $P(K[V]^W,x)$ of the ring of invariants $K[V]^W$, see
Shephard--Todd \cite[p.289]{ShTo}. It follows from Auslander's purity of the 
branch locus that this implication remains true in arbitrary characteristic 
(see Benson \cite[Th.~7.2.1]{Ben}, for example). The other direction was 
proved by Shephard and Todd as an application of their classification of 
complex reflection groups (see Section~(\ref{subsec:1.Class})). Chevalley
gave a general proof avoiding the classification which uses the combinatorics
of differential operators.

Let $W$ be an $n$-dimensional reflection group. By Theorem~\ref{Gthm:shtochev}
the ring of invariants is generated by $n$ algebraically independent
polynomials (so-called {\em basic invariants}), which may be taken to be
homogeneous. Although these polynomials are not uniquely determined in
general, their degrees $d_1\leq\ldots\leq d_n$ are. They
are called the {\em degrees of $W$}. Then $|W|=d_1\cdots d_n$, and the Molien
formula shows that $N(W):=\sum_{i=1}^n m_i$ is the number of reflections in $W$,
where $m_i:=d_i-1$ are the {\em exponents of $W$}. \par
The quotient $K[V]_W$ of $K[V]$ by the ideal generated by the invariants of
strictly positive degree is called the {\em coinvariant algebra} of $(V,W)$.
This is again a naturally graded $W$-module, whose structure is described by:

\begin{thm} [\rm Chevalley \cite{Ch}]
 \label{Gthm:chevalley}
 Let $V$ be a finite-dimen\-sional vector space over a field $K$ of
 characteristic~0 and $W\leq\GL(V)$ a reflection group. Then $K[V]_W$ carries a
 graded version of the regular representation of $W$. The grading is such that
 $$\sum_{i\ge0}\dim K[V]_W^i\ x^i=\prod_{i=1}^n\frac{x^{d_i}-1}{x-1},$$
 where $K[V]_W^i$ denotes the homogeneous component of degree~$i$.
\end{thm}

(See also Bourbaki \cite[V.5.2, Th.~2]{Bou}.) The polynomial
\[ P_W:=\sum_{i\ge0}\dim K[V]_W^i x^i\]
is called the {\em Poincar\'e-polynomial of $W$}.

\begin{para}{Parabolic subgroups} \label{subsec:1.parabolic}
Let $W$ be a reflection group on $V$. The {\em parabolic subgroups of $W$} are
by definition the pointwise stabilisers
$$W_{V'}:=\{g\in W\mid g\cdot v=v\text{ for all }v\in V'\}.$$
of subspaces $V'\leq V$. The following result is of big importance in the
theory of reflection groups:
\end{para}

\begin{thm} [\rm Steinberg \cite{Ste}]
 \label{Gthm:steinberg}
 Let $W\leq\GL(V)$ be a complex reflection group. For any subspace $V'\leq V$
 the parabolic subgroup $W_{V'}$ is generated by the reflections it contains,
 that is, by the reflections whose reflecting hyperplane contains $V'$.
 In particular, parabolic subgroups are themselves reflection groups.
\end{thm}

For the proof, Steinberg characterises reflection groups via eigenfunctions
of differential operators with constant coefficients that are invariant under
finite linear groups. For a generalisation to positive characteristic see
Theorem~\ref{Gthm:serre}.

\begin{para}{Exponents, coexponents and fake degrees} \label{subsec:1.fake}
Let $W$ be a complex reflection group. For $w\in W$ define
$k(w):=\dim V^{\langle w\rangle}$, the dimension of the fixed space of $w$
on $V$. Solomon \cite{So63} proved the following remarkable formula for the
generating function of $k$
$$\sum_{w\in W} x^{k(w)} = \prod_{i=1}^n (x+m_i),$$
by showing that the algebra of $W$-invariant differential forms with
polynomial coefficients is an exterior algebra of rank $n$ over the algebra
$K[V]^W$, generated by the differentials of a set of basic invariants (see
also Flatto \cite{Fla78}, Benson \cite[Th.~7.3.1]{Ben}). The formula was
first observed by Shephard--Todd \cite[5.3]{ShTo} using their
classification of irreducible complex reflection
groups. Dually, Orlik--Solomon \cite{OrSo80b} showed
$$\sum_{w\in W} \det{}_V(w)\, x^{k(w)} = \prod_{i=1}^n (x-m_i^*)$$
for some non-negative integers $m_1^*\leq\ldots\leq m_n^*$, the
{\em coexponents of $W$} (see Lehrer--Michel \cite{LeMi} for a generalization,
and Kusuoka \cite{Ku77}, Orlik--Solomon \cite{OrSo83, OrSo85} for versions
over finite fields). \par
Let $\chi$ be an irreducible character of $W$. The {\em fake degree of $\chi$}
is the polynomial
$$R_\chi:=\sum_{d\geq0} \langle K[V]_W^d,\chi^*\rangle_W\,x^d
  =\frac{1}{|W|}\sum_{w\in W}\frac{\chi(w)}{\det_V(xw-1)^*}
  \prod_{i=1}^n(x^{d_i}-1)\ \in\ZZ[x],$$
that is, the graded multiplicity of $\chi^*$ in the $W$-module $K[V]_W$. Thus,
in particular, $R_\chi$ specialises to the degree $\chi(1)$ at $x=1$. The
{\em exponents $(e_i(\chi)\mid 1\leq i\leq\chi(1))$ of an irreducible
character $\chi$ of $W$} are defined by the formula
$R_\chi=\sum_{i=1}^{\chi(1)} x^{e_i(\chi)}$. The exponents $m_i$ of $W$ are
now just the exponents of the contragradient of the reflection representation
$\rho^*:=\tr_{V^*}$, that is, $R_{\rho^*}=\sum_{i=1}^n x^{m_i}$. Dually, the
coexponents are the exponents of $\rho$. In particular, for real reflection
groups exponents and coexponents coincide.
In general $N^*(W):=\sum_{i=1}^n m_i^*$ equals the number of reflecting
hyperplanes of $W$. The $d_i^*:=m_i^*-1$ are sometimes called the
{\em codegrees of $W$}. \par
If $W$ is a Weyl group (see Section~\ref{Mse24}), the fake degrees constitute
a first approximation to
the degrees of principal series unipotent characters of finite groups of Lie
type with Weyl group~$W$. See also Section~\ref{subsec:4.fake} for further
properties.
\end{para}

\begin{para}{Reflection data} \label{subsec:1.refldata}
In the general theory of finite groups of Lie type (where an algebraic
group comes with an action of a Frobenius map) as well as in the study of 
Levi subgroups it is natural to consider reflection groups together with an 
automorphism $\phi$ normalising the reflection representation; see the 
survey article Brou\'e--Malle \cite{BrMa3}. This leads to the following 
abstract definition. 

A pair $(V,W\phi)$ is called a {\em reflection datum} if $V$ is a vector 
space over a subfield $K\subseteq \CC$ and $W\phi$ is a coset in $\GL(V)$
of a reflection group $W\subseteq \GL(V)$, where $\phi\in\GL(V)$ normalises 
$W$.

A {\em sub-reflection datum} of a reflection datum $(V,W\phi)$ is a reflection
datum of the form
$(V',W'(w\phi)|_{V'})$, where $V'$ is a subspace of $V$, $W'$ is a reflection
subgroup of $N_W(V')|_{V'}$ stabilising $V'$ (hence, a reflection subgroup of
$N_W(V')/W_{V'}$), and $w\phi$ is an element of $W\phi$ stabilising $V'$ and
normalising $W'$. A {\em Levi sub-reflection datum of $(V,W\phi)$} is a
sub-reflection datum of the form $(V,W_{V'}(w\phi))$ for some subspace
$V'\leq V$ (note that, by Theorem~\ref{Gthm:steinberg}, $W_{V'}$ is indeed a
reflection subgroup of $W$). A {\em torus of $\GG$} is a sub-reflection datum
with trivial reflection group. \par
Let $\GG=(V,W\phi)$ be a reflection datum. Then $\phi$ acts naturally on the
symmetric algebra $K[V]$. It is possible to choose basic invariants
$f_1,\ldots,f_n\in K[V]^W$, such that $f_i^\phi=\epsilon_i f_i$ for roots of
unity $\epsilon_1,\ldots,\epsilon_n$. The multiset $\{(d_i,\epsilon_i)\}$ of
{\em generalised degrees of $\GG$} then only depends on $W$ and $\phi$
(see for example Springer \cite[Lemma 6.1]{Spr}). The {\em polynomial order} 
of the reflection datum $\GG=(V,W\phi)$ is by definition the polynomial
$$|\GG|:=\frac{\eps_{\GG}x^{N(W)}}
  {\displaystyle\frac{1}{|W|}\sum_{w\in W}\frac{1}{\det_V(1-xw\phi)^*}}=
  x^{N(W)}\prod_{i=1}^n(x^{d_i}-\epsilon_i),$$
where $\eps_{\GG}:=(-1)^n\eps_1\cdots\eps_n$.
Let $\Phi(x)$ be a cyclotomic polynomial over $K$. A torus
$\TT=(V',(w\phi)|_{V'})$ of $\GG$ is called a {\em $\Phi$-torus} if the
polynomial order of $\TT$ is a power of $\Phi$. \par
Reflection data can be thought of as the skeletons of finite reductive groups.
\end{para}

\begin{para}{Regular elements} \label{subsec:1.regular}  
In this section we present results which show that certain subgroups
respectively subquotients of reflection groups are again reflection groups.
Let $(V,w\phi)$ be a reflection datum over $K=\CC$. For $w\phi\in W\phi$
and a root of unity $\zeta\in \CC^\times$ write
$$V(w\phi,\zeta):=\{v\in V\mid w\phi\cdot v=\zeta v\}$$
for the $\zeta$-eigenspace of $w\phi$. Note that $(V(w\phi,\zeta),w\phi)$ is
an $(x-\zeta)$-torus of $(V,w\phi)$ in the sense defined above. These
$(x-\zeta)$-tori for fixed $\zeta$ satisfy a kind of Sylow theory. Let
$f_1,\cdots,f_n$ be a set of basic invariants for $W$ and $H_i$ the surface
defined by $f_i=0$. Springer \cite{Spr} proves that
$$\bigcup_{w\phi\in W\phi} V(w\phi,\zeta)=
  \bigcap_{i: \epsilon_i\zeta^{d_i}=1} H_i\,,$$
the irreducible components of this algebraic set are just the maximal
$V(w\phi,\zeta)$, $W$ acts transitively on these components, and their common
dimension is just the number $a(d,\phi)$ of indices $i$ such that
$\epsilon_i\zeta^{d_i}=1$, where $d$ denotes the order of~$\zeta$. (Note that
$a(d,\phi)$ only depends on $d$, not on $\zeta$ itself.) From this
he obtains:
\end{para}

\begin{thm}[\rm Springer \protect{\cite[Th.~3.4 and~6.2]{Spr}}]
 \label{Gthm:springersyl}
 Let $(V,W\phi)$ be a reflection datum over $\CC$, $\zeta$ a primitive $d$th
 root of unity. Then:
 \begin{itemize}
 \item[(i)] $\max \{\dim V(w\phi,\zeta)\mid w\in W\} = a(d,\phi)$.
 \item[(ii)] For any $w\in W$ there exists a $w'\in W$ such that
 $V(w\phi,\zeta)\subseteq V(w'\phi,\zeta)$ and $V(w'\phi,\zeta)$ has maximal
 dimension.
 \item[(iii)] If $\dim V(w\phi,\zeta)=\dim V(w'\phi,\zeta)=a(d,\phi)$ then
 there exists a $u\in W$ with $u\cdot V(w\phi,\zeta)=V(w'\phi,\zeta)$.
 \end{itemize}
\end{thm}

This can be rephrased as follows: Let $K$ be a subfield of $\CC$, $\Phi$ a
cyclotomic polynomial over $K$. A torus $\TT$ of $\GG$ is called a
{\em $\Phi$-Sylow torus}, if its order equals the full $\Phi$-part of the
order of
$\GG$. Then $\Phi$-tori of $\GG$ satisfy the three statements of Sylow's
theorem. (For an analogue of the statement on the number of Sylow subgroups see
Brou\'e--Malle--Michel \cite[Th.~5.1(4)]{BMM99}). \par
This can in turn be used to deduce a Sylow theory for tori in finite groups of
Lie type (see Brou\'e--Malle \cite{BrMa}).

A vector $v\in V$ is called {\em regular (for $W$)} if it is not contained in
any reflecting hyperplane, i.e. (by Theorem~\ref{Gthm:steinberg}), if its
stabiliser $W_v$ is trivial. Let $\zeta\in\CC$ be a root of unity. An element
$w\phi\in W\phi$ is {\em $\zeta$-regular} if $V(w\phi,\zeta)$ contains a
regular vector. By definition, if $w\phi$ is regular for some root of unity,
then so is any power of
$w\phi$. If $\phi=1$, the Theorem~\ref{Gthm:steinberg} of Steinberg implies
that the orders of $w$ and $\zeta$ coincide. An integer $d$ is a
{\em regular number for $W$} if it is the order of a regular element of
$W$. 

\begin{thm} [\rm Springer \protect{\cite[Th.~6.4 and Prop.~4.5]{Spr}}]
 \label{Gthm:springerreg}
 Let $w\phi\in W\phi$ be $\zeta$-regular of order $d$. Then:
 \begin{itemize}
 \item[(i)] $\dim V(w\phi,\zeta) = a(d,\phi)$.
 \item[(ii)] The centraliser of $w\phi$ in $W$ is isomorphic to a reflection
 group in $V(w\phi,\zeta)$ whose degrees are the $d_i$ with
 $\eps_i\zeta^{d_i}=1$.
 \item[(iii)] The elements of $W\phi$ with property (i) form a single conjugacy
 class under $W$.
 \item[(iv)] Let $\phi=1$ and let $\chi$ be an irreducible character of $W$.
 Then the eigenvalues of $w$ in a representation with character $\chi$ are
  $(\zeta^{e_i(\chi)}\mid 1\leq i\leq \chi(1))$.
 \end{itemize}
\end{thm}

In particular, it follows from~(iv) that the eigenvalues of a $\zeta$-regular
element $w$ on $V$ are $(\zeta^{m_i^*}\mid 1\leq i\leq n)$. \par

Interestingly enough, the normaliser modulo centraliser of arbitrary Sylow
tori of reflection data are naturally reflection groups, as the following
generalisation of the previous result shows:

\begin{thm} [\rm Lehrer--Springer \cite{LeSp, LeSp2}]\label{Gthm:lehrerspringer}
 Let $w\in W$ and $\tilde V:=V(w\phi,\zeta)$ be such that $(\tilde V,w\phi)$
 is a $\Phi$-Sylow torus.
 Let $N:=\{w'\in W\mid w'\cdot\tilde V=\tilde V\}$ be the normaliser,
 $C:=\{w'\in W\mid w'\cdot v=v\text{ for all }v\in \tilde V\}$ the centraliser
 of $\tilde V$. \par
 \begin{itemize}
 \item[(i)] Then $N/C$ acts as a reflection group on $\tilde V$, with
 reflecting hyperplanes the intersections with $\tilde V$ of those of $W$.
 \item[(ii)] A set of basic invariants of $N/C$ is given by the restrictions
 to $\tilde V$ of those $f_i$ with $\epsilon_i\zeta^{d_i}=1$.
 \item[(iii)] If $W$ is irreducible on $V$, then so is $N/C$ on $\tilde V$.
 \end{itemize}
\end{thm}

In the case of regular elements, the second assertion of~(i) goes back to
Lehrer \cite[5.8]{Le95}, Denef--Loeser \cite{DeLo95}; see also 
Brou\'e--Michel \cite[Prop.~3.2]{BrMi}.

\begin{para}{The Shephard--Todd classification} \label{subsec:1.Class} 
Let $V$ be a finite-dimensional complex vector space and $W\leq\GL(V)$ a
reflection group. Since $W$ is finite, the representation on $V$ is completely
reducible, and $W$ is the direct product of irreducible reflection subgroups.
Thus, in order to determine all reflection groups over $\CC$, it is sufficient
to classify the irreducible ones. This was achieved by Shephard and Todd
\cite{ShTo}. \par
To describe this classification, first recall that a subgroup
$W\leq\GL(V)$ is called imprimitive if there exists a direct sum
decomposition $V=V_1\oplus\ldots\oplus V_k$
with $k>1$ stabilised by $W$ (that is, $W$ permutes the summands).
The bulk of irreducible complex reflection groups consists of imprimitive ones.
For any $d,e,n\geq1$ let $G(de,e,n)$ denote the group of monomial
$n\times n$-matrices (that is, matrices with precisely one non-zero entry in
each row and column) with non-zero entries in the set of $de$th roots of
unity, such that the product over these entries is a $d$th root of unity. \par
Explicit generators may be
chosen as follows: $G(d,1,n)$ is generated on $\CC^n$ with standard Hermitean
form) by the reflection $t_1$ of order~$d$ with root the first standard basis
vector $b_1$ and by the permutation matrices $t_2,\ldots,t_n$ for the
transpositions $(1,2),(2,3),\ldots,(n-1,n)$. For $d>1$ this is an irreducible
reflection group, isomorphic to the wreath product $C_d\wr\fS_n$ of the cyclic
group of order~$d$ with the symmetric group $\fS_n$, where the base group is
generated by the reflections of order~$d$ with roots the standard basis
vectors, and a complement consists of all permutation matrices. \par
Let $\gamma_d: G(d,1,n)\longrightarrow\CC^\times$ be the linear character of
$G(d,1,n)$ obtained by tensoring the determinant on $V$ with the sign
character on the quotient $\fS_n$. Then for any $e>1$ we have
$$G(de,e,n):=\ker(\gamma_{de}^{d})\leq G(de,1,n).$$
This is an irreducible reflection subgroup of $G(de,1,n)$ for all $n\geq2$,
$d\geq1$, $e\geq2$, except for $(d,e,n)=(2,2,2)$. It is generated by the
reflections
$$t_1^{e},t_1^{-1}t_2t_1,t_2,t_3,\ldots,t_n,$$
where the first generator is redundant if $d=1$. Clearly,
$G(de,e,n)$ stabilises the decomposition $V=\CC b_1\oplus\ldots\oplus\CC b_n$
of $V$, so it is imprimitive for $n>1$. The order of $G(de,e,n)$ is given by
$d^ne^{n-1} n!$. Using the wreath product structure it is easy to show that the
only isomorphisms among groups in this series are $G(2,1,2)\cong G(4,4,2)$,
and $G(de,e,1)\cong G(d,1,1)$ for all $d,e$, while $G(2,2,2)$ is reducible.
All these are isomorphisms of reflection groups. \par
In its natural action on $\QQ^n$ the symmetric group $\fS_n$ stabilises
the 1-dimensional subspace consisting of vectors with all coordinates equal
and the $(n-1)$-dimensional subspace consisting of those vectors whose
coordinates add up to~0. In its action on the latter, $\fS_n$ is an irreducible
and primitive reflection group. The classification result may now be stated
as follows (see also Cohen \cite{Co}):
\end{para}

\begin{thm} [\rm Shephard--Todd \cite{ShTo}] \label{Gthm:shephard}
 The irreducible complex reflection groups are the groups $G(de,e,n)$, for
 $de\geq2$, $n\geq1$, $(de,e,n)\ne(2,2,2)$, the groups $\fS_n$ ($n\geq2$) in
 their
 $(n-1)$-dimensional natural representation, and 34 further primitive groups.
 \par
 Moreover, any irreducible $n$-dimensional complex reflection group has
 a generating set of at most $n+1$ reflections.
\end{thm}

The primitive groups are usually denoted by $G_4,\ldots,G_{37}$, as in the
original article \cite{ShTo} (where the first three indices were reserved for
the families of imprimitive groups $G(de,e,n)$ ($de,n\geq2$), cyclic groups
$G(d,1,1)$ and symmetric groups $\fS_{n+1}$). 
An $n$-dimensional irreducible reflection groups generated by $n$ of its
reflections is called {\em well-generated}. The groups for which this fails
are the imprimitive groups $G(de,e,n)$, $d,e,n>1$, and the primitive
groups
$$G_i\qquad\text{with }i\in\{7,11,12,13,15,19,22,31\}.$$
By construction $G(de,e,n)$
contains the well-generated group $G(de,de,n)$. More generally, the
classification implies the following, for which no a priori proof is known: 

\begin{cor} \label{Gcor:irrwell}
 Any irreducible complex reflection group $W\leq\GL(V)$ contains a
 well-generated reflection subgroup $W'\leq W$ which is still irreducible
 on $V$.
\end{cor}

\begin{table}[htbp] \caption{Irreducible complex reflection groups}
\label{Gtab1} \begin{center}
$\begin{array}{c}
\renewcommand{\arraystretch}{1.1}
\begin{array}{lll@{\hspace{0pt}}c} \multicolumn{4}{c}{\text{Infinite series}}\\
\hline W     &\text{degrees}   &\text{codegrees}  &K_W \\ \hline
 G(d,1,n) \; (d{\geq}2,n{\geq}1) &d,2d,\dots,nd &* &\QQ(\zeta_{d})\\
G(de,e,n)\;(d,e,n{\geq}2)  &ed,2ed,\dots,(\!n\!{-}\!1\!)ed,nd
        &0,ed,\dots,(\!n\!{-}\!1\!)ed\!\! &\QQ(\zeta_{de}) \\
 G(e,e,n)\; (e{\geq} 2,n{\geq}3)  &e,2e,\dots, (\!n\!{-}\!1\!)e,n
&* &\QQ(\zeta_e) \\
G(e,e,2)\; (e{\geq} 3) &2,e &* &\QQ(\zeta_e\!{+}\!\zeta_e^{-1})\\
\fS_{n{+}1}\ (n{\geq}1)&2,3,\dots,n{+}1  &* &\QQ \\\hline\end{array}\\ \\
\begin{array}{llllc} \multicolumn{5}{c}{\text{Exceptional groups}}\\
\hline
  W     &\text{degrees}   &\text{codegrees}  &K_W  &W/Z(W)\cr
 \hline
 G_4    &{\bf4,6}   &*     &\QQ(\zeta_3)    &\fA_4 \cr
 G_5    &6,{\bf12}  &*     &\QQ(\zeta_3)    &\fA_4 \cr
 G_6    &4,{\bf12}  &*     &\QQ(\zeta_{12}) &\fA_4 \cr
 G_7    &12,{\bf12} &0,12  &\QQ(\zeta_{12}) &\fA_4 \cr
 G_8    &\bf{8,12}  &*     &\QQ(i)          &\fS_4 \cr
 G_9    &8,{\bf24}  &*     &\QQ(\zeta_{8})  &\fS_4 \cr
 G_{10} &12,{\bf24} &*     &\QQ(\zeta_{12}) &\fS_4 \cr
 G_{11} &24,{\bf24} &0,24  &\QQ(\zeta_{24}) &\fS_4 \cr
 G_{12} &{\bf6,8}   &0,10  &\QQ(\sqrt{-2})  &\fS_4 \cr
 G_{13} &8,{\bf12}  &0,16  &\QQ(\zeta_8)    &\fS_4 \cr
 G_{14} &6,{\bf24}  &*     &\QQ(\zeta_3,\sqrt{-2})  &\fS_4 \cr
 G_{15} &{\bf12},24 &0,24  &\QQ(\zeta_{24}) &\fS_4 \cr
 G_{16} &{\bf20,30} &*     &\QQ(\zeta_5)    &\fA_5 \cr
 G_{17} &20,{\bf60} &*     &\QQ(\zeta_{20}) &\fA_5 \cr
 G_{18} &30,{\bf60} &*     &\QQ(\zeta_{15}) &\fA_5 \cr
 G_{19} &60,{\bf60} &0,60  &\QQ(\zeta_{60}) &\fA_5 \cr
 G_{20} &{\bf12,30} &*     &\QQ(\zeta_3,\sqrt{5})   &\fA_5 \cr
 G_{21} &12,{\bf60} &*     &\QQ(\zeta_{12},\sqrt{5}) &\fA_5 \cr
 G_{22} &{\bf12,20} &0,28  &\QQ(i,\sqrt{5}) &\fA_5 \cr
 G_{23} &2,{\bf6,10} &*    &\QQ(\sqrt{5})  &\fA_5 \cr
 G_{24} &4,{\bf6,14}  &*   &\QQ(\sqrt{-7}) &\GL_3(2) \cr
 G_{25} &6,\bf{9,12}  &*   &\QQ(\zeta_3)   &3^2\colon \SL_2(3)  \cr
 G_{26} &6,12,{\bf18} &*   &\QQ(\zeta_3) &3^2\colon \SL_2(3) \cr
 G_{27} &6,12,{\bf30} &*   &\QQ(\zeta_3,\sqrt{5}) &\fA_6 \cr
 G_{28} &2,6,{\bf8,12} &*  &\QQ & 2^4\colon(\fS_3\times\fS_3)  \cr
 G_{29} &4,8,12,{\bf20} &* &\QQ(i) &2^4\colon\fS_5  \cr
 G_{30} &2,{\bf12,20,30} &*  &\QQ(\sqrt{5}) & \fA_5\wr 2 \cr
 G_{31} &8,12,{\bf20,24} & 0,12,16,28 &\QQ(i) & 2^4\colon\fS_6  \cr
 G_{32} &12,18,{\bf24,30} &* &\QQ(\zeta_3) &\UU_4(2) \cr
 G_{33} &4,6,{\bf10},12,{\bf18} &* &\QQ(\zeta_3) &\OO_5(3) \cr
 G_{34} &6,12,18,24,30,{\bf42}  &* &\QQ(\zeta_3) &\OO_6^-(3).2 \cr
 G_{35} &2,5,6,{\bf8,9,12} &* &\QQ & \OO_6^-(2) \cr
 G_{36} &2,6,8,10,12,{\bf14,18} &* &\QQ &\OO_7(2)  \cr
 G_{37} &2,8,12,14,18,{\bf20,24,30}  &*  &\QQ  &\OO_8^+(2).2 \cr
 \hline
\end{array}\end{array}$
\end{center}
\end{table}

The primitive groups
$G_4,\ldots,G_{37}$ occur in dimensions 2 up to~8. In Table~\ref{Gtab1} 
we collect some data on the irreducible complex reflection groups. (These and
many more data for complex reflection groups have been implemented by Jean
Michel into the {\sf CHEVIE}-system \cite{chev}.) In the first part,
the dimension is always equal to $n$, in the second it can be read of from the
number of degrees. We give the degrees, the codegrees in case they are not
described by the following Theorem~\ref{Gthm:orlik} (that is, if $W$ is not
well-generated), and the character field $K_W$ of the
reflection representation. For the exceptional groups we also give the
structure of $W/Z(W)$ and we indicate the regular degrees by boldface (that
is, those degrees which are regular numbers for $W$, see
Cohen \cite[p.~395 and p.~412]{Co} and Springer \cite[Tables~1--6]{Spr}). The 
regular degrees for the infinite series are: $n,n+1$ for $\fS_{n+1}$, $dn$ 
for $G(de,e,n)$ with $d>1$, $(n-1)e$ for $G(e,e,n)$ with $n|e$, and 
$(n-1)e,n$ for $G(e,e,n)$ with $n{\not|}e$. Lehrer and Michel 
\cite[Th.~3.1]{LeMi} have shown that an integer is a regular number if and only
if it divides as many degrees as codegrees. \par

Fundamental invariants for most types are given in Shephard--Todd \cite{ShTo}
as well as defining relations and further information on
parabolic subgroups (see also Coxeter \cite{Cox57}, Shephard \cite{Shep} 
and Brou\'e--Malle--Rouquier \cite{BMR} for presentations, and the tables 
in Cohen \cite{Co} and Brou\'e--Malle--Rouquier \cite[Appendix~2]{BMR}).

If the irreducible reflection group $W$ has an invariant of degree~2, then it
leaves invariant a non-degenerate quadratic form, so the representation may be
realised over the reals. Conversely, if $W$ is a real reflection group, then
it leaves invariant a quadratic form. Thus the real irreducible reflection
groups are precisely those with $d_1=2$, that is, the infinite series
$G(2,1,n)$, $G(2,2,n)$, $G(e,e,2)$ and $\fS_{n+1}$, and the six exceptional
groups $G_{23},G_{28},G_{30},G_{35},G_{36},G_{37}$ (see Section~\ref{Msec23}).

The degrees and codegrees of a finite complex reflection group satisfy some
remarkable identities. As an example, let us quote the following result, for
which at present only a case-by-case proof is known:

\begin{thm} [\rm Orlik--Solomon \cite{OrSo80b}] \label{Gthm:orlik}
 Let $W$ be an irreducible complex reflection group in dimension~$n$. Then the
 following are equivalent:
 \begin{itemize}
 \item[(i)] $d_i+d_{n-i+1}^*=d_n$ for $i=1,\ldots,n$,
 \item[(ii)] $N+N^* = nd_n$,
 \item[(iii)] $d_i^*<d_n$ for $i=1,\ldots,n$,
 \item[(iv)] $W$ is well-generated.
 \end{itemize}
\end{thm}

See also Terao--Yano \cite{TeYa} for a partial explanation.\par

\begin{table}[htbp] \caption{Exceptional twisted reflection data} 
\label{Gtab3}
\begin{center}
$\begin{array}{llllcc} \hline
  W     &d_i   &\epsilon_i   & \text{field}& \text{origin of }\phi\cr
 \hline
 G(4,2,2) &4,4      & 1,\zeta_3  &\QQ(i)  &<G_6\cr
 G(3,3,3) &3,6,3    & 1,1,-1     &\QQ(\zeta_3)   &<G_{26}\cr
 G(2,2,4) &2,4,4,6  & 1,\zeta_3,\zeta_3^2,1  &\QQ    &<G_{28}\cr
 G_5      &6,12     & 1,-1       &\QQ(\zeta_3,\sqrt{-2}) &<G_{14}\cr
 G_7      &12,12    & 1,-1       &\QQ(\zeta_{12})   &<G_{10}\cr
 G_{28}   &2,6,8,12 & 1,-1,1,-1  &\QQ(\sqrt{2})    &\text{graph aut.}\cr
 \hline
\end{array}$
\end{center}
\end{table}

From the Shephard--Todd classification, it is straightforward to obtain a
classification of reflection data. An easy argument
allows to reduce to the case where $W$ acts irreducibly on $V$. Then either up
to scalars $\phi$ can be chosen to be a reflection, or $W=G_{28}$ is the real
reflection group of type $F_4$, and $\phi$ induces the graph automorphism on
the $F_4$-diagram (see Brou\'e--Malle--Michel \cite[Prop.~3.13]{BMM99}). An
infinite series of examples is obtained from the embedding of $G(de,e,n)$ into
$G(de,1,n)$ (which is the full projective normaliser in all but finitely many
cases). Apart from this, there are only six further cases, which we list in
Table~\ref{Gtab3}.

\section{Real reflection groups} \label{sec12}
In this section we discuss in more detail the special case where $W$ is
a {\em real} reflection group. This is a well-developped theory, and 
there are several good places to learn about real reflection groups: the 
classical Bourbaki volume \cite{Bou}, the very elementary text by 
Benson--Grove \cite{BeGr}, the relevant chapters in Curtis--Reiner \cite{CR2},
Hiller \cite{Hiller82}, and Humphreys \cite{Humphreys2}. Various pieces of 
the theory have also been recollected in a concise way in  articles by 
Steinberg \cite{SteinColl}. The exposition here partly follows Geck--Pfeiffer 
\cite[Chap.~1]{GePf}. We shall only present the most basic results 
and refer to the above textbooks and our bibliography for further reading.

\begin{para}{Coxeter groups} \label{Msecc21} Let $S$ be a finite non-empty 
index set and $M=(m_{st})_{s,t\in S}$ be a symmetric matrix such that 
$m_{ss}=1$ for all $s\in S$ and $m_{st}\in \{2,3,4,\ldots\}\cup 
\{\infty\}$ for all $s\neq t$ in $S$.  Such a matrix is called a 
{\em Coxeter matrix}. Now let $W$ be a group containing $S$ as a subset. 
($W$ may be finite or infinite.)
Then the pair $(W,S)$ is called a {\em Coxeter system}, and $W$ is called
a {\em Coxeter group}, if $W$ has a presentation with generators $S$ and 
defining relations of the form
\[ (st)^{m_{st}}=1 \qquad \mbox{for all $s,t \in S$ with $m_{st}<\infty$};\]
in particular, this means that $s^2=1$ for all $s\in S$. Therefore, the
above relations (for $s\neq t$) can also be expressed in the form 
\[ \underbrace{sts\cdots}_{\text{$m_{st}$ times}}=\underbrace{tst
\cdots}_{\text{$m_{st}$ times}} \qquad \mbox{for all $s,t\in S$ with $2\leq 
m_{st}<\infty$}.\]
We say that $C$ is of finite type and that $(W,S)$ is a finite Coxeter system 
if $W$ is a finite group. The information contained in $M$ can be visualised 
by a corresponding 
{\em Coxeter graph}, which is defined as follows. It has vertices labelled 
by the elements of $S$, and two vertices labelled by $s \neq t$ are joined 
by an edge if $m_{st} \geq 3$. Moreover, if $m_{st} \geq 4$, we label the 
edge by~$m_{st}$.  The standard example of a finite Coxeter system  is the 
pair $(\fS_n, \{s_1,\ldots,s_{n-1}\})$ where $s_i=(i,i+1)$ for $1\leq i\leq 
n-1$. The corresponding graph is 
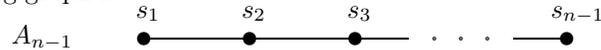
\begin{center}
\begin{picture}(240,20)
\put(  0,08){$A_{n{-}1}$}
\put( 50,10){\circle*{5}}
\put( 47,18){$s_1$}
\put( 50,10){\line(1,0){40}}
\put( 90,10){\circle*{5}}
\put( 87,18){$s_2$}
\put( 90,10){\line(1,0){40}}
\put(130,10){\circle*{5}}
\put(127,18){$s_3$}
\put(130,10){\line(1,0){20}}
\put(160,10){\circle{1}}
\put(170,10){\circle{1}}
\put(180,10){\circle{1}}
\put(190,10){\line(1,0){20}}
\put(210,10){\circle*{5}}
\put(204,18){$s_{n{-}1}$}
\end{picture}
\end{center}
\end{para}

Coxeter groups have a rich combinatorial structure. A basic tool is the 
{\em length function} $l \colon W \rightarrow {\NN}_0$, which is defined  as 
follows. Let $w \in W$. Then $l(w)$ is the length of a shortest possible
expression $w=s_1 \cdots s_k$ where $s_i\in S$. An expression of $w$
of length $l(w)$ is called a {\em reduced expression} for~$w$. We have
$l(1)=0$ and $l(s)=1$ for $s\in S$. Here is a key result  about Coxeter
groups.

\begin{thm}[Matsumoto \cite{Matsum}; see also Bourbaki \cite{Bou}] 
\label{Mmatsum} Let $(W,S)$ be a Coxeter system and $\cM$ be a monoid, with 
multiplication $\star \colon \cM \times \cM \rightarrow \cM$. Let $f \colon 
S \rightarrow \cM$ be a map such that
\[\underbrace{f(s)\star f(t)\star f(s)\star \cdots}_{\text{$m_{st}$ times}}=
\underbrace{f(t)\star f(s)\star f(t)\star \cdots}_{\text{$m_{st}$ times}}\]
for all $s\neq t$ in $S$ such that $m_{st} <\infty$. Then there exists a
unique map $F \colon W \rightarrow \cM$ such that $F(w)=f(s_1) \star  
\cdots \star f(s_k)$ whenever $w=s_1 \cdots s_k$ ($s_i \in S$) is reduced.
\end{thm}

Typically, Matsumoto's theorem can be used to show that certain
constructions with reduced expressions of elements of $W$ actually
do not depend on the choice of the reduced expressions. We give two 
examples.

\medskip
{\em (1) Let $w \in W$ and take a reduced expression $w=s_1 \cdots s_k$ 
with $s_i\in S$. Then the set $\{s_1,\ldots,s_k\}$ does not depend on
the choice of the reduced expression.}
\medskip

(Indeed, let $\cM$ be the monoid whose elements are the subsets of~$S$ and 
product given by $A\star B:=A\cup B$. Then the assumptions of Matsumoto's 
theorem are satisfied for the map 
$f \colon S \rightarrow \cM$, $s \mapsto \{s\}$, and this yields the
required assertion.)

\medskip
{\em (2)  Let $w \in W$ and fix a reduced expression $w=s_1 \cdots s_k$ ($s_i 
\in S$). Consider the set of all subexpressions:
\[ \cS(w):=\{y\in W\mid y=s_{i_1} \cdots s_{i_l} \mbox{where $l \geq 0$ 
and $1 \leq i_1 < \dots < i_l \leq k$}\}.\] 
Then $\cS(w)$ does not depend on the choice of the reduced expression for $w$.}
\medskip

(Indeed, let $\cM$ be the monoid whose elements are the subsets
of~$W$ and product given by $A \star B:=\{ab \mid a \in A, 
b \in B\}$ (for $A,B \subseteq W$). Then the assumptions of Matsumoto's 
theorem are satisfied for the map $f \colon S \rightarrow \cM$, 
$s \mapsto \{1,s\}$, and this yields the required assertion.)

We also note that the so-called {\em Exchange Condition} and the 
{\em Cancellation Law} are further consequences of the above results.
The ``Cancellation Law'' states that, given $w\in W$ and an expression 
$w=s_1 \cdots s_k$ ($s_i\in S$) which is not reduced, one can obtain a
reduced expression of $w$ by simply cancelling some of the factors in
the given expression.  This law together with (2) yields that the relation 
\[ y\leq w \quad \stackrel{\text{def}}{\Longleftrightarrow} \quad y
\in \cS(w)\]
is a partial order on $W$, called the {\em Bruhat--Chevalley order}. 
This ordering has been extensively studied; see, for example, 
Verma \cite{Verma1}, Deodhar \cite{Deodhar},  Bj\"orner \cite{Bjo}, 
Lascoux--Sch\"utzenberger \cite{LaSch} and Geck--Kim \cite{GeKim}. 
By Chevalley \cite{Chev1}, it is related to the Bruhat decomposition 
in algebraic groups; we will explain this result in (\ref{Mbnpair})
below. 

\begin{para}{Cartan matrices} \label{Msec22} Let $M=(m_{st})_{s,t\in S}$ be a 
Coxeter matrix as above. We can also associate with $M$ a group generated
by reflections. This is done as follows. Choose a matrix $C=(c_{st})_{s,t 
\in S}$ with entries in~$\RR$ such that the following conditions are satisfied:
\begin{itemize}
\item[(C1)] For $s \neq t$ we have $c_{st} \leq 0$; furthermore, $c_{st} 
\neq 0$ if and only if $c_{ts} \neq 0$.
\item[(C2)] We have $c_{ss}=2$ and, for $s \neq t$, we have $c_{st} 
c_{ts}=4\cos^2(\pi/m_{st})$. 
\end{itemize}
Such a matrix $C$ will be called a {\em Cartan matrix} associated with $M$.
For example, we could simply take $c_{st}:=-2\cos(\pi/m_{st})$ for all
$s,t\in S$; this may be called the {\em standard Cartan matrix} associated
with $M$. We always have $0 \leq c_{st}c_{ts} \leq 4$. Here are some
values for the product $c_{st}c_{ts}$:
\[ \begin{array}{l|ccccccc} \hline m_{st} & 2 & 3 & 4 & 5 & 6 & 8 & \infty \\ 
\hline c_{st}c_{ts} & 0 & 1 & 2 & (3+\sqrt{5})/2 &3 & 2+\sqrt{2} &4 \\ \hline 
\end{array}\]
Now let $V$ be an ${\RR}$-vector space of dimension $|S|$, with a fixed 
basis $\{\alpha_s \mid s \in S\}$. We define a linear action of the 
elements in $S$ on $V$ by the rule:
\[s \colon V \rightarrow V, \quad \alpha_t \mapsto \alpha_t-c_{st} 
\alpha_s \quad (t \in S).\]
It is easily checked that $s\in \GL(V)$ has order~$2$ and precisely one
eigenvalue $-1$ (with eigenvector $\alpha_s$). Thus, $s$ is a reflection
with root $\alpha_s$. We then define 
\[ W=W(C):=\langle S \rangle \subset \GL(V);\]
thus, if $|W|<\infty$ is finite, then $W$ will be a real reflection group.
Now we can state the following basic result.
\end{para}

\begin{thm}[Coxeter \cite{Cox34,Cox35}] \label{Mthm1}
 Let $M=(m_{st})_{s,t\in S}$ be a Coxeter matrix and $C$ be a Cartan matrix
 associated with $M$. Let $W(C)= \langle S\rangle \subseteq \GL(V)$ be the
 group constructed as in (\ref{Msec22}). Then the pair $(W(C),S)$ is a Coxeter
 system. Furthermore, the group $W(C)$ is finite if and only if 
 \begin{equation*}
  \mbox{the matrix } \left(-\cos(\pi/m_{st})\right)_{s,t \in S} 
  \mbox{ is positive-definite},\tag{$*$}
 \end{equation*}
 i.e., we have $\det (-\cos(\pi/m_{st}))_{s,t \in J}>0$ 
 for every subset $J \subseteq S$. All finite real reflection groups
 arise in this way.
\end{thm}

The fact that $(W(C),S)$ is a Coxeter system is proved in 
\cite[1.2.7]{GePf}. The finiteness condition can be found in
Bourbaki \cite[Chap.~V, \S 4, no.~8]{Bou}. Finally, the fact that all
finite real reflection groups arise in this way is established in
\cite[Chap.~V, \S 3, no.~2]{Bou}. The ``note historique'' in \cite{Bou} 
contains a detailed account of the history of the above result. 

\begin{para}{Classification of finite Coxeter groups}\label{Msec23}
(See also Section~\ref{subsec:1.Class}.)
Let $M=(m_{st})_{s,t \in S}$ be a Coxeter matrix. We say that $M$ is
{\em decomposable} if there is a partition $S=S_1 \amalg S_2$ with $S_1,S_2 
\neq \varnothing$ and such that $m_{st}=2$ whenever $s \in S_1$, 
$t \in S_2$. If $C=(c_{st})$ is any Cartan matrix associated with $M$,
then this condition translates to: $c_{st}=c_{ts}=0$ whenever $s \in S_1$, 
$t \in S_2$. Correspondingly, we also have a direct sum decomposition 
$V=V_1 \oplus V_2$ where $V_1$ has basis $\{\alpha_s \mid s \in S_1\}$ and 
$V_2$ has basis $\{\alpha_s \mid s \in S_2\}$. Then it easily follows that 
we have an isomorphism 
\[W(C) \stackrel{\sim}{\longrightarrow} W(C_1) \times W(C_2), \quad w \mapsto 
(w|_{V_1},w|_{V_2}).\]
In this way, the study of the groups $W(C)$ is reduced to the case 
where~$C$ is {\em indecomposable} (i.e., there is no partition 
$S=S_1\amalg S_1$ as above). If this holds, we call the corresponding 
Coxeter system $(W,S)$ an {\em irreducible Coxeter system}.
\end{para}

\begin{thm} \label{classdyn} The Coxeter graphs of the indecomposable 
Coxeter matrices $M$ such that condition ($*$) in Theorem~\ref{Mthm1} 
holds are precisely the graphs in Table~\ref{Mtab1}.
\end{thm}

\begin{table}[htbp] \caption{Coxeter graphs of irreducible finite 
Coxeter groups} \label{Mtab1} 
\begin{center}
\makeatletter
\vbox{\begin{picture}(345,160)
\put( 10, 25){$E_7$}
\put( 50, 25){\@dot{5}}
\put( 48, 30){$1$}
\put( 50, 25){\line(1,0){20}}
\put( 70, 25){\@dot{5}}
\put( 68, 30){$3$}
\put( 70, 25){\line(1,0){20}}
\put( 90, 25){\@dot{5}}
\put( 88, 30){$4$}
\put( 90, 25){\line(0,-1){20}}
\put( 90,  5){\@dot{5}}
\put( 95,  3){$2$}
\put( 90, 25){\line(1,0){20}}
\put(110, 25){\@dot{5}}
\put(108, 30){$5$}
\put(110, 25){\line(1,0){20}}
\put(130, 25){\@dot{5}}
\put(128, 30){$6$}
\put(130, 25){\line(1,0){20}}
\put(150, 25){\@dot{5}}
\put(148, 30){$7$}

\put(190, 25){$E_8$}
\put(220, 25){\@dot{5}}
\put(218, 30){$1$}
\put(220, 25){\line(1,0){20}}
\put(240, 25){\@dot{5}}
\put(238, 30){$3$}
\put(240, 25){\line(1,0){20}}
\put(260, 25){\@dot{5}}
\put(258, 30){$4$}
\put(260, 25){\line(0,-1){20}}
\put(260,  5){\@dot{5}}
\put(265,  3){$2$}
\put(260, 25){\line(1,0){20}}
\put(280, 25){\@dot{5}}
\put(278, 30){$5$}
\put(280, 25){\line(1,0){20}}
\put(300, 25){\@dot{5}}
\put(298, 30){$6$}
\put(300, 25){\line(1,0){20}}
\put(320, 25){\@dot{5}}
\put(318, 30){$7$}
\put(320, 25){\line(1,0){20}}
\put(340, 25){\@dot{5}}
\put(338, 30){$8$}

\put( 10, 84){$I_2(m)$}
\put( 10, 74){$\scriptstyle{m \geq 5}$}
\put( 50, 85){\@dot{5}}
\put( 48, 90){$1$}
\put( 56, 89){$\scriptstyle{m}$}
\put( 50, 85){\line(1,0){20}}
\put( 70, 85){\@dot{5}}
\put( 68, 90){$2$}

\put(103, 84){$F_4$}
\put(130, 85){\@dot{5}}
\put(128, 90){$1$}
\put(130, 85){\line(1,0){20}}
\put(150, 85){\@dot{5}}
\put(148, 90){$2$}
\put(150, 85){\line(1,0){20}}
\put(158, 88){$\scriptstyle{4}$}
\put(170, 85){\@dot{5}}
\put(168, 90){$3$}
\put(170, 85){\line(1,0){20}}
\put(190, 85){\@dot{5}}
\put(188, 90){$4$}

\put( 10, 55){$H_3$}
\put( 50, 55){\@dot{5}}
\put( 48, 60){$1$}
\put( 50, 55){\line(1,0){20}}
\put( 59, 58){$\scriptstyle{5}$}
\put( 70, 55){\@dot{5}}
\put( 68, 60){$2$}
\put( 70, 55){\line(1,0){20}}
\put( 90, 55){\@dot{5}}
\put( 88, 60){$3$}

\put(117, 55){$H_4$}
\put(150, 55){\@dot{5}}
\put(148, 60){$1$}
\put(150, 55){\line(1,0){20}}
\put(159, 58){$\scriptstyle{5}$}
\put(170, 55){\@dot{5}}
\put(168, 60){$2$}
\put(170, 55){\line(1,0){20}}
\put(190, 55){\@dot{5}}
\put(188, 60){$3$}
\put(190, 55){\line(1,0){20}}
\put(210, 55){\@dot{5}}
\put(208, 60){$4$}

\put(230, 75){$E_6$}
\put(260, 75){\@dot{5}}
\put(258, 80){$1$}
\put(260, 75){\line(1,0){20}}
\put(280, 75){\@dot{5}}
\put(278, 80){$3$}
\put(280, 75){\line(1,0){20}}
\put(300, 75){\@dot{5}}
\put(298, 80){$4$}
\put(300, 75){\line(0,-1){20}}
\put(300, 55){\@dot{5}}
\put(305, 53){$2$}
\put(300, 75){\line(1,0){20}}
\put(320, 75){\@dot{5}}
\put(318, 80){$5$}
\put(320, 75){\line(1,0){20}}
\put(340, 75){\@dot{5}}
\put(338, 80){$6$}

\put( 10,120){$B_n$}
\put( 10,110){$\scriptstyle{n \geq 2}$}
\put( 50,115){\@dot{5}}
\put( 50,120){$1$}
\put( 50,115){\line(1,0){20}}
\put( 59,118){$\scriptstyle{4}$}
\put( 70,115){\@dot{5}}
\put( 68,120){$2$}
\put( 70,115){\line(1,0){30}}
\put( 90,115){\@dot{5}}
\put( 88,120){$3$}
\put(110,115){\@dot{1}}
\put(120,115){\@dot{1}}
\put(130,115){\@dot{1}}
\put(140,115){\line(1,0){10}}
\put(150,115){\@dot{5}}
\put(147,120){$n$}

\put( 10,150){$A_{n}$}
\put( 10,140){$\scriptstyle{n \geq 1}$}
\put( 50,145){\@dot{5}}
\put( 48,150){$1$}
\put( 50,145){\line(1,0){20}}
\put( 70,145){\@dot{5}}
\put( 68,150){$2$}
\put( 70,145){\line(1,0){30}}
\put( 90,145){\@dot{5}}
\put( 88,150){$3$}
\put(110,145){\@dot{1}}
\put(120,145){\@dot{1}}
\put(130,145){\@dot{1}}
\put(140,145){\line(1,0){10}}
\put(150,145){\@dot{5}}
\put(147,150){$n$}

\put(210,127){$D_n$}
\put(210,117){$\scriptstyle{n \geq 4}$}
\put(240,145){\@dot{5}}
\put(245,145){$1$}
\put(240,105){\@dot{5}}
\put(246,100){$2$}
\put(240,145){\line(1,-1){21}}
\put(240,105){\line(1,1){21}}
\put(260,125){\@dot{5}}
\put(258,130){$3$}
\put(260,125){\line(1,0){30}}
\put(280,125){\@dot{5}}
\put(278,130){$4$}
\put(300,125){\@dot{1}}
\put(310,125){\@dot{1}}
\put(320,125){\@dot{1}}
\put(330,125){\line(1,0){10}}
\put(340,125){\@dot{5}}
\put(337,130){$n$}
\end{picture}}
\makeatother
\end{center}
\footnotesize The numbers on the vertices correspond to a chosen labelling 
of the elements of~$S$. 
\end{table}

For the proof of this classification, see 
\cite[Chap.~VI, no.~4.1]{Bou}. The identification with the  groups 
occurring in the Shephard--Todd classification (see 
Theorem~\ref{Gthm:shephard}) is given in the following table.
\[ \begin{array}{cc}\hline
\mbox{Coxeter graph} & \mbox{Shephard--Todd}\\ \hline
A_{n-1} & \fS_n \\
B_n     & G(2,1,n) \\
D_n     & G(2,2,n)\\
I_2(m)  &  G(m,m,2)  \\
H_3     &  G_{23}        \\
H_4     &  G_{30}  \\
F_4     &  G_{28}   \\
E_6     & G_{35}   \\
E_7     &  G_{36}   \\
E_8     &  G_{37} \\ \hline \end{array}\]
Thus, any finite irreducible real reflection group
is the reflection group arising from a Cartan matrix associated with one of 
the graphs in Table~\ref{Mtab1}.

Now, there are a number of results on finite Coxeter groups which can
be formulated in general terms but whose proof requires a case-by-case
verification using the above classification. We mention two such results,
concerning conjugacy classes. 

\begin{thm}[Carter \cite{Carter72}] \label{Min} Let $(W,S)$ be a finite
Coxeter system. Then every element in $W$ is conjugate to its inverse.
More precisely, given $w\in W$, there exist $x,y\in W$ such that
$w=xy$ and $x^2=y^2=1$.
\end{thm}

Every element $x\in W$ such that $x^2=1$ is a product of pairwise 
commuting reflections in $W$. Given $w\in W$ and an expression $w=xy$ as 
above, the geometry of the roots involved in the reflections determining
$x,y$ yields a diagram which can be used to label the conjugacy class of~$w$.
Complete lists of these diagrams can be found in \cite{Carter72}.

\begin{para}{Conjugacy classes and the length function}\label{Mfurther}
Let $(W,S)$ be a Coxeter system and $C$ be a conjugacy class in $W$. We
will be interested in studying how conjugation inside $C$ relates to
the length function on $W$. For this purpose, we introduce two relations,
following Geck--Pfeiffer \cite{GePf93}.

Given $x,y\in W$ and $s\in S$, we write $x\stackrel{s}{\longrightarrow} y$
if $y=sxs$ and $l(y)\leq l(x)$. We shall write $x\longrightarrow y$ if
there are sequences $x_0,x_1,\ldots,x_n\in W$ and $s_1, \ldots,s_n\in S$ 
(for some $n\geq 0$) such that 
\[x=x_0\stackrel{s_1}{\longrightarrow} x_1 \stackrel{s_2}{\longrightarrow} 
x_2 \stackrel{s_3}{\longrightarrow} \quad \cdots \quad 
\stackrel{s_{n}}{\longrightarrow} x_n=y.\]
Thus, we have $x\longrightarrow y$ if we can go from $x$ to $y$ by
a chain of conjugation with generators in $S$ such that, at each step,
the length of the elements either remains the same or decreases.

In a slightly different direction, let us now consider two elements
$x,y\in W$ such that $l(x)=l(y)$. We write $x\stackrel{w}{\sim} y$ (where
$w\in W$) if $wx=yw$ and $l(wx)=l(w)+l(x)$ or $xw=wy$ and $l(wy)=l(w)+l(y)$. 
We write $x\sim y$ if there are sequences $x_0,x_1,\ldots,x_n\in W$ and 
$w_1, \ldots,w_n\in W$ (for some $n\geq 0$) such that
\[x=x_0\stackrel{w_1}{\sim} x_1 \stackrel{w_2}{\sim} x_2 
\stackrel{w_3}{\sim} \quad \cdots \quad \stackrel{w_n}{\sim} x_n=y.\]
Thus, we have $x\sim y$ if we can go from $x$ to $y$ by a chain of conjugation
with elements of $W$ such that, at each step, the length of the
elements remains the same and an additional length condition involving
the conjugating elements is satisfied. This additional condition has
the following significance. Consider a group $\cM$ with multiplication 
$\star$ and assume that we have a map $f\colon S \rightarrow \cM$ which 
satisfies the requirements in Matsumoto's theorem~\ref{Mmatsum}. Then we have 
a canonical extension of $f$ to a map $F\colon W \rightarrow \cM$ such that
$F(ww')=F(w)\star F(w')$ whenever $l(ww')=l(w)+l(w')$. Hence, in this setting,
we have 
\[x\sim y\quad\Rightarrow\quad F(x),F(y)\mbox{ are conjugate in $\cM$}.\]
Thus, we can think of the relation ``$\sim$'' as ``universal conjugacy''.
\end{para}

\begin{thm}[Geck--Pfeiffer \cite{GePf93}, \cite{GePf}] \label{Mgp} Let 
$(W,S)$ be a finite Coxeter system and $C$ be a conjugacy class in $W$. 
We set $l_{\operatorname{min}}(C):= \min\{l(w)\mid w\in C\}$ and
\[ C_{\operatorname{min}}:=\{w\in C\mid l(w)=l_{\operatorname{min}}(C)\}.\]
Then the following hold:
\begin{itemize}
\item[(a)] For every $x\in C$, there exists some
$y\in C_{\operatorname{min}}$ such that $x\longrightarrow y$.
\item[(b)] For any two elements $x,y\in C_{\operatorname{min}}$, we 
have $x\sim y$.
\end{itemize}
\end{thm}

Precursors of the above result for type $A$ have been found much earlier
by Starkey \cite{Sta}; see also Ram \cite{Ram}. The above result allows
to define the {\em character table} of the Iwahori--Hecke algebra 
associated with $(W,S)$. This is discussed in more detail in the chapter
on Hecke algebras (??). See Richardson \cite{Richardson82},  Geck--Michel
\cite{GeMi}, Geck--Pfeiffer \cite[Chap.~3]{GePf}, Geck--Kim--Pfeiffer 
\cite{GeKiPf} and Shi \cite{Shi1} for further results on conjugacy classes. 
Krammer \cite{Kram} studies the conjugacy problem for arbitrary (infinite) 
Coxeter groups.

\begin{para}{The crystallographic condition}\label{Mse24} Let 
$M=(m_{st})_{s,t\in S}$ be a Coxeter matrix such that the connected components
of the corresponding Coxeter graph occur in Table~\ref{Mtab1}. We say that 
$M$ satisfies the {\em crystallographic condition} if there exists a Cartan 
matrix $C$ associated with $M$ which has integral coefficients. In this case, 
the corresponding reflection group $W=W(C)$ is called a {\em Weyl group}. 
The significance of this notion is that there exists a corresponding 
semisimple Lie algebra over $\CC$; see (\ref{Mlie}).

Now assume that $M$ is crystallographic. By condition (C2) this 
implies $m_{st}\in \{2,3,4,6\}$ for all $s\neq t$. Conversely, if $m_{st}$ 
satisfies this condition, then we have 
\begin{align*}
c_{st}c_{ts} & =0  \qquad \mbox{if $m_{st}=2$},\\
c_{st}c_{ts} & =1  \qquad \mbox{if $m_{st}=3$},\\
c_{st}c_{ts} & =2  \qquad \mbox{if $m_{st}=4$},\\
c_{st}c_{ts} & =3  \qquad \mbox{if $m_{st}=6$}.
\end{align*}
Thus, in each of these cases, we see that there are only two choices for
$c_{st}$ and $c_{ts}$: we must have $c_{st}=-1$ or $c_{ts}=-1$ (and then
the other value is determined). We encode this additional information in 
the Coxeter graph, by putting an arrow on the edge between the nodes 
labelled by $s,t$ according to the following scheme:
\begin{center}
$\begin{array}{ccll} \hline & \mbox{Edge between $s \neq t$} &
\multicolumn{2}{l}{\mbox{Values for $c_{st}$, $c_{ts}$}} \\  \hline 
\mbox{$m_{st}=3$} & \mbox{no arrow} & \multicolumn{2}{l}{c_{st}=c_{ts}=
-1}\\ m_{st}= 4 & 
\begin{picture}(32,15) 
\put(5,03){\circle*{5}}
\put(5,08){$\scriptstyle{s}$}
\put(5,05){\line(1,0){25}}
\put(5,1){\line(1,0){25}}
\put(13,.5){$>$}
\put(30,03){\circle*{5}}
\put(30,08){$\scriptstyle{t}$}
\end{picture} & c_{st}=-1& c_{ts}=-2 \\  m_{st}=6 & 
\begin{picture}(32,14) 
\put(5,03){\circle*{6}}
\put(5,09){$\scriptstyle{s}$}
\put(5,06){\line(1,0){25}}
\put(5,03){\line(1,0){25}}
\put(5,0){\line(1,0){25}}
\put(13,.5){$>$}
\put(30,03){\circle*{6}}
\put(30,09){$\scriptstyle{t}$}
\end{picture} & c_{st}=-1& c_{ts}=-3 \\ 
\hline \end{array}$
\end{center}
The Coxeter graph of $M$ equipped with this additional information will
be called a {\em Dynkin diagram}; it uniquely determines an integral Cartan 
matrix $C$ associated with $M$ (if such a Cartan matrix exists). The 
complete list of connected components of Dynkin diagrams is given in 
Table~\ref{Mtab2}. We see that all irreducible finite Coxeter groups are 
Weyl groups, except for those of type $H_3$, $H_4$ and $I_2(m)$ where 
$m=5$ or $m\geq 7$.  Note that type $B_n$ is the only case where we
have two different Dynkin diagrams associated with the same Coxeter graph.
\end{para}

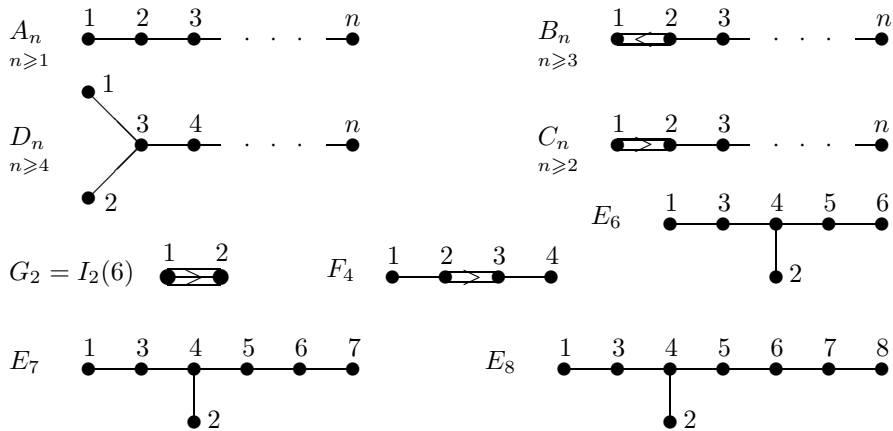
\begin{table}[htbp] \caption{Dynkin diagrams of Cartan
matrices of finite type} 
\label{Mtab2} \begin{center} \makeatletter
\vbox{
\begin{picture}(345,170)
\put( 10, 25){$E_7$}
\put( 40, 25){\@dot{5}}
\put( 38, 30){$1$}
\put( 40, 25){\line(1,0){20}}
\put( 60, 25){\@dot{5}}
\put( 58, 30){$3$}
\put( 60, 25){\line(1,0){20}}
\put( 80, 25){\@dot{5}}
\put( 78, 30){$4$}
\put( 80, 25){\line(0,-1){20}}
\put( 80, 05){\@dot{5}}
\put( 85, 03){$2$}
\put( 80, 25){\line(1,0){20}}
\put(100, 25){\@dot{5}}
\put( 98, 30){$5$}
\put(100, 25){\line(1,0){20}}
\put(120, 25){\@dot{5}}
\put(118, 30){$6$}
\put(120, 25){\line(1,0){20}}
\put(140, 25){\@dot{5}}
\put(138, 30){$7$}

\put(190, 25){$E_8$}
\put(220, 25){\@dot{5}}
\put(218, 30){$1$}
\put(220, 25){\line(1,0){20}}
\put(240, 25){\@dot{5}}
\put(238, 30){$3$}
\put(240, 25){\line(1,0){20}}
\put(260, 25){\@dot{5}}
\put(258, 30){$4$}
\put(260, 25){\line(0,-1){20}}
\put(260, 05){\@dot{5}}
\put(265, 03){$2$}
\put(260, 25){\line(1,0){20}}
\put(280, 25){\@dot{5}}
\put(278, 30){$5$}
\put(280, 25){\line(1,0){20}}
\put(300, 25){\@dot{5}}
\put(298, 30){$6$}
\put(300, 25){\line(1,0){20}}
\put(320, 25){\@dot{5}}
\put(318, 30){$7$}
\put(320, 25){\line(1,0){20}}
\put(340, 25){\@dot{5}}
\put(338, 30){$8$}

\put( 10, 59){$G_2=I_2(6)$}
\put( 70, 60){\@dot{6}}
\put( 68, 66){$1$}
\put( 70, 57){\line(1,0){20}}
\put( 70, 60){\line(1,0){20}}
\put( 70, 63){\line(1,0){20}}
\put( 76, 57.5){$>$}
\put( 90, 60){\@dot{6}}
\put( 88, 66){$2$}

\put(130, 60){$F_4$}
\put(155, 60){\@dot{5}}
\put(153, 65){$1$}
\put(155, 60){\line(1,0){20}}
\put(175, 60){\@dot{5}}
\put(173, 65){$2$}
\put(175, 58){\line(1,0){20}}
\put(175, 62){\line(1,0){20}}
\put(181, 57.5){$>$}
\put(195, 60){\@dot{5}}
\put(193, 65){$3$}
\put(195, 60){\line(1,0){20}}
\put(215, 60){\@dot{5}}
\put(213, 65){$4$}

\put(230, 80){$E_6$}
\put(260, 80){\@dot{5}}
\put(258, 85){$1$}
\put(260, 80){\line(1,0){20}}
\put(280, 80){\@dot{5}}
\put(278, 85){$3$}
\put(280, 80){\line(1,0){20}}
\put(300, 80){\@dot{5}}
\put(298, 85){$4$}
\put(300, 80){\line(0,-1){20}}
\put(300, 60){\@dot{5}}
\put(305, 58){$2$}
\put(300, 80){\line(1,0){20}}
\put(320, 80){\@dot{5}}
\put(318, 85){$5$}
\put(320, 80){\line(1,0){20}}
\put(340, 80){\@dot{5}}
\put(338, 85){$6$}

\put( 10,110){$D_n$}
\put( 10,100){$\scriptstyle{n \geq 4}$}
\put( 40,130){\@dot{5}}
\put( 45,130){$1$}
\put( 40,130){\line(1,-1){21}}
\put( 40, 90){\@dot{5}}
\put( 46, 85){$2$}
\put( 40, 90){\line(1,1){21}}
\put( 60,110){\@dot{5}}
\put( 58,115){$3$}
\put( 60,110){\line(1,0){30}}
\put( 80,110){\@dot{5}}
\put( 78,115){$4$}
\put(100,110){\@dot{1}}
\put(110,110){\@dot{1}}
\put(120,110){\@dot{1}}
\put(130,110){\line(1,0){10}}
\put(140,110){\@dot{5}}
\put(137,115){$n$}

\put(210,110){$C_n$}
\put(210,100){$\scriptstyle{n \geq 2}$}
\put(240,110){\@dot{5}}
\put(238,115){$1$}
\put(240,108){\line(1,0){20}}
\put(240,112){\line(1,0){20}}
\put(246,107.5){$>$}
\put(260,110){\@dot{5}}
\put(258,115){$2$}
\put(260,110){\line(1,0){30}}
\put(280,110){\@dot{5}}
\put(278,115){$3$}
\put(300,110){\@dot{1}}
\put(310,110){\@dot{1}}
\put(320,110){\@dot{1}}
\put(330,110){\line(1,0){10}}
\put(340,110){\@dot{5}}
\put(337,115){$n$}

\put( 10,150){$A_n$}
\put( 10,140){$\scriptstyle{n \geq 1}$}
\put( 40,150){\@dot{5}}
\put( 38,155){$1$}
\put( 40,150){\line(1,0){20}}
\put( 60,150){\@dot{5}}
\put( 58,155){$2$}
\put( 60,150){\line(1,0){30}}
\put( 80,150){\@dot{5}}
\put( 78,155){$3$}
\put(100,150){\@dot{1}}
\put(110,150){\@dot{1}}
\put(120,150){\@dot{1}}
\put(130,150){\line(1,0){10}}
\put(140,150){\@dot{5}}
\put(137,155){$n$}

\put(210,150){$B_n$}
\put(210,140){$\scriptstyle{n \geq 3}$}
\put(240,150){\@dot{5}}
\put(238,155){$1$}
\put(240,148){\line(1,0){20}}
\put(240,152){\line(1,0){20}}
\put(246,147.5){$<$}
\put(260,150){\@dot{5}}
\put(258,155){$2$}
\put(260,150){\line(1,0){30}}
\put(280,150){\@dot{5}}
\put(278,155){$3$}
\put(300,150){\@dot{1}}
\put(310,150){\@dot{1}}
\put(320,150){\@dot{1}}
\put(330,150){\line(1,0){10}}
\put(340,150){\@dot{5}}
\put(338,155){$n$}
\end{picture}}
\makeatother
\end{center}
\end{table}

The following discussion of root systems associated with finite 
reflection groups follows the appendix on finite reflection groups
 in Steinberg \cite{Steinberg0}.

\begin{para}{Root systems}\label{Msec24} Let $V$ be a finite dimensional
real vector space and let $(\;,\;)$ be a positive-definite scalar product
on $V$. Given a non-zero vector $\alpha\in V$, the corresponding 
reflection $w_\alpha\in \GL(V)$ is defined by
\[ w_\alpha(v)=v-2\frac{(v,\alpha)}{(\alpha,\alpha)}\alpha \qquad
\mbox{for all $v\in V$}.\]
Note that, for any $w\in \GL(V)$, we have $ww_\alpha w^{-1}=w_{w(\alpha)}$.
A finite subset $\phi \subseteq V\setminus\{0\}$ is called a {\em root system}
if the following conditions are satisfied:
\begin{itemize}
\item[(R1)] For any $\alpha\in \Phi$, we have $\Phi\cap {\RR}\alpha=\{\pm 
\alpha\}$;
\item[(R2)] For every $\alpha,\beta\in \Phi$, we have $w_\alpha(\beta)\in 
\Phi$.
\end{itemize}
Let $W(\Phi)\subseteq \GL(V)$ be the subgroup generated by the reflections
$w_\alpha$ ($\alpha \in \Phi$). A subset $\Pi \subseteq \Phi$ is called
a {\em simple system} if $\Pi$ is linearly independent and if any root 
in $\Phi$ can be written as a linear combination of the elements of $\Pi$ 
in which all non-zero coefficients are either all positive or all negative. 
It is known that simple systems always exist, and that any two simple systems
can be transformed into each other by an element of $W(\Phi)$. Let 
now fix a simple system $\Pi\subseteq \Phi$. Then it is also known that
\[W(\Phi)=\langle w_\alpha \mid \alpha \in \Pi\rangle \qquad \mbox{and}
\qquad (\alpha,\beta)\geq 0\quad\mbox{for all $\alpha\neq \beta$ in $\Pi$}.\]
For $\alpha \neq \beta$ in $\Pi$, let $m_{\alpha\beta}\geq 2$ be the order 
of $w_\alpha w_\beta$ in $\GL(V)$. Then $w_\alpha w_\beta$
is a rotation in $V$ through the angle $2\pi/m_{\alpha\beta}$ and that we 
have the relation
\[\frac{(\alpha,\beta)}{(\alpha,\alpha)}\frac{(\beta,\alpha)}{(\beta,\beta)}
=\cos^2(\pi/m_{\alpha\beta})\quad \mbox{for all $\alpha \neq \beta$ in 
$\Pi$}.\]
Thus, if we set $M:=(m_{\alpha\beta})_{\alpha,\beta\in \Pi}$ (where
$m_{\alpha\alpha}=1$ for all $\alpha\in \Pi$) and define 
\[ C:=(c_{\alpha\beta})_{\alpha,\beta\in \Pi}\quad \mbox{where} \quad 
c_{\alpha\beta}:= 2\frac{(\alpha,\beta)}{(\alpha,\alpha)} 
\quad \mbox{for $\alpha,\beta\in S$},\]
then $M$ is a Coxeter matrix, $C$ is an associated Cartan matrix and
$W(\Phi)=W(C)$ is the Coxeter group with Coxeter matrix $M$; see 
Theorem~\ref{Mthm1}. Thus, every root system leads to a finite Coxeter
group. 

Conversely, let $(W,S)$ be a Coxeter system associated to a Coxeter 
matrix $M=(m_{st})_{s,t\in S}$. Let $C$ be a corresponding Cartan matrix 
such that $c_{st}=c_{ts}=-2\cos(\pi/m_{st})$ for every $s\neq t$ in $S$ such 
that $m_{st}$ is odd. Then one can show that 
\[ \Phi=\Phi(C):=\{ w(\alpha_s) \mid w \in W, s \in S\} \subseteq V\]
is a root system with simple system $\{\alpha_s\mid s\in S\}$ and
that $W(\Phi)=W$; see, e.g., Geck--Pfeiffer \cite[Chap.~1]{GePf}. Thus, every 
Coxeter group leads to a root system.
\end{para}

Given a root system $\Phi$ as above, there is a strong link between the 
combinatorics of the Coxeter presentation of $W=W(\Phi)$ and the geometry 
of $\Phi$. To state the following basic result, we set 
\[ \Phi^+:=\{\alpha \in \Pi \mid \alpha=\sum_{\beta \in \Pi}
  x_\beta \beta \mbox{ where $x_\beta \in {\RR}_{\geq 0}$ for all 
  $\beta \in \Pi$}\};\]
the roots in $\Phi^+$ will be called {\em positive roots}. Similarly,
$\Phi^-:=-\Phi^+$ will be called the set of  {\em negative roots}.  By
the definition of a simple system, we have $\Phi=\Phi^+\amalg \Phi^-$. 

\begin{prop} \label{Mlengthroots} Given $\alpha \in \Pi$ and $w \in W$, we have
\begin{align*} 
w^{-1}(\alpha) \in \Phi^+ &\quad\Leftrightarrow\quad l(w_\alpha w)=l(w)+1,\\
w^{-1}(\alpha) \in \Phi^- &\quad\Leftrightarrow\quad l(w_\alpha w)=l(w)-1.
\end{align*}
Furthermore, for any $w\in W$, we have $l(w)=|\{\alpha \in \Phi^+\mid
w(\alpha)\in \Phi^-\}|$.
\end{prop}

The root systems associated with finite Weyl groups are explicitly 
described in Bourbaki \cite[p.~251--276]{Bou}. For type $H_3$ and $H_4$, 
see Humphreys \cite[2.13]{Humphreys2}. \par
Bremke--Malle \cite{BreMa1,BreMa2} have studied suitable generalisations of
root systems and length functions for the infinite series $G(d,1,n)$ and
$G(e,e,n)$, which have subsequently been extended in weaker form by
Rampetas--Shoji \cite{RaSh} to arbitrary imprimitive reflection groups. For
investigations of root systems see also Nebe \cite{Neb} and Hughes--Morris
\cite{HuMo}. 
But there is no general theory of root systems and length functions for
complex reflection groups (yet).

\begin{para}{Torsion primes} \label{Mtor} Assume that $\Phi \subseteq
V$ is a root system as above, with a set of simple roots $\Pi\subseteq 
\Phi$. Assume that the corresponding Cartan matrix $C$ is indecomposable
and has integral coefficients. Thus, its Dynkin diagram is one of the
graphs in Table~\ref{Mtab2}. Following Springer--Steinberg 
\cite[\S I.4]{SprSt} we shall now discuss ``bad primes'' and ``torsion
primes'' with respect to $\Phi$.

For every $\alpha\in \Phi$, the corresponding coroot is defined by 
$\alpha^*:=2\alpha/(\alpha,\alpha)$. Then $\Phi^*:=\{\alpha^*\mid 
\alpha\in \Phi\}$ also is a root system, the dual of $\Phi$. The 
Dynkin diagram of $\Phi^*$ is obtained from that of $\Phi$ by reversing
the arrows. (For example, the dual of a root
system of type $B_n$ is of type $C_n$.) 

Let $L(\Phi)$ denote the lattice spanned by $\Phi$ in $V$. A prime number
$p>0$ is called {\em bad} for $\Phi$ if $L(\Phi)/L(\Phi_1)$ has $p$-torsion
for some (integrally) closed subsystem $\Phi_1$ of $\Phi$. The prime
$p$ is called a {\em torsion prime} if $L(\Phi^*)/L(\Phi_1^*)$ has
$p$-torsion for some closed subsystem $\Phi_1$ of $\Phi$. Note that $\Phi_1^*$
need not be closed in $\Phi^*$, and so the torsion primes for $\Phi$ 
and the bad primes for $\Phi^*$ need not be the same. The bad primes
can be characterised as follows. Let $\alpha_0=\sum_{\alpha\in \Pi} 
m_\alpha \alpha$ be the unique positive root of maximal height.
(The {\em height} of a root is the sum of the coefficients in the expression
of that root as a linear combination of simple roots.) Then we have:
\[\mbox{$p$ bad} 
\quad \Leftrightarrow \quad \begin{array}{c} p=m_\alpha \\ 
\mbox{ for some $\alpha$}\end{array} 
\quad \Leftrightarrow \quad \begin{array}{c} p \mbox{ divides }m_\alpha \\ 
\mbox{ for some $\alpha$}\end{array} 
\quad \Leftrightarrow \quad \begin{array}{c} p\leq m_\alpha \\ 
\mbox{ for some $\alpha$}\end{array}.\]
Now let $\alpha_0^*=\sum_{\alpha \in \Pi} m_\alpha^* \alpha^*$. Then 
$p$ is a torsion prime if and only if $p$ satisfies one of the above
conditions, with $m_\alpha$ replaced by $m_\alpha^*$. For the 
various roots systems, the bad primes and the torsion primes are given
as follows.
\[\begin{array}{l|ccccccccc} \hline \mbox{Type} & A_n & B_n & C_n & 
D_n &G_2 & F_4 & E_6 & E_7 &E_8\\
& & (n\geq 2) & (n \geq 2) & (n\geq 4) & & & & & \\ \hline
\mbox{Bad} & \mbox{none} & 2 & 2 & 2 & 2,3 & 2,3 & 2,3 & 2,3 & 2,3,5\\
\hline
\mbox{Torsion} & \mbox{none} & 2 & \mbox{none} & 2 & 2 & 2,3 & 2,3 & 2,3
& 2,3,5 \\ \hline \end{array}\]
The bad primes and torsion primes play a role in various questions
related to sub-root systems, centralizers of semisimple elements
in algebraic groups, the classification of unipotent classes in
simple algebraic groups and so on; see \cite{SprSt} and also the survey
in \cite[\S\S 1.14--1.15]{Carter2}. 
\end{para}

\begin{para}{Affine Weyl groups} \label{affgrp} Let $\Phi
\subseteq V$ be a root system as above, with Weyl group $W$.
Let $L(\Phi):=\sum_{\alpha \in \Pi} {\ZZ}\alpha \subseteq V$ be
the lattice spanned by the roots in $V$. Then $W$ leaves $L(\Phi)$
invariant and we have a natural group homomorphism 
$W\rightarrow \mbox{Aut}(L(\Phi))$. The semidirect product
\[ W_a(\Phi):=L(\Phi) \rtimes W\]
is called the affine Weyl group associated with the root system $W$;
see Bourbaki \cite[Chap.~VI, \S 2]{Bou}. The group $W_a(\Phi)$ itself 
is a Coxeter group. The corresponding presentation can also be encoded 
in a graph, as follows. Let  $\alpha_0$ be the unique positive root
of maximal height in $\Phi$. We define an 
extended Cartan matrix $\tilde{C}$ by similar rules as before:
\[ \tilde{c}_{\alpha\beta}:=2\frac{(\alpha,\beta)}{(\alpha,\alpha)}
\qquad \mbox{for $\alpha,\beta \in \Pi \cup \{-\alpha_0\}$}.\]
The extended Dynkin diagrams encoding these matrices for irreducible $W$ are
given in  Table~\ref{Mtab2a}.  They are obtained from the diagrams in 
Table~\ref{Mtab2} by adjoining an additional node (corresponding to 
$-\alpha_0$) and putting edges according to the same rules as before.
\end{para}

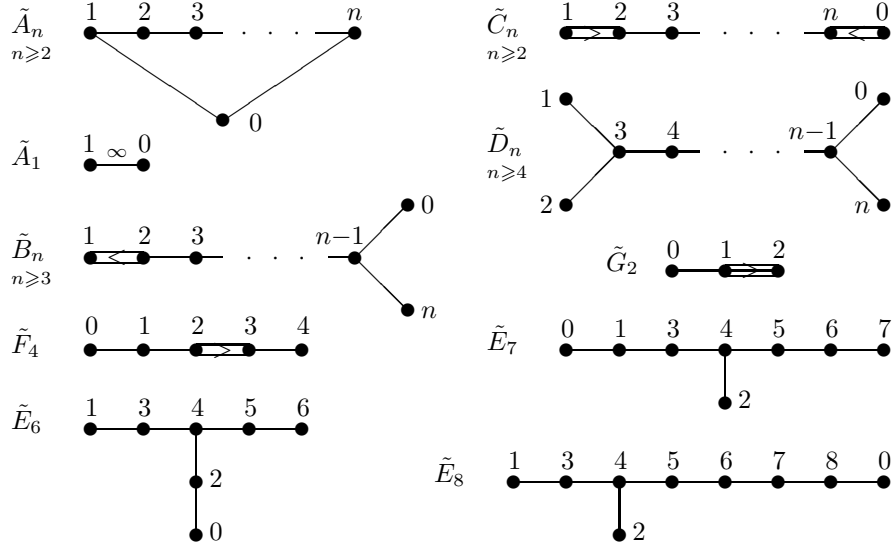
\begin{table}[htbp] \caption{Extended Dynkin diagrams} 
\label{Mtab2a} \begin{center} \makeatletter
\vbox{
\begin{picture}(345,205)
\put( 10, 45){$\tilde{E}_6$}
\put( 40, 45){\@dot{5}}
\put( 38, 50){$1$}
\put( 40, 45){\line(1,0){20}}
\put( 60, 45){\@dot{5}}
\put( 58, 50){$3$}
\put( 60, 45){\line(1,0){20}}
\put( 80, 45){\@dot{5}}
\put( 78, 50){$4$}
\put( 80, 45){\line(0,-1){20}}
\put( 80, 25){\@dot{5}}
\put( 85, 23){$2$}
\put( 80, 25){\line(0,-1){20}}
\put( 80,  5){\@dot{5}}
\put( 85,  3){$0$}
\put( 80, 45){\line(1,0){20}}
\put(100, 45){\@dot{5}}
\put( 98, 50){$5$}
\put(100, 45){\line(1,0){20}}
\put(120, 45){\@dot{5}}
\put(118, 50){$6$}

\put(170, 25){$\tilde{E}_8$}
\put(200, 25){\@dot{5}}
\put(198, 30){$1$}
\put(200, 25){\line(1,0){20}}
\put(220, 25){\@dot{5}}
\put(218, 30){$3$}
\put(220, 25){\line(1,0){20}}
\put(240, 25){\@dot{5}}
\put(238, 30){$4$}
\put(240, 25){\line(0,-1){20}}
\put(240, 05){\@dot{5}}
\put(245, 03){$2$}
\put(240, 25){\line(1,0){20}}
\put(260, 25){\@dot{5}}
\put(258, 30){$5$}
\put(260, 25){\line(1,0){20}}
\put(280, 25){\@dot{5}}
\put(278, 30){$6$}
\put(280, 25){\line(1,0){20}}
\put(300, 25){\@dot{5}}
\put(298, 30){$7$}
\put(300, 25){\line(1,0){20}}
\put(320, 25){\@dot{5}}
\put(318, 30){$8$}
\put(320, 25){\line(1,0){20}}
\put(340, 25){\@dot{5}}
\put(338, 30){$0$}

\put( 10, 74){$\tilde{F}_4$}
\put( 40, 75){\@dot{5}}
\put( 38, 81){$0$}
\put( 40, 75){\line(1,0){20}}
\put( 60, 75){\@dot{5}}
\put( 58, 81){$1$}
\put( 60, 75){\line(1,0){20}}
\put( 80, 75){\@dot{5}}
\put( 78, 81){$2$}
\put( 80, 73){\line(1,0){20}}
\put( 80, 77){\line(1,0){20}}
\put( 86, 72.5){$>$}
\put(100, 75){\@dot{5}}
\put( 98, 81){$3$}
\put(100, 75){\line(1,0){20}}
\put(120, 75){\@dot{5}}
\put(118, 81){$4$}

\put(190, 75){$\tilde{E}_7$}
\put(220, 75){\@dot{5}}
\put(218, 80){$0$}
\put(220, 75){\line(1,0){20}}
\put(240, 75){\@dot{5}}
\put(238, 80){$1$}
\put(240, 75){\line(1,0){20}}
\put(260, 75){\@dot{5}}
\put(258, 80){$3$}
\put(280, 75){\line(0,-1){20}}
\put(280, 55){\@dot{5}}
\put(285, 53){$2$}
\put(260, 75){\line(1,0){20}}
\put(280, 75){\@dot{5}}
\put(278, 80){$4$}
\put(280, 75){\line(1,0){20}}
\put(300, 75){\@dot{5}}
\put(298, 80){$5$}
\put(300, 75){\line(1,0){20}}
\put(320, 75){\@dot{5}}
\put(318, 80){$6$}
\put(320, 75){\line(1,0){20}}
\put(340, 75){\@dot{5}}
\put(338, 80){$7$}

\put(235,105){$\tilde{G}_2$}
\put(260,105){\@dot{5}}
\put(258,110){$0$}
\put(260,105){\line(1,0){20}}
\put(280,105){\@dot{5}}
\put(278,110){$1$}
\put(280,103){\line(1,0){20}}
\put(280,105){\line(1,0){20}}
\put(280,107){\line(1,0){20}}
\put(286,102.5){$>$}
\put(300,105){\@dot{5}}
\put(298,110){$2$}

\put( 10,110){$\tilde{B}_n$}
\put( 10,100){$\scriptstyle{n \geq 3}$}
\put( 40,110){\@dot{5}}
\put( 37,115){$1$}
\put( 46,107.5){$<$}
\put( 40,112){\line(1,0){20}}
\put( 40,108){\line(1,0){20}}
\put( 60,110){\@dot{5}}
\put( 58,115){$2$}
\put( 60,110){\line(1,0){30}}
\put( 80,110){\@dot{5}}
\put( 78,115){$3$}
\put(100,110){\@dot{1}}
\put(110,110){\@dot{1}}
\put(120,110){\@dot{1}}
\put(130,110){\line(1,0){10}}
\put(140,110){\@dot{5}}
\put(125,115){$n{-}1$}
\put(140,110){\line(1,1){21}}
\put(140,110){\line(1,-1){21}}
\put(160,130){\@dot{5}}
\put(160, 90){\@dot{5}}
\put(165, 87){$n$}
\put(165,127){$0$}

\put(190,150){$\tilde{D}_n$}
\put(190,140){$\scriptstyle{n \geq 4}$}
\put(220,170){\@dot{5}}
\put(210,167){$1$}
\put(220,170){\line(1,-1){21}}
\put(220,130){\@dot{5}}
\put(210,127){$2$}
\put(220,130){\line(1,1){21}}
\put(240,150){\@dot{5}}
\put(238,155){$3$}
\put(240,150){\line(1,0){30}}
\put(260,150){\@dot{5}}
\put(258,155){$4$}
\put(280,150){\@dot{1}}
\put(290,150){\@dot{1}}
\put(300,150){\@dot{1}}
\put(310,150){\line(1,0){10}}
\put(320,150){\@dot{5}}
\put(304,155){$n{-}1$}
\put(320,150){\line(1,1){21}}
\put(320,150){\line(1,-1){21}}
\put(340,170){\@dot{5}}
\put(340,130){\@dot{5}}
\put(329,127){$n$}
\put(329,170){$0$}

\put(190,195){$\tilde{C}_n$}
\put(190,185){$\scriptstyle{n \geq 2}$}
\put(220,195){\@dot{5}}
\put(218,200){$1$}
\put(220,193){\line(1,0){20}}
\put(220,197){\line(1,0){20}}
\put(226,192.5){$>$}
\put(240,195){\@dot{5}}
\put(238,200){$2$}
\put(240,195){\line(1,0){30}}
\put(260,195){\@dot{5}}
\put(258,200){$3$}
\put(280,195){\@dot{1}}
\put(290,195){\@dot{1}}
\put(300,195){\@dot{1}}
\put(310,195){\line(1,0){10}}
\put(320,195){\@dot{5}}
\put(317,200){$n$}
\put(320,193){\line(1,0){20}}
\put(320,197){\line(1,0){20}}
\put(326,192.5){$<$}
\put(340,195){\@dot{5}}
\put(337,200){$0$}

\put( 10,195){$\tilde{A}_n$}
\put( 10,185){$\scriptstyle{n \geq 2}$}
\put( 40,195){\@dot{5}}
\put( 38,200){$1$}
\put( 40,195){\line(1,0){20}}
\put( 60,195){\@dot{5}}
\put( 58,200){$2$}
\put( 60,195){\line(1,0){30}}
\put( 80,195){\@dot{5}}
\put( 78,200){$3$}
\put( 98,195){\@dot{1}}
\put(108,195){\@dot{1}}
\put(118,195){\@dot{1}}
\put(125,195){\line(1,0){15}}
\put(140,195){\@dot{5}}
\put(137,200){$n$}
\put( 40,195){\line(3,-2){50}}
\put( 90,162){\@dot{5}}
\put(140,195){\line(-3,-2){50}}
\put(100,158){$0$}

\put(10,145){$\tilde{A}_1$}
\put(40,145){\@dot{5}}
\put(37,150){$1$}
\put(40,145){\line(1,0){20}}
\put(46,148){$\scriptstyle{\infty}$}
\put(60,145){\@dot{5}}
\put(58,150){$0$}
\end{picture}}
\makeatother
\end{center}
\end{table}

In the following subsections, we describe some situations where Coxeter 
groups and root systems arise ``in nature''.

\begin{para}{Kac--Moody algebras} \label{Mlie} Here we briefly discuss how 
Coxeter groups and root systems arise in the theory of Lie algebras or, more 
generally, Kac--Moody algebras. We follow the exposition in Kac \cite{Kac}.
Let $C=(c_{st})_{s,t\in S}$ be a Cartan matrix all of whose coefficients 
are integers. We also assume that $C$ is symmetrisable, i.e., there 
exists a diagonal invertible matrix $D$ and a symmetric matrix $B$ such 
that $C=DB$. A realisation of $C$ is a triple $(\fh,\Pi, \Pi^\vee)$ where 
$\fh$ is a complex vector space, $\Pi=\{\alpha_s\mid s\in S\}\subseteq
\fh^*:= \mbox{Hom}(\fh,\CC)$ and $\Pi^\vee=\{\alpha_s^\vee\mid s\in S\}$ are 
subsets of $\fh^*$ and $\fh$, respectively, such that the following conditions
hold. 
\begin{itemize}
\item[(a)] Both sets $\Pi$ and $\Pi^\vee$ are linearly independent;
\item[(b)] we have $\langle \alpha_s^\vee,\alpha_t\rangle:=
\alpha_t(\alpha_s^\vee)=c_{st}$ for all $s,t\in S$;
\item[(c)] $|S|-\mbox{rank}(C)=\dim \fh-|S|$.
\end{itemize}
Let $\fg(C)$ be the corresponding Kac--Moody algebra. Then $\fg(C)$ is a 
Lie algebra which is generated by $\fh$ together with two collections of 
elements $\{e_s\mid s\in S\}$ and $\{f_s\mid s\in S\}$, where the following
relations hold:
\begin{alignat*}{2}
[e_s,f_t] &= \delta_{st}\alpha_s^\vee &&\qquad (s,t\in S),\\
[h,h'] &= 0 &&\qquad (h,h'\in \fh),\\
[h,e_s] &= \langle h,\alpha_s,\rangle e_s &&\qquad (s\in S, h\in \fh),\\
[h,f_s] &= -\langle h,\alpha_s\rangle f_s &&\qquad (s\in S, h\in \fh),\\
(\mbox{ad}\ e_s&)^{1-c_{st}}e_t = 0 &&\qquad (s,t\in S, s\neq t),\\
(\mbox{ad}\ f_s&)^{1-c_{st}}f_t = 0 &&\qquad (s,t\in S, s\neq t).
\end{alignat*}
(By \cite[9.11]{Kac}, this is a set of defining relations for $\fg(C)$.) We 
have a direct sum decomposition 
\[\fg(C)=\fh\oplus \bigoplus_{0\neq \alpha \in Q} \fg_\alpha(C)\qquad 
\mbox{where}\qquad Q:=\sum_{s\in S}  {\ZZ}\alpha_s \subseteq \fh^*\]
and $\fg_\alpha(C):=\{x\in \fg(C) \mid [h,x]=\alpha(h)x \mbox{ for 
all $h\in\fh$}\}$ for all $\alpha\in Q$; here, $\fh=\fg_0$. The
set of all $0\neq \alpha\in Q$ such that $\fg_\alpha(C)\neq \{0\}$ will
be denoted by $\Phi$ and called the {\em root system} of $\fg(C)$. 

For each $s\in S$, we define a linear map $\sigma_s\colon \fh^*\rightarrow
\fh^*$ by the formula
\[ \sigma_s(\lambda)=\lambda-\langle \lambda,\alpha_s^\vee\rangle\,
\alpha_s \qquad \mbox{for $\lambda \in \fh^*$}.\]
Then it is easily checked that $\sigma_s$ is a reflection where 
$\sigma_s(\alpha_s)=-\alpha_s$. We set  
\[W=W(C)=\langle \sigma_s \mid s\in S\rangle \subseteq \GL(\fh^*).\] 
Now we can state (see \cite[3.7, 3.11 and 3.13]{Kac}):
\begin{itemize}
\item[(a)] The pair $(W,\{\sigma_s\mid s\in S\})$ is a Coxeter
system; the corresponding Coxeter matrix is the one associated to $C$.
\item[(b)] The root system  $\Phi$ is invariant under the action of $W$
and we have $l(\sigma_s w)=l(w)+1$ if and only if $w^{-1}(\alpha_s)\in 
\Phi^+$, where $\Phi^{+}$ is defined as in (\ref{Msec24}).
\end{itemize}
Thus, $W$ and $\Phi$ have similar properties as before. Note, however,
that here we did not make any assumption on $C$ (except that it is
symmetrisable with integer entries) and so $W$ and $\Phi$ may be infinite.
The finite case is characterised as follows:
\begin{align*}
|W|<\infty \quad &\Leftrightarrow \quad |\Phi|<\infty \quad  
\Leftrightarrow \quad \dim \fg(C)<\infty \\ \quad &\Leftrightarrow \quad 
\mbox{all connected components of $C$ occur in Table~\ref{Mtab2}}.
\end{align*}
(This follows from \cite[3.12]{Kac} and the characterisation of 
finite Coxeter groups in Theorem~\ref{Mthm1}.) In fact, the 
finite dimensional Kac--Moody algebras are precisely the ``classical''
semisimple complex Lie algebras (see, for example, Humphreys 
\cite{Humphreys0}). The Kac--Moody algebras and the root systems 
associated to so-called Cartan matrices of {\em affine type} have an 
extremely rich structure and many applications in other branches of
mathematics and mathematical physics; see Kac \cite{Kac}.
\end{para}

\begin{para}{Groups with a $BN$-pair} \label{Mbnpair}
Let $G$ be an abstract group. We say that $G$ is a {\em group with a
$BN$-pair} or that $G$ admits  a {\em Tits system} if there are subgroups 
$B,N\subseteq G$ such that the following conditions are satisfied.
\begin{itemize}
\item[(BN1)] $G$ is generated by $B$ and $N$.
\item[(BN2)] $T:=B \cap N$ is normal in $N$ and the quotient $W:=N/T$ is a
finite group generated by a set $S$ of elements of order~$2$. 
\item[(BN3)] $n_sB n_s\neq B$ if $s\in S$ and $n_s$ is a 
representative of $s$ in $N$.
\item[(BN4)] $n_sBn\subseteq Bn_snB\cup BnB$ for any $s\in S$
and $n\in N$.
\end{itemize}
The group $W$ is called the {\em Weyl group} of $G$.  In fact, it is
a consequence of the above axioms that the pair $(W,S)$ is a Coxeter
system; see \cite[Chap.~IV, \S 2, Th\'eor\`eme~2]{Bou}. The notion of
groups with a $BN$-pair has been invented by Tits; see \cite{Tits}.
The standard example of a group with a $BN$-pair is the general
linear group $G=\GL_n(K)$, where $K$ is any field and
\begin{align*}
B &:=\mbox{subgroup of all upper triangular matrices in $G$},\\
N &:=\mbox{subgroup of all monomial matrices in $G$}.
\end{align*}
(A matrix is called monomial if it has exactly one non-zero entry in each
row and and each column.) We have
\[ T:=B\cap N=\mbox{subgroup of all diagonal matrices in $G$}\]
and $W=N/T\cong \fS_n$, Thus, $\fS_n$ is the Weyl group of $G$. 
More generally, the Chevalley groups (and their twisted analogues) 
associated with the semisimple complex Lie algebras all have $BN$-pairs;
see Chevalley \cite{Chev0}, Carter \cite{Carter1} and Steinberg 
\cite{Steinberg0}.

The above set of axioms imposes very strong conditions on the structure
of a group $G$ with a $BN$-pair. For example, we have the following
{\em Bruhat decomposition}, which gives the decomposition of $G$
into double cosets with respect to~$B$:
\[ G=\coprod_{w\in W} BwB.\]
(More accurately, we should write $Bn_wB$ where $n_w$ is a 
representative of $w\in W$ in $N$. But, since any two representatives 
of $w$ lie in the same coset of $T\subseteq B$, the double coset 
$Bn_wB$ does not depend on the choice of the representative.)

Furthermore, the proof of the simplicity of the Chevalley groups  and
their twisted analogues is most economically performed using the 
simplicity criterion for abstract groups with a $BN$-pair in Bourbaki 
\cite[Chap.~IV, \S 2, no.~7]{Bou}. 

Groups with a $BN$-pair play an important r\^ole in finite group theory.
In fact, it is known that every finite simple group possesses a $BN$-pair, 
except for the cyclic groups of prime order, the alternating groups of 
degree $\geq 5$, and the $26$ sporadic simple groups; see Gorenstein et 
al. \cite{GoLySo}. Given a finite group $G$ with a $BN$-pair, the
irreducible factors of the Weyl group $W$ are of type $A_n$, $B_n$, $D_n$,
$G_2$, $F_4$, $E_6$, $E_7$, $E_8$ or $I_2(8)$. This follows from the 
classification by Tits \cite{Tits} (rank $\geq 3$), Hering, Kantor, 
Seitz \cite{HeKaSe}, \cite{KaSe} (rank~$1$) and Fong--Seitz \cite{FoSe} 
(rank~$2$). Note that there is only one case where $W$ is not 
crystallographic: this is the case where $W$ has a component of type 
$I_2(8)$ (the dihedral group of order $16$), which corresponds to the 
twisted groups of type $F_4$ discovered by Ree (see Carter \cite{Carter1}
or Steinberg \cite{Steinberg0}).

In another direction, $BN$-pairs with infinite Weyl groups arise naturally
in the theory of $p$-adic groups; see Iwahori and Matsumoto \cite{iwamat}.
\end{para}

\begin{para}{Connected reductive algebraic groups} \label{bnpair1}
Here, we assume that the reader has some familiarity with the 
theory of linear algebraic groups; see Borel \cite{Borel}, Humphreys
\cite{Humphreys1} or Springer \cite{Spr2}. Let $G$ be a connected 
reductive algebraic group over an algebraically closed field $K$. Let 
$B\subseteq G$ be a Borel subgroup. Then we have a semidirect product 
decomposition $B=U\,T$ where $U$ is the unipotent radical of $B$ and $T$ 
is a maximal torus. Let $N=\mbox{N}_G(T)$, the normaliser of $T$ in $G$. 
Then the groups $B,N$ form a $BN$-pair in $G$; furthermore, $W$ must be a 
finite Weyl group (and not just a Coxeter group as for general groups with 
a $BN$-pair). This is a deep, important result whose proof goes back to 
Chevalley \cite{Chev}; detailed expositions can be found in  the monographs
by Borel \cite[Chap.~IV, 14.15]{Borel}, Humphreys \cite[\S 29.1]{Humphreys1} 
or Springer \cite[Chap.~8]{Spr2}. 

For example, in $G=\GL_n(K)$, the subgroup $B$ of all upper triangular 
matrices is a Borel subgroup by the Lie--Kolchin Theorem (see, for example, 
Humphreys \cite[17.6]{Humphreys1}). Furthermore, we have a semidirect 
product decomposition $B=U\,T$ where $U\subseteq B$ is the normal subgroup 
consisting of all upper triangular matrices with $1$ on the diagonal and 
$T$ is the group of all diagonal matrices in $B$. Since $K$ is infinite, 
it is easily checked that $N=\mbox{N}_G(T)$, the group of all monomial
matrices.  

Returning to the general case, let us consider the Bruhat cells $BwB$ 
($w\in W$). These are locally closed subsets of $G$ since they are orbits
of $B\times B$ on $G$ under left and right multiplication. The Zariski 
closure of $BwB$ is given by 
\[ \overline{BwB}=\bigcup_{y \in \cS(w)} ByB.\]
This yields the promised geometric description of the Bruhat--Chevalley 
order $\leq$ on $W$ (as defined in the remarks following 
Theorem~\ref{Mmatsum}.) 
The proof (see, for example, Springer \cite[\S 8.5]{Spr2}) relies in 
an essential way on the fact that $G/B$ is a projective variety. 
\end{para}

\section{Braid groups} \label{sec13}

\begin{para}{The braid group of a complex reflection group}
\label{subsec:3.braid}
For a complex reflection group $W\leq\GL(V)$, $V=\CC^n$, denote by $\cA$ the
set of its reflecting hyperplanes in $V$. The topological space
$$\regV:=V\setminus\bigcup_{H\in\cA} H$$
is (pathwise) connected in its inherited complex topology. For a fixed base
point $x_0\in\regV$ we define the {\em pure braid group of $W$} as the
fundamental group $P(W):=\pi_1(\regV,x_0)$. Now $W$ acts on $\regV$, and by the
Theorem of Steinberg (Theorem~\ref{Gthm:steinberg}) the covering
$\bar{\ }:\regV\rightarrow \regV/W$ is Galois, with group $W$. This induces a
short exact sequence
\begin{equation}\label{GEq:braid}
1\longrightarrow P(W)\longrightarrow B(W)\longrightarrow W
  \longrightarrow1
\end{equation}
for the {\em braid group} $B(W):=\pi_1(\regV/W,\bar x_0)$ of $W$. \par
If $W=\fS_n$ in its natural permutation representation, the group $B(W)$ is
just the classical Artin braid group on $n$ strings \cite{Ar}. \par
We next describe some natural generators of $B(W)$. Let $H\in\cA$ be a
reflecting hyperplane. Let $x_H\in H$ and $r>0$ such that the open ball
$B(x_H,2r)$ around $x_H$ does not intersect any other reflecting hyperplane
and $x_0\notin B(x_H,2r)$. Choose a path $\gamma:[0,1]\rightarrow V$ from the
base point $x_0$ to $x_H$, with $\gamma(t)\in\regV$ for $t<1$. Let $t_0$ be
minimal subject to $\gamma(t)\in B(x_H,r)$ for all $t>t_0$. Then
$\gamma':=\gamma(t/t_0)$ is a path from $x_0$ to $\gamma(t_0)$. Then
$$\lambda:[0,1]\rightarrow B(x_H,2r),\qquad
  t\mapsto \gamma(t_0)\exp(2\pi it/e_H),$$
where $e_H=|W_H|$ is the order of the fixator of $H$ in $W$, defines a closed
path in the quotient $\regV/W$. The homotopy class in $B(W)$ of the composition
$\gamma'\circ\lambda\circ\gamma'^{-1}$ is then called a {\em braid reflection}
(see Brou\'e \cite{Br01}) or {\em generator of the monodromy around $H$}. Its
image in $W$ is a reflection $s_H$ generating $W_H$, with non-trivial
eigenvalue $\exp(2\pi i/e_H)$. It can be shown that $B(W)$ is generated by all
braid reflection, when $H$ varies over the reflecting hyperplanes of $W$
(Brou\'e--Malle--Rouquier \cite[Th.~2.17]{BMR}).
\end{para}

Assume from now on that $W$ is irreducible. Recall the definition of $N$, $N^*$
in Sections~\ref{subsec:1.Inv} and~\ref{subsec:1.fake} as the number of
reflections respectively of reflecting hyperplanes. The following can be shown
without recourse to the classification of irreducible complex reflection
groups:

\begin{thm} [\rm Bessis \cite{Bes3}] \label{Gthm:bessis}
 Let $W\leq\GL_n(\CC)$ be an irreducible complex reflection group with braid
 group $B(W)$. Let $d$ be a degree of $W$ which is a regular number for $W$
 and let $r:=(N+N^*)/d$. Then $r\in\NN$, and there exists a subset
 $\bS=\{\bs_1,\ldots,\bs_r\}\subset B(W)$ with:
 \begin{itemize}
 \item[(i)] $\bs_1,\ldots,\bs_r$ are braid reflections, so their
 images $s_1,\ldots,s_r\in W$ are reflections.
 \item[(ii)] $\bS$ generates $B(W)$, and hence $S:=\{s_1,\ldots,s_r\}$
 generates $W$.
 \item[(iii)] There exists a finite set $\cR$ of relations of the form
 $w_1=w_2$, where $w_1,w_2$ are words of equal length in
 $\bs_1,\ldots,\bs_r$, such that $\langle \bs_1,\ldots,\bs_r\mid \cR\rangle$
 is a presentation for $B(W)$.
 \item[(iv)] Let $e_s$ denote the order of $s\in S$. Then
 $\langle s_1,\ldots,s_r\mid \cR; s^{e_s}=1\ \forall s\in S\rangle$
 is a presentation for $W$, where now $\cR$ is viewed as a set of relations
 on $S$.
 \item[(v)] $(\bs_1\cdots\bs_r)^d$ is central in $B(W)$ and lies in $P(W)$.
 \item[(vi)] The product $c:=s_1\cdots s_r$ is a $\zeta:=\exp(2\pi i/d)$-regular
 element of $W$ (hence has eigenvalues $\zeta^{-m_1},\ldots,\zeta^{-m_r}$).
 \end{itemize}
\end{thm}

It follows from the classification (see Table~\ref{Gtab1}) that there
always exists a regular degree. In many cases, for example if $W$ is
well-generated, the number $(N+N^*)/r$ is regular, when $r$ is chosen as the
minimal number of generating reflections for $W$ (so
$n\leq r\leq n+1$). Thus, in those cases $B(W)$ is finitely presented on the
same minimal number of generators as $W$. Under the assumptions~(i) or (ii) of
Theorem~\ref{Gthm:orlik}, the largest degree $d_n$ is regular, whence
Theorem~\ref{Gthm:orlik}~(iv) is a consequence of the previous theorem. \par
At present, presentations of the type described in Theorem~\ref{Gthm:bessis}
have been found for all but six irreducible types, by case-by-case
considerations, see Bannai \cite{Ban}, Naruki \cite{Nar},
Brou\'e--Malle--Rouquier \cite{BMR}. For the remaining six groups, conjectural
presentations have been found by Bessis and Michel using computer calculations.
\par
For the case of real reflection groups, Brieskorn \cite{Br71} and Deligne
\cite{De} determined the structure of $B(W)$ by a nice geometric argument.
They show that the
generators in Theorem~\ref{Gthm:bessis} (with $r=n$) can be taken as suitable
preimages of the Coxeter generators, and the relations $\cR$ as the Coxeter
relations. For the case of $W(A_n)=\fS_{n+1}$ of the classical braid group,
this was first shown by Artin \cite{Ar}.

A topological space $X$ is called $K(\pi,1)$ if all homotopy groups $\pi_i(X)$
for $i\ne1$ vanish. The following is conjectured by Arnol'd to be true for all
irreducible complex reflection groups:

\begin{thm} \label{Gthm:Kpi1}
 Assume that $W$ is not of type $G_i$, $i\in\{24,27,29,31,33,34\}$. Then
 $\regV$ and $\regV/W$ are $K(\pi,1)$-spaces.
\end{thm}

This was proved by a general argument for Coxeter groups by Deligne \cite{De},
after Fox and Neuwirth \cite{FoNe} showed it for type $A_n$ and Brieskorn
\cite{Br73} for those of type different from $H_3, H_4, E_6, E_7, E_8$. For
the non-real Shephard groups (non-real groups with Coxeter braid diagrams),
it was proved by Orlik and Solomon \cite{OrSo88b}. The case of the infinite
series $G(de,e,r)$ has been solved by Nakamura \cite{Nak83}. In that case, there
exists a locally trivial fibration 
$$\regV(G(de,e,n))\longrightarrow \regV(G(de,e,n-1)),$$
with fiber isomorphic to $\CC$ minus $m(de,e,n)$ points, where
$$m(de,e,n):=\left\{\begin{array}{cc} (n-1)de+1& \text{for } d\ne1,\cr
  (n-1)(e-1)& \text{for }d=1.\cr\end{array}\right.$$
This induces a split exact sequence
$$1\longrightarrow F_m\longrightarrow P(G(de,e,n))\longrightarrow
  P(G(de,1,n-1))\longrightarrow1$$
for the pure braid group, with a free group $F_m$ of rank~$m=m(de,e,n)$.
In particular, the pure braid group has the structure of an iterated semidirect
product of free groups (see Brou\'e--Malle--Rouquier
\cite[Prop.~3.37]{BMR}).

\begin{para}{The center and regular elements} \label{subsec:3.regular}
Denote by $\bpi$ the class in $P(W)$ of the loop
$$[0,1]\rightarrow\regV,\qquad t\mapsto x_0\exp(2\pi it).$$
Then $\bpi$ lies in the center $Z(P(W))$ of the pure braid group. Furthermore,
$$[0,1]\rightarrow\regV,\qquad t\mapsto x_0\exp(2\pi it/|Z(W)|),$$
defines a closed path in $\regV/W$, so an element $\bbeta$ of $B(W)$, which is
again central. Clearly $\bpi=\bbeta^{|Z(W)|}$.
\end{para}

The following was shown independently by Brieskorn--Saito \cite{BrSa} and
Deligne \cite{De} for Coxeter groups, and by Brou\'e--Malle--Rouquier
\cite[Th.~2.24]{BMR} for the other groups:

\begin{thm} \label{Gthm:center}
 Assume that $W$ is not of type $G_i$, $i\in\{24,27,29,31,33,34\}$. Then
 the center of $B(W)$ is infinite cyclic generated by $\bbeta$, the
 center of $P(W)$ is infinite cyclic generated by $\bpi$, and the exact
 sequence~(\ref{GEq:braid}) induces an exact sequence
 $$1\longrightarrow Z(P(W))\longrightarrow Z(B(W))\longrightarrow Z(W)
   \longrightarrow1.$$
\end{thm}

In their papers, Brieskorn--Saito \cite{BrSa} and Deligne \cite{De} also
solve the word problem and the conjugation problem for braid
groups attached to real reflection groups.

For each $H\in\cA$ choose a linear form $\alpha_H:V\rightarrow\CC$ with kernel
$H$. Let $e_H:=|W_H|$, the order of the minimal parabolic subgroup fixing $H$.
The {\em discriminant of $W$}, defined as
$$\delta:=\delta(W):=\prod_{H\in\cA} \alpha_H^{e_H},$$
is then a $W$-invariant element of the symmetric algebra $S(V^*)$ of $V^*$,
well-defined up to non-zero scalars (Cohen \cite[1.8]{Co}).
It thus induces a continuous function
$\delta:\regV/W\rightarrow\CC^\times$, hence by functoriality a group
homomorphism $\pi_1(\delta):B(W)\rightarrow\pi_1(\CC^\times,1)\cong\ZZ$. For
$\bb\in B(W)$ let $l(\bb):=\pi_1(\delta)(\bb)$ denote the {\em length of
$\bb$}. For example, every braid reflection $\bs$ has length~$l(\bs)=1$, and
we have
$$l(\bbeta)=(N+N^*)/|Z(W)|\quad\text{ and hence }\quad l(\bpi)=N+N^*$$
by Brou\'e--Malle--Rouquier \cite[Cor.~2.21]{BMR}. \par
The elements $\bb\in B(W)$ with $l(\bb)\geq0$ form the {\em braid monoid}
$B^+(W)$. A {\em $d$th root of $\bpi$} is by definition an element
$\bw\in B^+$ with $\bw^d=\bpi$. 

Let $d$ be a regular number for $W$, and $\bw$ a $d$th root of $\bpi$.
Assume that the image $w$ of $\bw$ in $W$ is $\zeta$-regular for some $d$th
root of unity $\zeta$ (in the sense of~\ref{subsec:1.regular}). (This is for
example the case if $W$ is a Coxeter group by Brou\'e-Michel
\cite[Thm.~3.12]{BrMi}.)
By Theorem~\ref{Gthm:springerreg}(ii) the centraliser $W(w):=C_W(w)$ is a
reflection group on $V(w,\zeta)$, with reflecting hyperplanes the
intersections of $V(w,\zeta)$ with the hyperplanes in $\cA$ by
Theorem~\ref{Gthm:lehrerspringer}(i). Thus the hyperplane complement of $W(w)$
on $V(w,\zeta)$ is just $\regV(w):=\regV\cap V(w,\zeta)$.
Assuming that the base point $x_0$ has been chosen in $\regV(w)$, this defines
natural maps
$$P(W(w))\rightarrow P(W)\qquad\text{and}\qquad
  \psi_w:B(W(w))\rightarrow B(W).$$
By Brou\'e--Michel \cite[3.4]{BrMi} the image of $B(W(w))$ in $B(W)$ 
centralises $\bw$. It is conjectured (see Bessis--Digne--Michel
\cite[Conj.~0.1]{BDM}) that $\psi_w$ defines an
isomorphism $B(W(w)) \cong C_{B(W)}(\bw)$. The following partial
answer is known:

\begin{thm} [\rm Bessis--Digne--Michel \protect{\cite[Th.~0.2]{BDM}}]
 \label{Gthm:braidiso}
 Let $W$ be an irreducible reflection group of type $\fS_n$, $G(d,1,n)$ or
 $G_i$, $i\in \{4,5,8,10,16,18,25,26,32\}$, and let $w\in W$ be regular. Then
 $\psi_w$ induces an isomorphism $B(W(w)) \cong C_{B(W)}(\bw)$
\end{thm}

This has also been proved by Michel \cite[Cor.~4.4]{Mi99} in the case that
$W$ is a Coxeter group and $w$ acts on $W$ by a diagram automorphism. The
injectivity of $\psi_w$ was shown for all but finitely many types of $W$ by
Bessis \cite[Th.~1.3]{Bes2}.

The origin of Artin's work on the braid group  associated with the
symmetric group lies in the theory of knots and links. We shall now
briefly discuss this connection and explain the construction of
the ``HOMFLY--PT'' invariant of knots and links (which includes
the famous Jones polynomial as a special case). We follow the exposition 
in Geck--Pfeiffer \cite[\S 4.5]{GePf}.

\begin{para}{Knots and links, Alexander and Markov theorem}
\label{subsec:3.knots} If $n$ is a positive integer, an oriented
$n$-{\em link} is an embedding of~$n$ copies of the interval $[0,1] 
\subset \RR$ into ${\RR}^3$ such that $0$ and~$1$ are mapped to the same
point (the orientation is induced by the natural ordering of $[0,1]$); a
$1$-link is also called a {\em knot}. We are only interested in knots and
links modulo isotopy, i.e., homeomorphic transformations which preserve
the orientation. We refer to Birman \cite{Birman74}, Crowell--Fox 
\cite{CroFox} or Burde--Zieschang \cite{BuZi} for precise versions of 
the above definitions. 

By Artin's classical interpretation of $B(\fS_n)$ as the braid group on 
$n$ strings, each generator of $B(\fS_n)$ can be represented by oriented 
diagrams as indicated below; writing any $g \in B(\fS_n)$ as a product of 
the generators and their inverses, we also obtain a diagram for~$g$, by 
concatenating the diagrams for the generators.  ``Closing'' such a diagram 
by joining the end points, we obtain the plane projection of an oriented 
link in~${\RR}^3$:
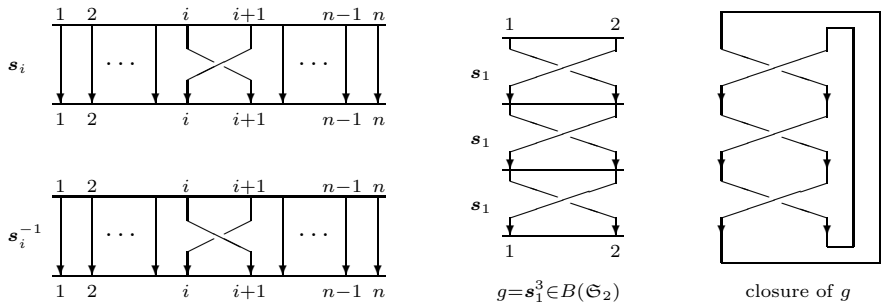
\begin{center} \unitlength1pt
\vbox{\begin{picture}(330,115)
\put(  0,23){$\scriptstyle{\braid{s}_i^{{-}1}}$}
\put( 17,10){\line(1,0){126}}
\put( 17,40){\line(1,0){126}}
\put( 20,40){\vector(0,-1){30}}
\put( 32,40){\vector(0,-1){30}}
\put( 37,25){$\ldots$}
\put( 56,40){\vector(0,-1){30}}
\put(104,40){\vector(0,-1){30}}
\put(110,25){$\ldots$}
\put(128,40){\vector(0,-1){30}}
\put(140,40){\vector(0,-1){30}}
\put( 68,40){\line(0,-1){9}}
\put( 92,40){\line(0,-1){9}}
\put( 68,19){\vector(0,-1){9}}
\put( 92,19){\vector(0,-1){9}}
\put( 92,31){\line(-2,-1){10}}
\put( 68,19){\line(2,1){10}}
\put( 68,31){\line(2,-1){24}}
\put( 18, 2){$\scriptstyle{1}$}
\put( 30, 2){$\scriptstyle{2}$}
\put( 66, 2){$\scriptstyle{i}$}
\put( 85, 2){$\scriptstyle{i{+}1}$}
\put(119, 2){$\scriptstyle{n{-}1}$}
\put(138, 2){$\scriptstyle{n}$}
\put( 18,42){$\scriptstyle{1}$}
\put( 30,42){$\scriptstyle{2}$}
\put( 66,42){$\scriptstyle{i}$}
\put( 85,42){$\scriptstyle{i{+}1}$}
\put(119,42){$\scriptstyle{n{-}1}$}
\put(138,42){$\scriptstyle{n}$}

\put(  0, 88){$\scriptstyle{\braid{s}_i}$}
\put( 17, 75){\line(1,0){126}}
\put( 17,105){\line(1,0){126}}
\put( 20,105){\vector(0,-1){30}}
\put( 32,105){\vector(0,-1){30}}
\put( 37, 90){$\ldots$}
\put( 56,105){\vector(0,-1){30}}
\put(104,105){\vector(0,-1){30}}
\put(110, 90){$\ldots$}
\put(128,105){\vector(0,-1){30}}
\put(140,105){\vector(0,-1){30}}
\put( 68,105){\line(0,-1){9}}
\put( 92,105){\line(0,-1){9}}
\put( 68, 84){\vector(0,-1){9}}
\put( 92, 84){\vector(0,-1){9}}
\put( 68, 96){\line(2,-1){10}}
\put( 92, 84){\line(-2,1){10}}
\put( 92, 96){\line(-2,-1){24}}
\put( 18, 67){$\scriptstyle{1}$}
\put( 30, 67){$\scriptstyle{2}$}
\put( 66, 67){$\scriptstyle{i}$}
\put( 85, 67){$\scriptstyle{i{+}1}$}
\put(119, 67){$\scriptstyle{n{-}1}$}
\put(138, 67){$\scriptstyle{n}$}
\put( 18,107){$\scriptstyle{1}$}
\put( 30,107){$\scriptstyle{2}$}
\put( 66,107){$\scriptstyle{i}$}
\put( 85,107){$\scriptstyle{i{+}1}$}
\put(119,107){$\scriptstyle{n{-}1}$}
\put(138,107){$\scriptstyle{n}$}
\put(185,  2){$\scriptstyle{g{=}\braid{s}_1^3\in B(\fS_2)}$}
\put(175, 35){$\scriptstyle{\braid{s}_1}$}
\put(175, 60){$\scriptstyle{\braid{s}_1}$}
\put(175, 85){$\scriptstyle{\braid{s}_1}$}
\put(187, 25){\line(1,0){46}}
\put(187, 50){\line(1,0){46}}
\put(187, 75){\line(1,0){46}}
\put(187,100){\line(1,0){46}}
\put(190,100){\line(0,-1){4}}
\put(230,100){\line(0,-1){5}}
\put(190, 82){\vector(0,-1){7}}
\put(230, 82){\vector(0,-1){7}}
\put(190, 75){\line(0,-1){4}}
\put(230, 75){\line(0,-1){5}}
\put(190, 57){\vector(0,-1){7}}
\put(230, 57){\vector(0,-1){7}}
\put(190, 50){\line(0,-1){4}}
\put(230, 50){\line(0,-1){5}}
\put(190, 32){\vector(0,-1){7}}
\put(230, 32){\vector(0,-1){7}}
\put(190, 96){\line(3,-1){18}}
\put(230, 82){\line(-3,1){17}}
\put(190, 82){\line(3,1){40}}
\put(190, 71){\line(3,-1){18}}
\put(230, 57){\line(-3,1){17}}
\put(190, 57){\line(3,1){40}}
\put(190, 46){\line(3,-1){18}}
\put(230, 32){\line(-3,1){17}}
\put(190, 32){\line(3,1){40}}
\put(188, 17){$\scriptstyle{1}$}
\put(228, 17){$\scriptstyle{2}$}
\put(188,103){$\scriptstyle{1}$}
\put(228,103){$\scriptstyle{2}$}

\put(279,  2){$\scriptstyle{\text{closure of $g$}}$}
\put(270,100){\line(0,1){10}}
\put(270,110){\line(1,0){60}}
\put(270, 15){\line(0,1){10}}
\put(270, 15){\line(1,0){60}}
\put(330, 15){\line(0,1){95}}
\put(310,100){\line(0,1){4}}
\put(310,104){\line(1,0){10}}
\put(310, 25){\line(0,-1){4}}
\put(310, 21){\line(1,0){10}}
\put(320, 21){\line(0,1){83}}
\put(270,100){\line(0,-1){4}}
\put(310,100){\line(0,-1){5}}
\put(270, 82){\vector(0,-1){7}}
\put(310, 82){\vector(0,-1){7}}
\put(270, 75){\line(0,-1){4}}
\put(310, 75){\line(0,-1){5}}
\put(270, 57){\vector(0,-1){7}}
\put(310, 57){\vector(0,-1){7}}
\put(270, 50){\line(0,-1){4}}
\put(310, 50){\line(0,-1){5}}
\put(270, 32){\vector(0,-1){7}}
\put(310, 32){\vector(0,-1){7}}
\put(270, 96){\line(3,-1){18}}
\put(310, 82){\line(-3,1){17}}
\put(270, 82){\line(3,1){40}}
\put(270, 71){\line(3,-1){18}}
\put(310, 57){\line(-3,1){17}}
\put(270, 57){\line(3,1){40}}
\put(270, 46){\line(3,-1){18}}
\put(310, 32){\line(-3,1){17}}
\put(270, 32){\line(3,1){40}}
\end{picture}}
\end{center}
By Alexander's theorem (see Birman \cite{Birman74} or, for a more recent
proof,  Vogel \cite{Vogel90}), every oriented link in ${\RR}^3$ is isotopic
to the closure of an element in $B(\fS_n)$,
for some $n \geq 1$. The question of when two links in ${\RR}^3$ are isotopic
can also be expressed algebraically. For this purpose, we consider the
infinite disjoint union
\[ B_{\infty}:=\coprod_{n \geq 1} B(\fS_n).\]
Given $g,g' \in B_{\infty}$, we write $g \sim g'$ if one of the following
relations is satisfied:
\begin{itemize}
\item[(I)] We have $g,g' \in B(\fS_n)$ and $g'=x^{-1}gx$ for some $x \in 
B(\fS_n)$.
\item[(II)] We have $g \in B(\fS_n)$, $g' \in B(\fS_{n+1})$ and
$g'=g\braid{s}_n$ or $g'=g\braid{s}_n^{-1}$.
\end{itemize}
The above two relations are called {\em Markov relations}. By a classical
result due to Markov (see Birman \cite{Birman74} or, for a more recent
proof, Traczyk \cite{TracBanach}), two elements of $B_{\infty}$ are equivalent 
under the equivalence relation generated by~$\sim$ if and only if the 
corresponding links obtained by closure are isotopic. Thus, to define an 
invariant of oriented links is the same as to define a map on $B_{\infty}$
which takes equal values on elements $g,g' \in B_{\infty}$ satisfying~(I) 
or~(II).
\end{para}

We now consider the Iwahori--Hecke algebra $H_{\CC}(\fS_n)$ of the
symmetric group $\fS_n$ over $\CC$. By definition, $H_{\CC}(\fS_n)$ is 
a quotient of the group algebra of $B(\fS_n)$, where we factor by an 
ideal generated by certain quadratic relations depending on two parameters
$u,v\in \CC$. This is done such that
\[ T_{s_i}^2=uT_1+vT_{s_i} \qquad \mbox{for $1\leq i\leq n-1$},\]
where $T_{s_i}$ denotes the image of  the generator $\bs_i$ of
$B(\fS_n)$ and $T_1$ denotes the identity element. For each $w\in \fS_n$, we 
have a well-defined element $T_w$ such that 
\[T_{w}T_{w'} = T_{ww'} \qquad \mbox{whenever $l(ww')=l(w)+l(w')$}.\]
This follows easily from Matsumoto's theorem~\ref{Mmatsum}. In fact, one can
show
that the elements $\{T_w\mid w\in \fS_n\}$ form a $\CC$-basis of
$H_{\CC}(\fS_n)$. (For more details, see the chapter on Hecke
algebras.) The map $\braid{w}\mapsto T_w$ ($w\in\fS_n$) extends to a 
well-defined algebra homomorphism from the group algebra of $B(\fS_n)$ 
over $\CC$ onto $H_{\CC}(\fS_n)$. Furthermore, the inclusion $\fS_{n-1} 
\subseteq \fS_n$ also defines an inclusion of algebras $H_{\CC}(\fS_{n-1}) 
\subseteq H_{\CC} (\fS_n)$. 

\begin{thm}[Jones, Ocneanu \protect{\cite{Jones87}}] \label{Mjones} 
There is a unique family of $\CC$-linear maps $\tau_n \colon H_{\CC}(\fS_n) 
\rightarrow \CC$ ($n \geq 1$) such that the following conditions hold:
\begin{itemize}
\item[(M1)] $\tau_1(T_1)=1$;
\item[(M2)] $\tau_{n+1}(hT_{s_n})=\tau_{n+1}(hT_{s_n}^{-1})=\tau_{n}(h)$ for
all $n \geq 1$ and $h \in H_{\CC}(\fS_n)$; 
\item[(M3)] $\tau_n(hh')=\tau_n(h'h)$ for all $n \geq 1$ and $h,h' \in 
H_{\CC}(\fS_n)$.
\end{itemize}
Moreover, we have $\tau_{n+1}(h)=v^{-1}(1-u)\tau_n(h)$ for all $n \geq 1$
and $h \in H_{\CC}(\fS_n)$.
\end{thm}

In \cite{Jones87}, Jones works with an Iwahori--Hecke algebra of $\fS_n$ 
where the parameters are related by $v=u-1$. The different formulation above
follows a suggestion by J.~Michel. It results in a simplification of the
construction of the link invariants below. (The simplification arises from 
the fact that, due to the presence of two different parameters in the 
quadratic relations, the ``singularities'' mentioned in 
 \cite[p.~349, Notes~(1)]{Jones87} simply disappear.) Generalisations of 
Theorem~\ref{Mjones} to types $G(d,1,n)$ and $D_n$ have been found in 
Geck--Lambropoulou \cite{GeLa}, Lambropoulou \cite{Lamb2} and Geck 
\cite{GeckBanach}. 

\begin{para}{The HOMFLY-PT polynomial} \label{Mhomflpt} We can now
construct a two-variable invariant of oriented knots and links as follows.
Consider an oriented link $L$ and assume that it is isotopic to the closure
of $g \in B(\fS_n)$ for $n \geq 1$. Then we set 
\[ X_{L}(u,v):=\tau_n(\bar{g}) \in {\CC} \qquad 
\mbox{with $\tau_n$ as in Theorem~\ref{Mjones}}.\] 
Here, $\bar{g}$ denotes the image of~$g$ under the natural map ${\CC}
[B(\fS_n)] \rightarrow H_{\CC}(\fS_n)$, $\braid{w} \mapsto T_w$ ($w \in 
\fS_n$). It is easily checked that $X_L(u,v)$ can be expressed as a
Laurent polynomial in $u$ and $v$; the properties (M2) and (M3) make sure 
that $\tau_n(\bar{g})$ does not depend on the choice of~$g$. If we make the 
change of variables $u=t^2$ and $v=tx$, we can identify the above invariant 
with the HOMFLY-PT polynomial $P_{L}(t,x)$ discovered by Freyd et al.
\cite{HOMFLY} and Przytycki--Traczyk \cite{PT}; see also Jones 
\cite[(6.2)]{Jones87}. Furthermore, the {\em Jones polynomial} $J_L(t)$ is 
obtained by setting $u=t^2$, $v=\sqrt{t}(t-1)$ (see \cite[\S 11]{Jones87}). 
Finally, setting $u=1$ and $v=\sqrt{t}-1/\sqrt{t}$, we obtain the classical 
{\em Alexander polynomial} $A_L(t)$ whose definition can be found in 
Crowell--Fox \cite{CroFox}.
\end{para}

For a survey about recent developments in the theory of knots and links, 
especially since the discovery of the Jones polynomial, see Birman 
\cite{Birman93}.

\begin{para}{Further aspects of braid groups}  \label{linbraid}
One of the old problems concerning braid groups is the question
whether or not they are linear, i.e., whether there exists a 
faithful linear representation on a finite dimensional vector
space.  Significant progress has been made recently on this
problem. Krammer \cite{Kram2} and Bigelow \cite{Big2} proved that the 
classical Artin braid group is linear. Then Digne \cite{Digne} 
and Cohen--Wales \cite{CoWa} extended this result and showed that all 
Artin groups of crystallographic type have a faithful representation of 
dimension equal to the number of reflections of the associated Coxeter 
group. 

On the other hand, there is one particular representation of the braid 
group associated with $\fS_n$, the so-called Burau representation (see 
Birman \cite{Birman74}, for which it has been a long-standing problem to 
determine for which values of $n$ it is faithful. Moody \cite{Moody93}
showed that it is not faithful for $n\geq 10$; this bound was improved
by Long and Paton \cite{LoPa} to $6$. Recently, Bigelow \cite{Big1} showed 
that the Burau representation is not faithful already for $n=5$. (It is 
an old result of Magnus and Peluso that the Burau representation is
faithful for $n=3$.)

In a different direction, Deligne's and Brieskorn--Saito's solution of
the word and conjugacy problem in braid groups led to new developments 
in combinatorial group and monoid theory; see, for example, Dehornoy and 
Paris \cite{DehPar} and Dehornoy \cite{Dehor}.
\end{para}

\section{ Representation theory} \label{sec14}

In this section we report about the representation theory of finite complex
reflection groups.

\begin{para}{Fields of definition}
\label{subsec:4.deffields}
Let $W$ be a finite complex reflection group on $V$. Let $K_W$ denote the
character field of the reflection representation of $W$, that is, the field
generated by the traces $\tr_V(w)$, $w\in W$. It is easy to see that the
reflection representation can be realised over $K_W$ (see for example
\cite[Prop.~7.1.1]{Ben}). But we have a much stronger statement:
\end{para}

\begin{thm} [\rm Benard \cite{Ben76}, Bessis \cite{Bes}] \label{Gthm:BenBes}
 Let $W$ be a complex reflection group. Then the field $K_W$ is a splitting
 field for $W$.
\end{thm}

The only known proof for this result is case-by-case, treating the reflection
groups according to the Shephard--Todd classification.

The field $K_W$ has a nice description at least in the case of well-generated
groups, that is, irreducible groups generated by $\dim(V)$ reflections. In
this case, the largest degree
$d_n$ of $W$ is regular, so there exists an element $c:=s_1\cdots s_n$
as in the Theorem~\ref{Gthm:bessis}(vi) of Bessis, called a {\em Coxeter
element} of $W$, with eigenvalues $\zeta^{-m_1},\ldots,\zeta^{-m_n}$ in the
reflection representation, where $\zeta:=\exp(2\pi i/d_n)$ and the $m_i$ are
the exponents of $W$.

If $W$ is a real reflection group, then $W$ is a Coxeter group associated with
some Coxeter matrix $M$ (see Theorem~\ref{Mthm1}) and we have
\[ K_W={\QQ}(\cos(2\pi/m_{st}) \mid s,t \in S) \subset \RR.\]
In particular, this shows that $K_W=\QQ$ if $W$ is a finite Weyl group.

For well-generated irreducible complex reflection groups $W\leq\GL(V)$, the
field of definition $K_W$ is generated over $\QQ$ by the coefficients of the
characteristic polynomial on $V$ of a Coxeter element, see Malle 
\cite[Th.~7.1]{MaR}. This characterisation is no longer true for non-well
generated reflection groups.

\begin{para}{Macdonald--Lusztig--Spaltenstein induction}
\label{subsec:4.j-induction}
Let $W$ be a complex reflection group on $V$. Recall from
Section~\ref{subsec:1.fake} the definition of the fake degree $R_\chi$ of an
irreducible character $\chi\in\Irr(W)$. The {\em $b$-invariant $b_\chi$ of
$\chi$} is defined as the order of vanishing of $R_\chi$ at $x=0$, that is,
as the minimum of the exponents $e_i(\chi)$ of $\chi$. The coefficient of
$x^{b_\chi}$ in $R_\chi$ is denoted by $\gamma_\chi$. 
\end{para}

\begin{thm} [\rm Macdonald \cite{Mac72}, Lusztig--Spaltenstein \cite{LuSp79}]
 \label{Gthm:j-induction}
 Let $W$ be a complex reflection group on $V$, $W'$ a reflection subgroup
 (on $V_1:=V/C_V(W')$). Let $\psi$ be an irreducible character of $W'$ such
 that $\gamma_\psi=1$. Then $\Ind_{W'}^W(\psi)$ has a unique irreducible
 constituent $\chi\in\Irr(W)$ with $b_\chi=b_\psi$. This satisfies
 $\gamma_\chi=1$. All other constituents have $b$-invariant bigger than
 $b_\psi$
\end{thm}

The character $\chi\in\Irr(W)$ in Theorem~\ref{Gthm:j-induction} is called
the {\em $j$-induction} $j_{W'}^W(\psi)$ of the character $\psi\in\Irr(W')$.
Clearly, $j$-induction is transitive; it is also compatible with direct
products. \par
An important example of characters $\psi$ with $\gamma_\psi=1$ is given by the
determinant character $\det_V:W\rightarrow\CC^\times$ of a complex reflection
group (see Geck--Pfeiffer \cite[Th.~5.2.10]{GePf}).

\begin{para}{Irreducible characters}
\label{subsec:4.characters}
There is no general construction of all irreducible representations of a
complex reflection group known. Still, we have the following partial result:

\begin{thm} [\rm Steinberg]
 \label{Gthm:extpower}
 Let $W\leq\GL(V)$ be an irreducible complex reflection group. Then the
 exterior powers $\Lambda^i(V)$, $1\leq i\leq \dim V$, are irreducible,
 pairwise non-equivalent representations of $W$.
\end{thm}

A proof in the case of well-generated groups can be found in Bourbaki
\cite[V, \S2, Ex.~3(d)]{Bou} and Kane \cite[Thm.~24-3 A]{Kane}, for example.
The general case then follows with Corollary~\ref{Gcor:irrwell}.
\par

We now give some information on the characters of individual reflection groups. 
The irreducible characters of the symmetric group $\fS_n$ were determined by
Frobenius \cite{Frob}, see also Macdonald \cite{Mac95} and Fulton 
\cite{Fult}. Here we follow the exposition in Geck--Pfeiffer 
\cite[5.4]{GePf}. \par
Let $\lambda=(\lambda_1,\ldots,\lambda_r)\vdash n$ be a partition of $n$. The
corresponding {\em Young subgroup $\fS_\lambda$ of $\fS_n$} is the common
setwise stabiliser
$\{1,\ldots,\lambda_1\}$, $\{\lambda_1+1,\ldots,\lambda_1+\lambda_2\},\ldots$,
abstractly isomorphic to
$\fS_\lambda=\fS_{\lambda_1}\times\ldots\times\fS_{\lambda_r}$.
This is a parabolic subgroup of $\fS_n$ in the sense of
Section~\ref{subsec:1.parabolic}. For any $m\geq1$,
let $1_m$, $\eps_m$ denote the trivial respectively the sign character of
$\fS_m$. For each partition we have the two induced characters
$$\pi_\lambda:=\Ind_{\fS_\lambda}^{\fS_n}(1_{\lambda_1}\#\ldots\#1_{\lambda_r}),
  \qquad \theta_\lambda:=\Ind_{\fS_\lambda}^{\fS_n}
  (\eps_{\lambda_1}\#\ldots\#\eps_{\lambda_r}).$$
Then $\pi_\lambda$ and $\theta_{\lambda^*}$ have a unique irreducible
constituent $\chi_\lambda\in\Irr(\fS_n)$ in common, where $\lambda^*$ denotes
the partition dual to $\lambda$. This constituent can also
be characterised in terms of $j$-induction as
$$\chi_\lambda=j_{\fS_\lambda}^{\fS_n}(1_{\lambda_1}\#\ldots\#1_{\lambda_r}).$$
Then the $\chi_\lambda$ are mutually distinct and exhaust the irreducible
characters of $\fS_n$, so $\Irr(\fS_n)=\{\chi_\lambda\mid\lambda\vdash n\}$.
From the above construction it is easy to see that all $\chi_\lambda$ are
afforded by rational representations.
\par

The construction of the irreducible characters of the imprimitive group
$G(d,1,n)$ goes back at least to Osima \cite{Osi54} (see also Read
\cite{Read77}, Hughes \cite{Hu83}, Bessis \cite{Bes}) via their abstract
structure as wreath product $ C_d\wr\fS_n$. Let us fix $d\geq2$ and write
$W_n:=G(d,1,n)$. 
A $d$-tuple $\alpha=(\alpha_0,\alpha_1,\ldots,\alpha_{d-1})$ of partitions
$\alpha_i\vdash n_i$ with $\sum n_i=n$ is called a {\em $d$-partition of $n$}. 
We denote by $W_\alpha$ the natural subgroup
$W_{n_0}\times\ldots\times W_{n_{d-1}}$ of $W_n$, where $\alpha_j\vdash n_j$,
corresponding to the Young subgroup $\fS_{n_0}\times\ldots\times\fS_{n_{d-1}}$
of $\fS_n$. Via the natural projection $W_{n_j}\rightarrow\fS_{n_j}$ the
characters of $\fS_{n_j}$ may be regarded as characters of $W_{n_j}$. Thus
each $\alpha_j$ defines an irreducible character $\chi_{\alpha_j}$ of $W_{n_j}$.
For any $m$, let $\zeta_d:W_m\rightarrow\CC^\times$ be the linear character
defined by $\zeta_d(t_1)=\exp(2\pi i/d)$, $\zeta_d(t_i)=1$ for $i>1$,
with the standard generators $t_i$ from~\ref{subsec:1.Class}. Then for
any $d$-partition $\alpha$ of $n$ we can define a character
$\chi_\alpha$ of $W_n$ as the induction of the exterior product 
$$\chi_\alpha:=\Ind_{W_\alpha}^{W_n}\left(\chi_{\alpha_0}\#
  (\chi_{\alpha_1}\otimes\zeta_d)\# \ldots\#
  (\chi_{\alpha_{d-1}}\otimes\zeta_d^{d-1})\right)\,.$$
By Clifford-theory $\chi_\alpha$ is irreducible, $\chi_\alpha\ne\chi_\beta$ if
$\alpha\ne\beta$, and all irreducible characters of $W_n$ arise in this way,
so
$$\Irr(W_n)=\{\chi_\alpha\mid \alpha=(\alpha_0,\ldots,\alpha_{d-1})
{\vdash_d} n\}.$$

We describe the irreducible characters of $G(de,e,n)$ in terms of those of
$W_n:=G(de,1,n)$. Recall that the imprimitive reflection group $G(de,e,n)$ is
generated by the reflections
$t_2,\tilde t_2:=t_1^{-1}t_2t_1,t_3,\ldots,t_n$, and $t_1^e$.
Denote by $\pi$ the cyclic shift on $de$-partitions of $n$,
i.e.,
$$\pi(\alpha_0,\ldots,\alpha_{de-1})=(\alpha_1,\ldots,\alpha_{de-1},\alpha_0)
  \,.$$
By definition we then have $\chi_{\pi(\alpha)}\otimes\zeta_{de}=\chi_\alpha$.
Let $s_e(\alpha)$ denote the order of the stabiliser of $\alpha$ in the
cyclic group $\langle\pi^{d}\rangle$. Then upon restriction to $G(de,e,n)$
the irreducible character $\chi_\alpha$ of $W_n$ splits into $s_e(\alpha)$
different irreducible constituents, and this exhausts the set of irreducible
characters of $G(de,e,n)$. More precisely, let $\alpha$ be a $de$-partition of
$n$ with $\tilde e:=s_e(\alpha)$, $W_{\alpha,e}:=W_\alpha\cap G(de,e,n)$, and
$\psi_\alpha$ the restriction of
$$\chi_{\alpha_0}\# (\chi_{\alpha_1}\otimes\zeta_{de})\# \ldots\#
  (\chi_{\alpha_{de-1}}\otimes\zeta_{de}^{de-1})$$
to $W_{\alpha,e}$. Then $\psi_\alpha$
is invariant under the element $\sigma:=(t_2\cdots t_n)^{n/\tilde e}$
(note that $\tilde e=s_e(\alpha)$ divides $n$), and it extends to the
semidirect product $W_{\alpha,e}.\langle\sigma\rangle$. The different
extensions of $\psi_\alpha$ induced to $G(de,e,n)$ then exhaust the irreducible
constituents of the
restriction of $\chi_\alpha$ to $G(de,e,n)$. Thus, we may parametrise
$\Irr(G(de,e,n))$ by $de$-partitions of $n$ up to cyclic shift by
$\pi^{d}$ in such a way that any $\alpha$ stands for $s_e(\alpha)$
different characters. \par

In order to describe the values of the irreducible characters we need
the following definitions. We identify partitions with their Young diagrams.
A $d$-partition $\alpha$ is called a {\em hook} if it has just one
non-empty part, which is a hook (i.e., does not contain a $2\times2$-block).
The position of the non-empty part is then denoted by $\tau(\alpha)$. 
If $\alpha$, $\beta$ are $d$-partitions such that $\beta_i$ is
contained in $\alpha_i$ for all $0\leq i\leq d-1$, then $\alpha\setminus\beta$
denotes the $d$-partition
$(\alpha_i\setminus\beta_i\mid 0\leq i\leq d-1)$, where
$\alpha_i\setminus\beta_i$ is the set theoretic difference of $\alpha_i$ and
$\beta_i$. If $\alpha\setminus\beta$ is a hook, we denote by $l_\beta^\alpha$
the number of rows of the hook
$(\alpha\setminus\beta)_{\tau(\alpha\setminus\beta)}$ minus~1.
With these notations the values of the irreducible characters of $G(d,1,n)$
can be computed recursively with a generalised Murnaghan--Nakayama rule
(see also Stembridge \cite{Ste89}):

\begin{thm} [\rm Osima \cite{Osi54}]
 Let $\alpha$ and $\gamma$ be $d$-partitions of $n$, let $m$ be a
 part of $\gamma_t$ for some $1\leq t\leq d$, and denote by $\gamma'$ the
 $d$-partition of $n-m$ obtained from $\gamma$ by deleting the part
 $m$ from $\gamma_t$. Then the value of the irreducible character $\chi_\alpha$
 on an element of $G(d,1,n)$ with cycle structure $\gamma$ is given by
 $$\chi_\alpha(\gamma) = \sum_{\alpha\setminus\beta\vdash_d m}
   \zeta_d^{st}(-1)^{l_\beta^\alpha}\chi_\beta(\gamma')$$
 where the sum ranges over all $d$-partitions $\beta$ of $n-m$ such
 that $\alpha\setminus\beta$ is a hook and where
 $s=\tau(\alpha\setminus\beta)$.
\end{thm}

An overview of the irreducible characters of the exceptional groups (and of
their projective characters) is given in Humphreys \cite{Hum94}. All character
tables of irreducible complex reflection groups are also available in the
computer algebra system {\sf CHEVIE} \cite{chev}.
\end{para}

\begin{para}{Fake degrees}\label{subsec:4.fake}
The fake degrees (introduced in Section~\ref{subsec:1.fake}) of all complex
reflection groups are known. For the symmetric groups, they were first
determined by Steinberg \cite{St51}, as generic degrees of the unipotent
characters of the general linear groups over a finite field. Then Lusztig
\cite{Lu77} determined the fake degrees for finite Coxeter groups of
type $B_n$ and $D_n$. The fake degrees of arbitrary imprimitive reflection 
groups can easily be derived by similar arguments (see Malle \cite[Bem.~2.10
 and~5.6]{MaU}).

\begin{thm} [\rm Steinberg \cite{St51}, Lusztig \cite{Lu77}, Malle \cite{MaU}]
 The fake degrees of the irreducible complex reflection groups $G(de,e,n)$ are
 given as follows:
 \begin{itemize}
 \item[(i)] Let $\chi\in\Irr(G(d,1,n))$ be parameterised by the $d$-partition
 $(\alpha_0,\ldots,\alpha_{d-1})$, where
 $\alpha_i=(\alpha_{i1}\geq\ldots\geq\alpha_{im_i})\vdash n_i$, and let
 $(S_0,\ldots,S_{d-1})$, where $S_i=(\alpha_{i1}+m_i-1,\ldots,\alpha_{im_i})$,
 denote the corresponding tuple of $\beta$-numbers. Then
 $$R_\chi=\prod_{i=1}^n (x^{id}-1)\prod_{i=0}^{d-1}
   \frac{\Delta(S_i,x^d)x^{in_i}}
   {\Theta(S_i,x^d)x^{d\binom{m_i-1}{2}+d\binom{m_i-2}{2}+\ldots}}\,,$$
 where, for a finite subset $S\subset\NN$,
 $$\Delta(S,x):=\prod_{\atop{\lambda,\lambda'\in S}{\lambda'<\lambda}}
   (x^\lambda-x^{\lambda'}),\qquad
   \Theta(S,x):=\prod_{\lambda\in S}\prod_{h=1}^\lambda (x^h-1).$$
 \item[(ii)] The fake degree of $\chi\in\Irr(G(de,e,n))$ is obtained
 from the fake degrees in $G(de,1,n)$ as
 $$R_\chi=\frac{x^{nd}-1}{x^{nde}-1}\sum_{\psi\in\Irr(G(de,1,n))}
  \left\langle\chi,\psi|_{G(de,e,n)}\right\rangle\, R_{\psi}.$$
 \end{itemize}
\end{thm}

For exceptional complex reflection groups, the fake degrees can easily be
computed, for example in the computer algebra system {\sf CHEVIE} \cite{chev}.
In the case of exceptional Weyl groups, they were first studied by
Beynon--Lusztig \cite{BeLu}. \par
The fake degrees of reflection groups satisfy a remarkable palindromicity
property:

\begin{thm} [\rm Opdam \cite{Op95,Op00}, Malle \cite{MaR}]
 \label{Gthm:fakepalindrom}
 Let $W$ be a complex reflection group. There exists a permutation $\delta$ on
 $\Irr(W)$ such that for every $\chi\in\Irr(W)$ we have
 $$R_\chi(x) = x^c R_{\delta(\bar\chi)}(x^{-1}),$$
 where $c= \sum_r (1 - \chi(r)/\chi(1))$, the sum running over all reflections
 $r\in W$.
\end{thm}

The interest of this result also lies in the fact that the permutation $\delta$
is strongly related to the irrationalities of characters of the associated
Hecke algebra. Theorem~\ref{Gthm:fakepalindrom} was first observed empirically
by Beynon--Lusztig \cite{BeLu} in the case of Weyl groups. Here $\delta$ is
non-trivial only for characters such that the corresponding character of the
associated Iwahori-Hecke algebra is non-rational. An a priori proof in this
case was later given by Opdam \cite{Op95}.
If $W$ is complex, Theorem~\ref{Gthm:fakepalindrom} was verified by Malle 
\cite[Thm.~6.5]{MaR} in a case-by-case analysis. Again in all but possibly
finitely many cases, $\delta$ comes from the irrationalities of characters of
the associated Hecke algebra. Opdam \cite[Thm.~4.2 and Cor.~6.8]{Op00} gives
a general argument which proves Theorem~\ref{Gthm:fakepalindrom} under a
suitable assumption on the braid group $B(W)$.

\end{para}

The above discussion is exclusively concerned with representations over 
a field of characteristic $0$ (the ``semisimple case''). We close this
chapter with some remarks concerning the modular case.

\begin{para}{Modular representations of $\fS_n$} \label{modsn} Frobenius'
theory (as developed further by Specht, James and others) yields a 
parametrisation of $\Irr(\fS_n)$ and explicit formulas for the degrees 
and the values of all irreducible characters. As soon as we consider 
representations 
over a field of characteristic $p>0$, the situation changes drastically.
James \cite{James1} showed that the irreducible representations of $\fS_n$
still have a natural parametrisation, by so-called $p$-regular partitions. 
Furthermore, the decomposition matrix relating representations in 
characteristic~$0$ and in characteristic~$p$ has a lower triangular shape 
with~$1$s on the diagonal. This result shows that, in principle, a knowledge 
of the irreducible representations of $\fS_n$ in characteristic~$p$ is 
equivalent to the knowledge of the decomposition matrix.

There are a number of results known about certain entries of that 
decomposition matrix, but a general solution to this problem is completely 
open; see James \cite{James1} for a survey. Via the classical Schur algebras 
(see Green \cite{Green1}) it is known that the decomposition numbers of 
$\fS_n$ in characteristic~$p$ can be obtained from those of the finite 
general linear group $\mbox{GL}_n({\FF}_q)$ (where $q$ is a power of $p$). 
Now, at first sight, the problem of computing the decomposition numbers 
for $\mbox{GL}_n({\FF}_q)$ seems to be much harder than for $\fS_n$. 
However, Erdmann \cite{Erd} has shown that, if one knows the decomposition 
numbers of $\fS_m$ (for sufficiently large values of $m$), then one will
also know the decomposition numbers of $\mbox{GL}_n({\FF}_q)$. Thus, the 
problem of determining the decomposition numbers for symmetric groups 
appears to be as difficult as the corresponding problem for general linear 
groups. 

In a completely different direction, Dipper and James \cite{DJ1} showed
that the decomposition numbers of $\fS_n$ can also be obtained from
the so-called $q$-Schur algebra, which is defined in terms of the
Hecke algebra of $\fS_n$. Thus, the problem of determining the 
representations of $\fS_n$ in characteristic~$p$ is seen to be a special 
case of the more general problem of studying the representations of 
Hecke algebras associated with finite Coxeter groups. This is discussed 
in more detail in the chapter on Hecke algebras. In this context,
we just mention here that James' result on the triangularity of the 
decomposition matrix of $\fS_n$ is generalised to all finite Weyl groups 
in Geck \cite{mykl}.
\end{para}

\section{Hints for further reading} \label{sec15}
Here, we give some hints on topics which were not touched in the previous
sections.

\begin{para}{Crystallographic reflection groups}
\label{subsec:5.crystallographic}
Let $W$ be a discrete subgroup of the group of all affine transformations of
a finite-dimensional affine space $E$ over $K=\CC$ or $K=\RR$, generated by
affine reflections. If $W$ is finite, it necessarily fixes a point and $W$ is
a finite complex reflection group. If $W$ is infinite and $K=\RR$, the
irreducible examples are precisely the affine Weyl groups (see
Section~\ref{affgrp}). In the complex case $K=\CC$ there are two essentially
different cases. If $E/W$ is compact, the group $W$ is called
{\em crystallographic}. The noncrystallographic groups are now just the
complexifications of affine Weyl groups. The crystallographic reflection groups
have been classified by Popov \cite{Pop82} (see also Kaneko--Tokunaga--Yoshida
\cite{TY82,KTY82} for related results). As in the real case they are extensions
of a finite complex reflection group $W_0$ by an invariant lattice which is
generated by roots for $W_0$. \par
It turns out that presentations for these groups can be obtained as in the
case of affine Weyl groups by adding a further generating reflection
corresponding to a highest root in a root system for $W_0$ (see Malle
\cite{MaP}). As for the finite complex reflection groups in
Section~\ref{subsec:3.braid}, the braid group $B(W)$ of $W$ is defined as the
fundamental group of the hyperplane complement. For many of the irreducible
crystallographic complex reflection groups it is known that a presentation
for $B(W)$ can be obtained by ommitting the order relations from the
presentation of $W$ described above (see D\~ung \cite{Dung} for affine groups,
Malle \cite{MaP} for the complex case).
\end{para}

\begin{para}{Quaternionic reflection groups}
\label{subsec:5.quaternionic}
Cohen \cite{Co80} has obtained the classification of finite reflection groups
over the quaternions. This is closely related to finite linear groups over
$\CC$ generated by {\em bireflections}, that is, elements of finite order
which fix pointwise a subspace of codimension~2. Indeed, using the
identification of the quaternions as a certain ring of $2\times 2$-matrices
over the complex numbers, reflections over the quaternions become complex
bireflections. The primitive bireflection groups had been classified previously
by Huffman--Wales \cite{HuWa,Wal78}. One of the examples is a 3-dimensional
quaternionic representation of the double cover of the sporadic Hall--Janko
group $J_2$, see Wilson \cite{Wil86}.   \par
Presentations for these groups resembling the Coxeter presentations for
Weyl groups are given by Cohen \cite{Co91}. \par
In recent work on the McKay correspondence, quaternionic reflection groups
play an important r\^ole under the name of {\em symplectic reflection groups}
in the construction of so-called symplectic reflection algebras, see
for example Etingof--Ginzburg \cite{EtGi}.
\end{para}

\begin{para}{Reflection groups over finite fields}
\label{subsec:5.finite}
Many of the general results for complex reflection groups presented in
Section~1 are no longer true for reflection groups over fields of positive
characteristic. Most importantly, the ring of invariants of such a reflection
group is not necessarily a polynomial ring. Nevertheless we have the following
criterion due to Serre \cite[V.6, Ex.~8]{Bou} and Nakajima \cite{Nak79},
generalising Theorem~\ref{Gthm:steinberg}:

\begin{thm} [\rm Serre, Nakajima \cite{Nak79}]
 \label{Gthm:serre}
 Let $V$ be a finite-dimensional vector space over a field $K$ and
 $W\leq \GL(V)$ a finite group such that $K[V]^W$ is polynomial. Then
 the point-wise stabiliser of any subspace $U\leq V$ has polynomial ring of
 invariants $($and thus is generated by reflections$)$. 
 \end{thm}

The irreducible reflection groups over finite fields were classified by
Wagner \cite{Wag78,Wag80} and Zalesski{\u\i}--Sere\v zkin \cite{ZaSe80}, the
determination of transvection groups was completed by Kantor \cite{Kan79}. In
addition to the modular reductions of complex reflection groups, there arise
the infinite families of classical linear, symplectic, unitary and orthogonal
groups, as well as some further exceptional examples. For a complete list see
for example Kemper--Malle \cite[Section~1]{KM97}. \par
The results of Wagner \cite{Wag80} and Kantor \cite{Kan79} are actually
somewhat stronger, giving a classification of all indecomposable reflection
groups $W$ over finite fields of characteristic~$p$ for which the maximal
normal $p$-subgroup is contained in the intersection $W'\cap Z(W)$ of the
center with the derived group. \par

Using this classification, the {\em irreducible} reflection groups over finite
fields with polynomial ring of invariants could be determined, leading to
the following criterion:

\begin{thm} [\rm Kemper--Malle \cite{KM97}]
 \label{Gthm:kempermalle}
 Let $V$ be a finite-dimensional vector space over $K$, $W\leq \GL(V)$ a finite
 irreducible linear group. Then $K[V]^W$ is a polynomial ring if and only if
 $W$ is generated by reflections and the pointwise stabiliser in $W$ of any
 nontrivial subspace of $V$ has a polynomial ring of invariants.
\end{thm}

The list of groups satisfying this criterion can be found in
\cite[Thm.~7.2]{KM97}. That paper also contains some information on
indecomposable groups. \par
It is an open question whether at least the field of invariants of a reflection
group in positive characteristic is purely transcendental (by Kemper--Malle
\cite{KM99} the answer is positive in the irreducible case). \par

For further discussions of modular invariant theory of reflection groups see
also Derksen--Kemper \cite[3.7.4]{DK02}.
\end{para}

\begin{para}{$p$-adic reflection groups}
\label{subsec:5.padic}
Let $R$ be an integral domain, $L$ an $R$-lattice of finite rank, i.e., a
torsion-free finitely generated $R$-module, and $W$ a finite subgroup of
$\GL(L)$ generated by reflections. Again one can ask under which conditions
the invariants of $W$ on the symmetric algebra $R[L]$ of the dual $L^*$ are a
graded polynomial ring. In the case of Weyl groups Demazure shows the
following extension of Theorem~\ref{Gthm:shtochev}:

\begin{thm}  [\rm Demazure \cite{Dem}]
 \label{Gthm:demazure}
 Let $W$ be a Weyl group, $L$ the root lattice of $W$, and $R$ a ring in
 which all torsion primes of $W$ are invertible. Then the invariants of $W$
 on $R[L]$ are a graded polynomial algebra, and $R[L]$ is
 a free graded module over $R[L]^W$.
\end{thm}

In the case of general lattices for reflection groups, the following example
may be instructive: Let $W=\fS_3$ the symmetric group of degree~3. Then the
weight lattice $L$ of $S_3$, considered as $\ZZ_3$-lattice, yields a faithful
reflection representation of $\fS_3$ with the following property:
$\ZZ_3[L]^{\fS_3}$ is not polynomial, while both the reflection
representations over the quotient field $\QQ_3$ and over
the residue field $\FF_3$ have polynomial invariants, the first with
generators in degrees~2 and~3, the second with generators in degrees~1 and~6.

The list of all irreducible $p$-adic reflection groups, that is, reflection
groups over the field of $p$-adic numbers $\QQ_p$, was obtained by
Clark--Ewing \cite{ClEw}, building on the Theorem of Shephard--Todd. We
reproduce it in Table~\ref{Gtab4}.

\begin{table}[htbp] \caption{$p$-adic reflection groups}
\label{Gtab4}
$\begin{array}{c}
\begin{array}{ll}
\hline W     &\text{conditions} \\ \hline
G(de,e,n)\;((d,n)\ne(1,2))  &p\equiv1\pmod{de}; \text{ none when }de=2\\
G(e,e,2)\; (e{\geq} 3) &p\equiv\pm1\pmod e; \text{ none when }e=3,4,6\\
\fS_{n{+}1}& \text{---}  \\
\hline
\end{array}\\ \\
\begin{array}{ll}
\hline W     &\text{conditions} \\ \hline
 G_4    &p\equiv1\pmod3 \cr
 G_5    &p\equiv1\pmod3 \cr
 G_6    &p\equiv1\pmod{12} \cr
 G_7    &p\equiv1\pmod{12} \cr
 G_8    &p\equiv1\pmod4 \cr
 G_9    &p\equiv1\pmod8 \cr
 G_{10} &p\equiv1\pmod{12} \cr
 G_{11} &p\equiv1\pmod{24} \cr
 G_{12} &p\equiv1,3\pmod8 \cr
 G_{13} &p\equiv1\pmod8 \cr
 G_{14} &p\equiv1,19\pmod{24} \cr
 G_{15} &p\equiv1\pmod{24} \cr
 G_{16} &p\equiv1\pmod5 \cr
 G_{17} &p\equiv1\pmod{20} \cr
 G_{18} &p\equiv1\pmod{15} \cr
 G_{19} &p\equiv1\pmod{60} \cr
 G_{20} &p\equiv1,4\pmod{15} \cr
 \hline
\end{array}\qquad
\begin{array}{lll@{\hspace{0pt}}c}
\hline W     &\text{conditions} \\ \hline
 G_{21} &p\equiv1,49\pmod{60} \cr
 G_{22} &p\equiv1,9\pmod{20} \cr
 G_{23} &p\equiv1,4\pmod5 \cr
 G_{24} &p\equiv1,2,4\pmod7 \cr
 G_{25} &p\equiv1\pmod3  \cr
 G_{26} &p\equiv1\pmod3 \cr
 G_{27} &p\equiv1,4\pmod{15} \cr
 G_{28} &\text{---}  \cr
 G_{29} &p\equiv1\pmod4  \cr
 G_{30} &p\equiv1,4\pmod5 \cr
 G_{31} &p\equiv1\pmod4  \cr
 G_{32} &p\equiv1\pmod3 \cr
 G_{33} &p\equiv1\pmod3 \cr
 G_{34} &p\equiv1\pmod3 \cr
 G_{35} &\text{---} \cr
 G_{36} &\text{---} \cr
 G_{37} &\text{---} \cr
 \hline
\end{array}
\end{array}
$
\end{table}

Using a case-by-case argument based on the Clark--Ewing classification and
his own classification of $p$-adic lattices for reflection groups, Notbohm
\cite{Not99} was able to determine all finite reflection groups $W$ over the
ring of $p$-adic integers $\ZZ_p$, $p>2$, with polynomial ring of invariants.
This was subsequently extended by Andersen--Grodal--M\o ller--Viruel
\cite{AGMV} to include the case $p=2$. 
\end{para}

\begin{para}{$p$-compact groups}
The $p$-adic reflection groups play an important r\^ole in the theory of
so-called $p$-compact groups, which constitute a homotopy theoretic
analogue of compact Lie groups.  By definition, a {\em $p$-compact group} is a
$p$-complete topological space $BX$ such that the homology $H_*(X;\FF_p)$
of the loop space $X=\Omega BX$ is finite. Examples for $p$-compact groups are
$p$-completions of classifying spaces of compact Lie groups. Further examples
were constructed by Clark--Ewing \cite{ClEw}, Aguad\'e \cite{Agu},
Dwyer--Wilkerson \cite{DwWi93} and Notbohm \cite{Not98}. To each $p$-compact
group $X$ Dwyer--Wilkerson \cite{DwWi94} associate a maximal torus (unique up
to conjugacy) together with a \lq Weyl group\rq, which comes equipped with a
representation as a reflection group over the $p$-adic integers $\ZZ_p$, which
is faithful if $X$ is connected. Conversely, by a theorem of
Andersen et al. \cite{AGMV} a connected $p$-compact group,
for $p>2$, is determined up to isomorphism by its Weyl group data, that is, by
its Weyl group in a reflection representation on a $\ZZ_p$-lattice.  \par
It has been shown that at least for $p>2$ all $p$-adic reflection groups (as
classified by Clark--Ewing) and all their $\ZZ_p$-reflection representations
arise in that way (see Andersen et al. \cite{AGMV}, Notbohm \cite{Not98} and
also Adams--Wilkerson \cite{AdWi}, Aguad\'e \cite{Agu}).
\end{para}

\bigskip



\begin{thebibliography}{131}
 
\bibitem{AdWi}
{\sc J. F. Adams and C. W. Wilkerson}, Finite $H$-spaces and algebras over the
  Steenrod algebra. Ann. of Math. {\bf111} (1980), 95--143.

\bibitem{Agu}
{\sc J. Aguad\'e}, Constructing modular classifying spaces. Israel J. Math.
  {\bf66} (1989), 23--40.

\bibitem{AGMV}
{\sc K. Andersen, J. Grodal, J. M\o ller and A. Viruel}, The classification of
  $p$-compact groups for $p$ odd. In preparation.

\bibitem{Ar}
{\sc E.~Artin}, Theory of braids. Ann. of Math. {\bf48} (1947), 101--126.

\bibitem{Ban}
{\sc E.~Bannai}, Fundamental groups of the spaces of regular orbits of the
  finite unitary reflection groups of dimension $2$. J. Math. Soc. Japan
  {\bf28} (1976), 447--454.

\bibitem{Ben76}
{\sc M.~Benard}, Schur indices and splitting fields of the unitary reflection
  groups. J. Algebra {\bf38} (1976), 318--342.

\bibitem{Ben}
{\sc D.~Benson}, {\em Polynomial invariants of finite groups}. London
  Math. Soc. Lecture Note Series {\bf190}, Cambridge University Press,
  Cambridge, 1993. 

\bibitem{BeGr}
{\sc C.~T.~Benson and L.~C. Grove}, {\em Finite reflection groups 
  (2nd edition)}. Graduate Texts in Mathematics, vol.~99, Springer-Verlag, 1985.

\bibitem{Bes}
{\sc D.~Bessis}, Sur le corps de d\'efinition d'un groupe de r\'eflexions
  complexe. Comm. in Algebra {\bf25} (1997), 2703--2716.

\bibitem{Bes2}
{\sc D.~Bessis}, Groupes des tresses et \'el\'ements r\'eguliers. J. reine
  angew. Math. {\bf518} (2000), 1--40.

\bibitem{Bes3}
{\sc D.~Bessis}, Zariski theorems and diagrams for braid groups. Invent. Math.
  {\bf145} (2001),487--507.
  
\bibitem{BDM}
{\sc D.~Bessis, F.~Digne and J.~Michel}, Springer theory in braid groups and
  the Birman-Ko-Lee monoid. Pacific J. Math. {\bf 205} (2002), 287--310.

\bibitem{BeLu}
{\sc W. Beynon and G. Lusztig}, Some numerical results on the characters of
  exceptional Weyl groups. Math. Proc. Cambridge Philos. Soc. {\bf84} (1978),
  417--426.

\bibitem{Big1}
{\sc S.~Bigelow}, The Burau representation is not faithful for $n=5$.
  Geometry and Topology {\bf 3} (1999), 397--404.

\bibitem{Big2}
{\sc S.~Bigelow}, Braid groups are linear. J. Amer. Math. Soc. {\bf 14}
  (2001), 471--486.

\bibitem{Birman74}
{\sc J.~S.~Birman}, {\em Braids, links and mapping class groups}.
   Annals of Math. Stud., vol. 84, Princeton University Press, 1974.

\bibitem{Birman93}
{\sc J.~S.~Birman}, New points of view in knot theory. Bull. Amer. Math. Soc.
  {\bf28} (1993), 253--287.

\bibitem{BKL}
{\sc J.~Birman, K.~H.~Ko and S.~J.~Lee}, A new approach to the word and
  conjugacy problems in the braid groups. Adv. Math. {\bf139} (1998), 322--353. 

\bibitem{Bjo}
{\sc A.~Bj\"orner}, Orderings of Coxeter groups. {\em In: Combinatorics and
  Algebra}, Contemporary Math., vol. 34, Amer. Math. Soc., 1984, pp.~175--195.

\bibitem{BlLe}
{\sc J.~Blair and G.~I.~Lehrer}, Cohomology actions and centralisers in unitary
  reflection groups. Proc. London Math. Soc. {\bf83} (2001), 582--604.

\bibitem{BlGeKi97}
{\sc F.~Bleher, M.~Geck and W.~Kimmerle}, Automorphisms of integral group 
  rings of finite Coxeter groups and Iwahori--Hecke algebras. J. Algebra
  {\bf197} (1997), 615--655.

\bibitem{Bor98}
{\sc R.~E.~Borcherds}, Coxeter groups, Lorentzian lattices, and $K3$
  surfaces. Internat. Math. Res. Notices {\bf19} (1998), 1011--1031.

\bibitem{Borel}
{\sc A.~Borel}, {\em Linear algebraic groups. Second enlarged edition}.
  Graduate Texts in Mathematics vol.~126. Springer Verlag, 
  Berlin--Heidelberg--New York, 1991.

\bibitem{Boro}
{\sc A.~V.~Borovik, I.~M.~Gelfand and N.~White}, Coxeter matroid polytopes.
  Ann. Comb. {\bf 1} (1997), 123--134.

\bibitem{Bou}
{\sc N.~Bourbaki}, {\em Groupes et alg{\`e}bres de {L}ie, chap. 4, 5 et 6}.
 Hermann, Paris, 1968.

\bibitem{BreMa1}
{\sc K. Bremke and G. Malle}, Reduced words and a length function for
  $G(e,1,n)$. Indag. Mathem. {\bf8} (1997), 453--469.

\bibitem{BreMa2}
{\sc K. Bremke and G. Malle}, Root systems and length functions.
  Geom. Dedicata {\bf72} (1998), 83--97.

\bibitem{Br71}
{\sc E.~Brieskorn}, Die Fundamentalgruppe des Raumes der regul\"aren Orbits
  einer endlichen komplexen Spiegelungsgruppe. Invent. Math. {\bf12} (1971),
  37--61

\bibitem{Br73}
{\sc E.~Brieskorn}, Sur les groupes de tresses [d'apr\`es V.I.~Arnold].
  {\em In: S\'eminaire Bourbaki, 24\`eme ann\'ee, 1971/72}, Lecture Notes in
  Math. {\bf317}, Springer, Berlin, 1973.

\bibitem{BrSa}
{\sc E.~Brieskorn and K.~Saito}, Artin-Gruppen und Coxeter-Gruppen.
  Invent. Math. {\bf17} (1972), 245--271.

\bibitem{BriHow93}
{\sc B.~Brink and R.~B.~Howlett}, A finiteness property and an automatic
  structure for Coxeter groups. Math. Ann. {\bf296} (1993), 179--190.

\bibitem{BriHow99}
{\sc B.~Brink and R.~B.~Howlett}, Normalizers of parabolic subgroups in
  Coxeter groups. Invent. Math. {\bf136} (1999), 323--351.

\bibitem{Br01}
{\sc M.~Brou\'e}, Reflection groups, braid groups, Hecke algebras, finite
  reductive groups. {\em In: Current developments in mathematics}, Int. Press,
  Somerville, 2001, pp.~1--107.
  
\bibitem{BrMa}
{\sc M.~Brou\'e and G.~Malle}, Th\'eor\`emes de Sylow g\'en\'eriques pour les
 groupes r\'eductifs sur les corps finis. Math. Ann. {\bf292} (1992), 241--262.

\bibitem{BrMa93}
{\sc M.~Brou\'e and G.~Malle}, Zyklotomische Heckealgebren. 
  Ast\'erisque {\bf212} (1993), 119--189.

\bibitem{BrMa3}
{\sc M.~Brou\'e and G.~Malle}, Generalized Harish-Chandra theory.
  {\em In: Representations of reductive groups}, Cambridge
  University Press, Cambridge, 1998, pp.~85--103.

\bibitem{BMM99}
{\sc M.~Brou\'e, G.~Malle and J.~Michel}, Towards spetses I.
  Transform. Groups {\bf 4} (1999), 157--218.

\bibitem{BMR}
{\sc M. Brou\'e, G. Malle and R.~Rouquier}, Complex reflection groups,
  braid groups, Hecke algebras. J. reine angew. Math. {\bf 500} (1998),
  127--190.

\bibitem{BrMi} 
{\sc M.~Brou\'{e} and J.~Michel}, Sur certains \'el\'ements r\'eguliers 
  des groupes de Weyl et les vari\'et\'es de Deligne-Lusztig associ\'ees. 
  {\em In: Finite reductive groups, related structures 
  and representations}, Progress in Math. {\bf 141}, Birkh\"auser, 1997, 
  pp. 73--139.

\bibitem{BuZi}
{\sc G.~Burde and H.~Zieschang}, {\em Knots}. De Gruyter Studies in 
  Math., vol.~5; Walter de Gruyter \& Co., Berlin, 1985.

\bibitem{Carter72}
{\sc R.~W.~Carter}, Conjugacy classes in the Weyl group. Compositio Math.
  {\bf25} (1972), 1--59.

\bibitem{Carter1}
{\sc R.~W.~Carter}, {\em Simple groups of Lie type}. Wiley, New York, 1972;
  reprinted 1989 as Wiley Classics Library Edition.

\bibitem{Carter2}
{\sc R.~W.~Carter}, {\em Finite groups of Lie type: Conjugacy classes and
  complex characters}. Wiley, New York, 1985; reprinted 1993 as Wiley
  Classics Library Edition.

\bibitem{Casselman94}
{\sc W.~A.~Casselman}, Machine calculations in Weyl groups. Invent. Math.
  {\bf116} (1994), 95--108.

\bibitem{Chev0}
{\sc C.~Chevalley}, Sur certains groupes simples. T\^ohoku Math. J. (2)
  {\bf 7} (1955), 14--66.

\bibitem{Ch}
{\sc C.~Chevalley}, Invariants of finite groups generated by reflections.
  Amer. J. Math. {\bf 77} (1955), 778--782.

\bibitem{Chev}
{\sc C.~Chevalley}, {\em Classification des groupes de Lie alg\'ebriques}.
  S\'eminaire \'Ecole Normale Sup\'erieure, Mimeographed Notes, Paris,
  1956--1958.

\bibitem{Chev1}
{\sc C.~Chevalley}, Sur les d\'ecompositions cellulaires des espaces $G/B$.
  With a foreword by A.~Borel. Proc. Symp. Pure Math., vol.~56,
  Amer. Math. Soc., Providence, RI, 1994, pp.~1--23.

\bibitem{ClEw}
{\sc A. Clark and J. Ewing}, The realization of polynomial algebras as
  cohomology rings. Pacific J. Math. {\bf50} (1974), 425--434.

\bibitem{Co}
{\sc A.~M.~Cohen}, Finite complex reflection groups. Ann. scient. \'Ec. Norm.
  Sup. {\bf9} (1976), 379--436. {\em Erratum:} ibid. {\bf11} (1978), 613.

\bibitem{Co80}
{\sc A.~M.~Cohen}, Finite quaternionic reflection groups. J. Algebra {\bf64}
  (1980), 293--324.

\bibitem{Co91}
{\sc A.~M.~Cohen},  Presentations for certain finite quaternionic reflection
  groups. {\em In: Advances in finite geometries and designs},
  Oxford Univ. Press, New York, 1991, pp.~69--79.

\bibitem{CoWa}
{\sc A.~M.~Cohen and D.~Wales}, Linearity of Artin groups of finite
  type. preprint, ArXiv, math.GR/0010204.

\bibitem{Cox34}
{\sc H.~S.~M.~Coxeter}, Discrete groups generated by reflections.
  Annals of Math. {\bf 35} (1934), 588--621.

\bibitem{Cox35}
{\sc H.~S.~M.~Coxeter}, The complete enumeration of finite groups of 
  the form $R_i^2=(R_iR_j)^{k_{ij}}=1$. J. London Math. Soc. {\bf 10} (1935),
  21--25. 

\bibitem{Cox57}
{\sc H.~S.~M.~Coxeter}, Groups generated by unitary reflections of period two.
  Canad. J. Math. {\bf9} (1957), 243--272.

\bibitem{CroFox}
{\sc R.~H.~Crowell and R.~H.~Fox}, {\em Introduction to knot theory}. 
  Grad. Texts in Math., vol.~57, Springer-Verlag, 1977; reprint of the 1963 
  original.

\bibitem{CR2}
{\sc C.~W.~Curtis and I.~Reiner}, {\em Methods of representation theory, 
  vol.~{II}}. Wiley, New York, 1987; reprinted 1994 as Wiley Classics 
  Library Edition.

\bibitem{Dehor}
{\sc P.~Dehornoy}, {\em Braids and self-distributivity}. Progress
in Math., vol.~192, Birkh\"auser Verlag, Basel, 2000.

\bibitem{DehPar}
{\sc P.~Dehornoy and L.~Paris}, Gaussian groups and Garside groups, two
  generalizations of Artin groups. Proc. London Math. Soc. {\bf 79} (1999),
  569--604.

\bibitem{De}
{\sc P.~Deligne}, Les immeubles des groupes de tresses g\'en\'eralis\'es.
  Invent. Math. {\bf17} (1972), 273--302

\bibitem{Deligne97}
{\sc P.~Deligne}, Action du groupe de tresses sur une cat\'egorie. Invent.
  Math. {\bf128} (1997), 159--175.

\bibitem{Dem}
{\sc M. Demazure}, Invariants sym\'etriques entiers des groupes de Weyl et
  torsion. Invent. Math. {\bf21} (1973), 287--301. 

\bibitem{DeLo95}
{\sc J.~Denef and F.~Loeser},  Regular elements and monodromy of discriminants
  of finite reflection groups. Indag. Math. {\bf6} (1995), 129--143.

\bibitem{Deodhar}
{\sc V.~V.~Deodhar}, Some characterizations of Bruhat ordering on a 
  Coxeter group and determination of the relative M\"obius function.
  Invent. Math. {\bf39} (1977), 187--198.

\bibitem{Deodhar82}
{\sc V.~V.~Deodhar}, On the root system of a Coxeter group. Comm. Algebra
  {\bf10} (1982), 611--630.

\bibitem{Deo89}
{\sc V.~V.~Deodhar}, A note on subgroups generated by reflections in Coxeter
  groups. Arch. Math. (Basel) {\bf53} (1989), 543--546.

\bibitem{DK02}
{\sc H.~Derksen and G.~Kemper}, {\em Computational Invariant Theory}.
  Encyclopaedia of Mathematical Sciences {\bf130}, Springer-Verlag,
  Berlin Heidelberg New York, 2002.

\bibitem{Digne}
{\sc F.~Digne}, On the linearity of Artin braid groups. J. Algebra (to appear).

\bibitem{DJ1} 
{\sc R.~Dipper and G.~D.~James}, The $q$-Schur algebra. Proc.\ London 
  Math.\ Soc. {\bf 59} (1989), 23--50.

\bibitem{Cloux:Trans}
{\sc F.~du~Cloux}, A transducer approach to Coxeter groups. J. Symbolic
  Comput. {\bf27} (1999), 1--14.

\bibitem{Dung}
{\sc N. V. D\~ung, }, The fundamental groups of the spaces
  of regular orbits of the affine Weyl groups. Topology {\bf22} (1983),
  425--435.

\bibitem{DwWi93}
{\sc W. G. Dwyer and C. W. Wilkerson}, A new finite loop space at the prime
  two. J. Amer. Math. Soc. {\bf6} (1993), 37--64.

\bibitem{DwWi94}
{\sc W. G. Dwyer and C. W. Wilkerson}, Homotopy fixed-point methods for Lie
  groups and finite loop spaces. Ann. of Math. {\bf139} (1994), 395--442. 

\bibitem{DyerPhd}
{\sc M.~Dyer}, {\em Hecke algebras and reflections in Coxeter groups}. 
  Ph. {D}.  thesis, University of Sydney, 1987.

\bibitem{Dye90}
{\sc M.~Dyer}, {Reflection subgroups of {Coxeter} systems}. J. Algebra
  {\bf135} (1990), 57--73.

\bibitem{Erd}
{\sc K.~Erdmann}, Decomposition numbers for symmetric groups and
  composition factors of Weyl modules. J. Algebra {\bf 180} (1996), 316--320.

\bibitem{EtGi}
{\sc P. Etingof, V. Ginzburg}, Symplectic reflection algebras, Calogero-Moser
  space, and deformed Harish-Chandra homomorphism. Invent. Math. {\bf147}
  (2002), 243--348.

\bibitem{Fla78}
{\sc L.~Flatto}, Invariants of finite reflection groups. Enseign. Math. {\bf24}
  (1978), 237--292.

\bibitem{FlJa93}
{\sc P.~Fleischmann and I.~Janiszczak}, Combinatorics and Poincar\'e
  polynomials of hyperplane complements for exceptional Weyl groups. J. Combin.
  Theory Ser. A {\bf63} (1993), 257--274.

\bibitem{FoSe}
{\sc P.~Fong and G.~M.~Seitz}, Groups with a $BN$-pair of rank~$2$, I.
  Invent. Math. {\bf 21} (1973), 1--57; II, ibid. {\bf 24} (1974),
  191--239.

\bibitem{FoNe}
{\sc R.~H.~Fox and L.~Neuwirth}, The braid groups. Math. Scand. {\bf10} (1962),
  119--126.

\bibitem{HOMFLY}
{\sc D.~Freyd, D.~Yetter, J.~Hoste, W.~B.~R.~Lickorish, K.~Millet and
  A.~Ocneanu}, A new polynomial invariant of knots and links. Bull. Amer.
  Math. Soc. {\bf 12} (1985), 239--246.

\bibitem{Frob}
{\sc F.~G.~Frobenius}, \"Uber die Charaktere der symmetrischen Gruppe.
  Sitzungsber. Preuss. Akad. Wiss. Berlin (1900), 516--534.

\bibitem{Fult}
{\sc W.~Fulton}, {\em Young tableaux}. London Math. Soc. Student Texts {\bf35}, 
  Cambridge University Press, Cambridge, 1997.

\bibitem{Garside69}
{\sc F.~A.~Garside}, The braid group and other groups. Quart. J. Math. Oxford
  {\bf20} (1969), 235--254.

\bibitem{GeckBanach}
{\sc M.~Geck}, {Trace functions on {H}ecke algebras}. Banach Center Publ.,
  Polish Acad. Sci. {\bf42} (1998), 87--109.

\bibitem{mykl}
{\sc M.~Geck}, Kazhdan--Lusztig cells and decomposition numbers.
  Represent. Theory {\bf 2} (1998), 264--277 (electronic).

\bibitem{chev}
{\sc M.~Geck, G.~Hi{\ss}, F.~L\"ubeck, G.~Malle and G.~Pfeiffer},
{\sf CHEVIE} --- A system for computing and processing generic character
  tables. Appl. Algebra Engrg. Comm. Comput. {\bf 7} (1996), 175--210.

\bibitem{GeKim}
{\sc M.~Geck and S.~Kim}, Bases for the Bruhat--Chevalley order on all
  finite Coxeter groups. J.\ Algebra {\bf 197} (1997), 278--310.

\bibitem{GeKiPf}
{\sc M.~Geck, S.~Kim and G.~Pfeiffer}, Minimal length elements in
  twisted conjugacy classes of finite Coxeter groups. J. Algebra {\bf 229}
  (2000), 570--600.

\bibitem{GeLa}
{\sc M.~Geck and S.~Lambropoulou}, Markov traces and knot invariants
  related to Iwahori--Hecke algebras of type $B$. J. reine angew. Math.
  {\bf 482} (1997), 191--213.

\bibitem{GeMi}
{\sc M.~Geck and J.~Michel}, "Good" elements of finite Coxeter groups and
  representations of Iwahori-Hecke algebras. Proc. London Math. Soc. {\bf74}
  (1997), 275--305. 

\bibitem{GePf93}
{\sc M.~Geck and G.~Pfeiffer}, {On the irreducible characters of Hecke
  algebras}. Adv. Math. {\bf102} (1993), 79--94.

\bibitem{GePf}
{\sc M.~Geck and G.~Pfeiffer}, {\em Characters of finite Coxeter groups and
  Iwahori-Hecke algebras}. Oxford University Press, Oxford, 2000.

\bibitem{GoLySo}
{\sc D.~Gorenstein, R.~Lyons and R.~Solomon}, {\em The classification of the
  finite simple groups}. Math. Surveys and Monographs vol.~40, no.~1,
  Amer. Math. Soc., 1994.

\bibitem{Green1}
{\sc J.~A.~Green}, {\em Polynomial representations of  $\mbox{GL}_n$}.
  Lecture Notes in Math., vol.~830, Springer--Verlag, Berlin--New York, 1980.

\bibitem{HeKaSe}
{\sc C.~Hering, W.~M.~Kantor and G.~M.~Seitz}, Finite groups with a split
  $BN$-pair of rank~$1$. J.~Algebra {\bf 20} (1972), 435--475.

\bibitem{Hiller82} 
{\sc H.~Hiller}, {\em Geometry of Coxeter groups}. Res. Notes Mathematics,
  vol.~54, Pitman, Boston, Mass., 1982.
  
\bibitem{Howlett80}
{\sc R.~B. Howlett}, {Normalizers of parabolic subgroups of reflection 
  groups}. J. London Math. Soc. (2) {\bf21} (1980), 62--80.

\bibitem{HuWa}
{\sc W. Huffman and D. Wales}, Linear groups containing an involution with two
  eigenvalues $-1$. J. Algebra {\bf45} (1977), 465--515. 

\bibitem{Hu83}
{\sc M.~Hughes}, Representations of complex imprimitive reflection groups.
  Math. Proc. Cambridge Philos. Soc. {\bf94} (1983), 425--436.

\bibitem{HuMo}
{\sc M.~Hughes and A. Morris}, Root systems for two dimensional complex
  reflection groups. S\'em. Lothar. Combin. {\bf45}, 18 pp. (electronic).

\bibitem{Humphreys0} 
{\sc J.~E.~Humphreys}, {\em Introduction to Lie algebras and representation
  theory}. Graduate Texts in Mathematics, vol.~9, Springer-Verlag, 1972.

\bibitem{Humphreys2}
{\sc J.~E.~Humphreys}, {\em Reflection groups and Coxeter groups}. 
  Cambridge Stud. Adv. Math., vol.~29, Cambridge University Press, 1990.

\bibitem{Humphreys1}
{\sc J.~E.~Humphreys}, {\em Linear algebraic groups, second edition}.
  Graduate Texts in Mathematics vol.~21. Springer Verlag,
  Berlin--Heidelberg--New York, 1991.

\bibitem{Hum94}
{\sc J.~F.~Humphreys}, Character tables for the primitive finite unitary
  reflection groups. Comm. Algebra {\bf22} (1994), 5777--5802.

\bibitem{Iwa64}
{\sc N.~Iwahori}, {On the structure of a Hecke ring of a Chevalley group
  over a finite field}. J. Fac. Sci. Univ. Tokyo {\bf10} (1964), 215--236.

\bibitem{iwamat}
{\sc N.~Iwahori and H.~Matsumoto}, On some Bruhat decomposition and the
  structure of the Hecke ring of $\mathfrak{p}$-adic Chevalley groups.
  Publ. Math. I.H.E.S. {\bf 25} (1965), 5--48.

\bibitem{James1}
{\sc G.~D.~James}, {\em The representation theory of the symmetric groups}.
  Lecture Notes in Math., vol.~682, Springer--Verlag, Berlin, 1978.

\bibitem{JaKe}
{\sc G.~D.~James and A.~Kerber}, {\em The representation theory of the 
  symmetric group}. Encyclopedia of Mathematics, vol.~16, Addison-Wesley, 1981.

\bibitem{Jones87}
{\sc V.~F.~R.~Jones}, {{H}ecke algebra representations of braid groups and 
  link polynomials}. Ann. of Math. {\bf126} (1987), 335--388.

\bibitem{Kac}
{\sc V.~G.~Kac}, {\em Infinite dimensional {L}ie algebras}. Cambridge 
  University Press, 1985.

\bibitem{Kane}
{\sc R. Kane}, {\em Reflection groups and invariant theory}. CMS Books in
  Mathematics {\bf5}, Springer, New York, 2001.

\bibitem{KTY82}
{\sc J. Kaneko, S. Tokunaga and M. Yoshida}, Complex crystallographic groups.
  II. J. Math. Soc. Japan {\bf34} (1982), 595--605.

\bibitem{Kan79}
{\sc W. M. Kantor}, Subgroups of classical groups generated by long root
  elements. Trans. Amer. Math. Soc. {\bf248} (1979), 347--379. 

\bibitem{KaSe}
{\sc W.~M.~Kantor and G.~M.~Seitz}, Finite groups with a split
  $BN$-pair of rank~$1$. II. J.~Algebra {\bf 20} (1972), 476--494.

\bibitem{KM97}
{\sc G. Kemper and G. Malle}, The finite irreducible linear groups with
  polynomial ring of invariants. Transform. Groups {\bf2} (1997), 57--89.

\bibitem{KM99}
{\sc G. Kemper and G. Malle}, Invariant fields of finite irreducible
  reflection groups. Math. Ann. {\bf315} (1999), 569--586.

\bibitem{Kim98}
{\sc S.~Kim}, {Mots de longueur maximale dans une classe de conjugaison d'un
  groupe sym\'etrique, vu comme groupe de Coxeter}. C. R. Acad. Sci. Paris
  {\bf327} (1998), 617--622.

\bibitem{Kram}
{\sc D.~Krammer}, {\em The conjugacy problem for Coxeter groups}. Ph. D.
  thesis, University of Utrecht, 1994.

\bibitem{Kram2}
{\sc D.~Krammer}, Braid groups are linear. Annals of Math. {\bf 155} (2002),
  131--156.

\bibitem{Ku77}
{\sc S.~Kusuoka}, On a conjecture of L. Solomon. J. Fac. Sci. Univ. Tokyo
  Sect. IA Math. {\bf24} (1977), 645--655.

\bibitem{Lamb1}
{\sc S.~Lambropoulou}, Solid torus links and Hecke algebras of $B$-type.
  {\em In: Proceedings of the Conference on Quantum Topology} (Manhattan,
  KS, 1993), World Sci. Publishing, River Edge, NJ, 1994, pp.~225--245.

\bibitem{Lamb2}
{\sc S.~Lambropoulou}, Knot theory related to generalized and cyclotomic
  Hecke algebras of type $B$. J. Knot Theory Ramifications {\bf 8} (1999),
  621--658.

\bibitem{LaSch} 
{\sc A.~Lascoux and M.~P.~Sch{\"u}tzenberger}, Treillis et bases des groupes
  de Coxeter. Electron. J. Combin. {\bf 3} (1995), 35 pp.

\bibitem{Le87b}
{\sc G.~I.~Lehrer}, On the Poincar\'e series associated with Coxeter group
  actions on complements of hyperplanes. J. London Math. Soc. {\bf36} (1987),
  275--294.

\bibitem{Le95}
{\sc G.~I.~Lehrer}, Poincar\'e polynomials for unitary reflection groups.
  Invent. Math. {\bf120} (1995), 411--425.
 
\bibitem{LeMi}
{\sc G.~I.~Lehrer and J.~Michel}, Invariant theory and eigenspaces for
 unitary reflection groups. Preprint 2003. 

\bibitem{LeSp}
{\sc G.~I.~Lehrer and T.~A.~Springer}, Intersection multiplicities and
 reflection subquotients of unitary reflection groups I. {\em In: Geometric
 group theory down under} (Canberra, 1996), de Gruyter, Berlin, 1999,
 pp.~181--193. 

\bibitem{LeSp2}
{\sc G.~I.~Lehrer and T.~A.~Springer}, Reflection subquotients of unitary
 reflection groups. Canad. J. Math. {\bf51} (1999), 1175--1193. 

\bibitem{LoPa}
{\sc D.~D.~Long and M.~Paton}, The Burau representation is not faithful 
  for $n\geq 6$. Topology {\bf 32} (1993), 439--447.

\bibitem{Lu77}
{\sc G.~Lusztig}, Irreducible representations of finite classical groups. 
  Invent. Math. {\bf43} (1977), 125--175. 

\bibitem{LuSp79}
{\sc G.~Lusztig and N.~Spaltenstein}, Induced unipotent classes. 
  J. London Math. Soc. {\bf19} (1979), 41--52. 

\bibitem{Mac72}
{\sc I.~G.~Macdonald}, Some irreducible representations of Weyl groups. 
  Bull. London Math. Soc. {\bf4} (1972), 148--150. 

\bibitem{Mac95}
{\sc I.~G.~Macdonald}, {\em Symmetric functions and Hall polynomials}. 
  2nd edition, Oxford University Press, New York, 1995.

\bibitem{MaU}
{\sc G. Malle}, Unipotente Grade imprimitiver komplexer Spiegelungsgruppen.
  J. Algebra {\bf 177} (1995), 768--826.

\bibitem{MaP}
{\sc G. Malle}, Presentations for crystallographic complex reflection groups.
  Transform. Groups {\bf1} (1996), 259--277.

\bibitem{MaR}
{\sc G.~Malle}, On the rationality and fake degrees of characters of
  cyclotomic algebras. J. Math. Sci. Univ. Tokyo {\bf6} (1999), 647--677.

\bibitem{Matsum}
{\sc H.~Matsumoto}, G\'en\'erateurs et relations des groupes de Weyl
  g\'en\'eralis\'ees. C. R. Acad. Sci. Paris {\bf 258} (1964), 3419--3422.

\bibitem{Mi99}
{\sc J.~Michel}, A note on words in braid monoids. J. Algebra {\bf215} (1999),
 366--377.    

\bibitem{Moody93}
{\sc J.~A. Moody}, {The faithfulness question for the Burau representation}.
  Proc. Amer. Math. Soc. {\bf119} (1993), 671--679.

\bibitem{Nak79}
{\sc H. Nakajima}, Invariants of finite groups generated by pseudoreflections
 in positive characteristic. Tsukuba J. Math. {\bf3} (1979), 109--122. 

\bibitem{Nak83}
{\sc T.~Nakamura}, A note on the $K(\pi ,\,1)$ property of the orbit space of
 the unitary reflection group $G(m,l,n)$. Sci. Papers College Arts Sci.
 Univ. Tokyo {\bf33} (1983), 1--6.

\bibitem{Nar}
{\sc I.~Naruki}, The fundamental group of the complement for Klein's
 arrangement of twenty-one lines. Topology Appl. {\bf34} (1990), 167--181.

\bibitem{Neb}
{\sc G. Nebe}, The root lattices of the complex reflection groups.
 J. Group Theory {\bf2} (1999), 15--38.

\bibitem{Not96}
{\sc D. Notbohm}, $p$-adic lattices of pseudo reflection groups. {\em In:
  Algebraic topology: new trends in localization and periodicity}, Progr.
  Math. {\bf136}, Birkh\"auser, Basel, 1996, pp.~337--352.

\bibitem{Not98}
{\sc D. Notbohm}, Topological realization of a family of pseudoreflection
  groups. Fund. Math. {\bf155} (1998), 1--31.

\bibitem{Not99}
{\sc D. Notbohm}, For which pseudo-reflection groups are the $p$-adic
  polynomial invariants again a polynomial algebra? J. Algebra {\bf214} (1999),
  553--570. {\em Erratum:} ibid. {\bf218} (1999), 286--287.

\bibitem{Op95}
{\sc E. M. Opdam},  A remark on the irreducible characters and fake degrees
 of finite real reflection groups. Invent. Math. {\bf120} (1995), 447--454.

\bibitem{Op00}
{\sc E. M. Opdam}, {\em Lecture notes on Dunkl operators for real and complex
  reflection groups}. MSJ Memoirs {\bf8}, Math. Soc. Japan, Tokyo, 2000.

\bibitem{OrSo80a}
{\sc P.~Orlik and L.~Solomon}, Combinatorics and topology of complements
 of hyperplanes. Invent. Math. {\bf56} (1980), 167--189.

\bibitem{OrSo80b}
{\sc P.~Orlik and L.~Solomon}, Unitary reflection groups and cohomology.
  Invent. Math. {\bf59} (1980), 77--94.

\bibitem{OrSo83}
{\sc P.~Orlik and L.~Solomon}, A character formula for the unitary group over
 a finite field. J. Algebra {\bf84} (1983), 136--141.

\bibitem{OrSo85}
{\sc P.~Orlik and L.~Solomon}, Arrangements in unitary and orthogonal geometry
  over finite fields. J. Combin. Theory Ser. A {\bf38} (1985), 217--229.

\bibitem{OrSo88}
{\sc P.~Orlik and L.~Solomon}, Braids and discriminants. {\em In: Braids}
  (Santa Cruz, CA, 1986), Contemp. Math. {\bf78}, Amer. Math. Soc.,
  Providence, RI, 1988, pp.~605--613.

\bibitem{OrSo88b}
{\sc P.~Orlik and L.~Solomon}, Discriminants in the invariant theory of
  reflection groups. Nagoya Math. J. {\bf109} (1988), 23--45.

\bibitem{Osi54}
{\sc M.~Osima}, On the representations of the generalized symmetric group. 
  Math. J. Okayama Univ. {\bf4} (1954), 39--56. 

\bibitem{Pop82}
{\sc V. Popov}, {\em Discrete complex reflection groups}. Communications of
  the Mathematical Institute, 15. Rijksuniversiteit Utrecht, Utrecht, 1982.

\bibitem{PT}
{\sc J.~H.~Przytycki and  P.~Traczyk}, Invariants of links of Conway
  type. Kobe J. Math. {\bf 4} (1987), 115--139.

\bibitem{Ram} 
{\sc A.~Ram}, A Frobenius formula for the characters of the Hecke algebras.
  Invent. Math. {\bf 106} (1991), 461--488.

\bibitem{RaSh}
{\sc K. Rampetas and T. Shoji}, Length functions and Demazure operators for
  $G(e,1,n)$. I, II. Indag. Math. {\bf9} (1998), 563--580, 581--594.

\bibitem{Read77}
{\sc E.~W.~Read}, On the finite imprimitive unitary reflection groups.
  J. Algebra {\bf45} (1977), 439--452.

\bibitem{Richardson82}
{\sc R.~W.~Richardson}, Conjugacy classes of involutions in Coxeter 
  groups. Bull. Austral. Math. Soc. {\bf26} (1982), 1--15.

\bibitem{Shep}
{\sc G.~C.~Shephard}, Abstract definitions for reflection groups. Canad.
  J. Math. {\bf9} (1957), 273--276.

\bibitem{ShTo}
{\sc G.~C.~Shephard and J.~A.~Todd}, Finite unitary reflection groups.
  Canad. J. Math. {\bf 6} (1954), 274--304.

\bibitem{Shi1}
{\sc J.-Y.~Shi}, Conjugacy relation on Coxeter elements. Adv. Math. {\bf 161}
  (2001), 1--19.

\bibitem{So63}
{\sc L.~Solomon}, Invariants of finite reflection groups. Nagoya Math. J.
 {\bf22} (1963), 57--64.

\bibitem{Spr}
{\sc T.~A.~Springer}, Regular elements of finite reflection groups.
 Invent. Math. {\bf25} (1974), 159--198.

\bibitem{Spr2}
{\sc T.~A.~Springer}, {\em Linear algebraic groups, second edition}.
  Progress in Math. vol.~9. Birkh\"auser, Boston--Basel--Berlin, 1998.

\bibitem{SprSt} 
{\sc T.~A.~Springer and R.~Steinberg}, Conjugacy classes. 1970
Seminar on algebraic groups and related finite groups (The Institute
for Advanced Studies, Princeton, N.J., 1968/69), pp.~167--266,
Lecture Notes in Math., vol.~131, Springer, Berlin, 1970.

\bibitem{Sta}
{\sc A.~J. Starkey}, {\em Characters of the generic Hecke algebra of a 
  system of BN-pairs}. Ph.D. thesis, University of Warwick, July 1975.

\bibitem{St51}
{\sc R.~Steinberg}, A geometric approach to the representations of the full
  linear group over a Galois field. Trans. Amer. Math. Soc. {\bf71} (1951),
  274--282; see also \cite{SteinColl}, pp.~1--9.

\bibitem{Ste}
{\sc R.~Steinberg}, Differential equations invariant under finite reflection
  groups. Trans. Amer. Math. Soc. {\bf112} (1964), 392--400; see also
  \cite{SteinColl}, pp.~173--181.

\bibitem{Steinberg0}
{\sc R.~Steinberg}, Lectures on Chevalley groups. Mimeographed notes,
  Department of Mathematics, Yale University, 1967.

\bibitem{Steinberg1}
{\sc R.~Steinberg}, Endomorphisms of linear algebraic groups. Mem. Amer. Math.
  Soc. {\bf 80} (1968), 1--108;  see also \cite{SteinColl}, pp.~229--285.

\bibitem{SteinColl}
{\sc R.~Steinberg}, {\em Collected papers, with a foreword by J-P.~Serre}.
  Amer. Math. Soc., Providence, RI, 1997.

\bibitem{Ste89}
{\sc J. Stembridge}, On the eigenvalues of representations of reflection
  groups and wreath products. Pacific J. Math. {\bf140} (1989), 353--396.

\bibitem{TeYa}
{\sc H.~Terao and T.~Yano}, The duality of the exponents of free deformations
  associated with unitary reflection groups. {\em In: Algebraic groups and
  related topics} (Kyoto/Nagoya, 1983), Adv. Stud. Pure Math. {\bf6},
  North-Holland, Amsterdam, 1985, pp.~339--348.

\bibitem{Tits}
{\sc J.~Tits}, {\em Buildings of spherical types and finite $BN$-pairs}. 
  Lecture Notes in Math. {\bf 386}, Springer, Heidelberg, 1974. 

\bibitem{TY82}
{\sc S. Tokunaga and M. Yoshida}, Complex crystallographic groups. I. 
  J. Math. Soc. Japan {\bf34} (1982), 581--593.

\bibitem{TracBanach}
{\sc P.~Traczyk}, A new proof of Markov's braid theorem. Banach Center
  Publ., Polish Acad. Sci. {\bf42} (1998), 409--419.

\bibitem{Verma1}
{\sc D.~N.~Verma}, M\"obius inversion for the Bruhat ordering on a 
  Weyl group. Ann. Sci. \'Ecole Norm. Sup. {\bf 4} (1971), 393--398.

\bibitem{Verma}
{\sc D.~N.~Verma}, The r\^ole of affine Weyl groups in the representation
  theory of algebraic Chevalley groups and their Lie algebras. {\em In:
  Lie groups and their representations}, Halsted, New York, 1975, pp.~653--705.

\bibitem{Vogel90}
{\sc P.~Vogel}, {Representation of links by braids: a new algorithm}.
  Comment. Math. Helv. {\bf65} (1990), 104--113.

\bibitem{Wag78}
{\sc A. Wagner}, Collineation groups generated by homologies of order greater
  than $2$. Geom. Dedicata {\bf7} (1978), 387--398. 

\bibitem{Wag80}
{\sc A. Wagner}, Determination of the finite primitive reflection groups over
  an arbitrary field of characteristic not $2$. I. Geom. Dedicata {\bf9}
  (1980), 239--253; II ibid. {\bf10}, 191--203; III ibid. {\bf10}, 475--523.

\bibitem{Wal78}
{\sc D. Wales}, Linear groups of degree $n$ containing an involution with two
  eigenvalues $-1$. II.  J. Algebra {\bf53} (1978), 58--67.

\bibitem{White}
{\sc N.~White}, The Coxeter matroids of Gelfand et al. {\em In:} Contemp.
  Math., vol.~197, Amer. Math. Soc. Providence, RI, 1996, pp.~401--409.

\bibitem{Wil86}
{\sc R. A. Wilson}, The geometry of the Hall-Janko group as a quaternionic
  reflection group. Geom. Dedicata {\bf20} (1986), 157--173. 

\bibitem{ZaSe80}
{\sc A. E. Zalesski\u \i\ and V. N. Sere\v zkin}, Finite linear groups generated
  by reflections.  Izv. Akad. Nauk SSSR Ser. Mat. {\bf44} (1980),
  1279--1307 (in Russian). English translation: Math. USSR-Izv. {\bf17}
  (1981), 477--503.


\end{thebibliography}
\end{document}